\documentstyle[11pt,xypic,amssymb]{article}

\newlength{\quotesize}
\newcommand{\vvv}{$\,$}
\def\appendix{\par
              \setcounter{section}{0}%
	      \def\@chapapp{Appendix}
              \def\thechapter{}
              \chapter{}
              \def\thesection{\Alph{section}}
              \def\thesubsection{\thesection\arabic{subsection}}
              \def\thesubsubsection{\thesubsection.\arabic{subsubsection}} 
}
\newtheorem{prop}{Proposition}[subsection]
\newtheorem{lemma}[prop]{Lemma}
\newtheorem{Th}[prop]{Theorem}
\newtheorem{corol}[prop]{Corollary}
\newtheorem{Loc}[prop]{Localization theorem}
\newtheorem{remark}[prop]{Remark}
\newenvironment{rem}[1]{\par\smallskip\noindent{\bf Remark. }#1}
{\par\smallskip}
\def\endpr{\hfill $\Box$}
\newenvironment{proof}[2]{\par\medskip\noindent{\bf #1: }#2}
{\endpr\par\medskip}
\renewcommand{\thesection}{\arabic{section}}
\renewcommand{\thesubsection}{\thesection.\arabic{subsection}}
\renewcommand{\theenumi}{(\roman{enumi})}

\def\afterthmseparator{.}
\def\@begintheorem#1#2{\trivlist \item[\hskip \labelsep{\bf
#1\ #2\unskip\afterthmseparator}]} %
\def\@opargbegintheorem#1#2#3{\trivlist
      \item[\hskip \labelsep{\bf
#1\ #2\ (#3)\afterthmseparator}]} 
\def\@eqnnum{{\reset@font\rm (\theprop)}}

\expandafter\ifx\csname amssym.def\endcsname\relax \else\endinput\fi
\expandafter\edef\csname amssym.def\endcsname{%
       \catcode`\noexpand\@=\the\catcode`\@\space}
\catcode`\@=11
\def\undefine#1{\let#1\undefined}
\def\newsymbol#1#2#3#4#5{\let\next@\relax
 \ifnum#2=\@ne\let\next@\msafam@\else
 \ifnum#2=\tw@\let\next@\msbfam@\fi\fi
 \mathchardef#1="#3\next@#4#5}
\def\mathhexbox@#1#2#3{\relax
 \ifmmode\mathpalette{}{\m@th\mathchar"#1#2#3}%
 \else\leavevmode\hbox{$\m@th\mathchar"#1#2#3$}\fi}
\def\hexnumber@#1{\ifcase#1 0\or 1\or 2\or 3\or 4\or 5\or 6\or 7\or 8\or
 9\or A\or B\or C\or D\or E\or F\fi}

\font\tenmsa=msam10 scaled 1200
\font\sevenmsa=msam7 scaled 1200
\font\fivemsa=msam5 scaled 1200
\newfam\msafam
\textfont\msafam=\tenmsa
\scriptfont\msafam=\sevenmsa
\scriptscriptfont\msafam=\fivemsa
\edef\msafam@{\hexnumber@\msafam}
\mathchardef\dabar@"0\msafam@39
\def\dashrightarrow{\mathrel{\dabar@\dabar@\mathchar"0\msafam@4B}}
\def\dashleftarrow{\mathrel{\mathchar"0\msafam@4C\dabar@\dabar@}}

\def\ulcorner{\delimiter"4\msafam@70\msafam@70 }
\def\urcorner{\delimiter"5\msafam@71\msafam@71 }
\def\llcorner{\delimiter"4\msafam@78\msafam@78 }
\def\lrcorner{\delimiter"5\msafam@79\msafam@79 }
\def\yen{{\mathhexbox@\msafam@55 }}
\def\checkmark{{\mathhexbox@\msafam@58 }}
\def\circledR{{\mathhexbox@\msafam@72 }}
\def\maltese{{\mathhexbox@\msafam@7A }}

\font\tenmsb=msbm10 scaled 1200
\font\sevenmsb=msbm7 scaled 1200
\font\fivemsb=msbm5  scaled 1200
\newfam\msbfam
\textfont\msbfam=\tenmsb
\scriptfont\msbfam=\sevenmsb
\scriptscriptfont\msbfam=\fivemsb
\edef\msbfam@{\hexnumber@\msbfam}
\def\Bbb#1{{\fam\msbfam\relax#1}}
\def\widehat#1{\setbox\z@\hbox{$\m@th#1$}%
 \ifdim\wd\z@>\tw@ em\mathaccent"0\msbfam@5B{#1}%
 \else\mathaccent"0362{#1}\fi}
\def\widetilde#1{\setbox\z@\hbox{$\m@th#1$}%
 \ifdim\wd\z@>\tw@ em\mathaccent"0\msbfam@5D{#1}%
 \else\mathaccent"0365{#1}\fi}
\font\teneufm=eufm10 scaled 1200
\font\seveneufm=eufm7 scaled 1200
\font\fiveeufm=eufm5 scaled 1200
\newfam\eufmfam
\textfont\eufmfam=\teneufm
\scriptfont\eufmfam=\seveneufm
\scriptscriptfont\eufmfam=\fiveeufm

\let\goth\frak 

\csname amssym.def\endcsname

\def\oplusl{\mathop\oplus\limits}
\def\no{$n^\circ$~}
\def\oO{\overline O_\mu}
\def\oC{\overline C^\lambda } 
\def\Tr{{\rm Tr}}

\def\F{{\Bbb  F}}
\def\Q{{\Bbb  Q}}
\def\Z{{\Bbb  Z}}
\def\C{{\Bbb  C}}
\def\g{{\goth g}}
\def\a{{\goth a}}
\def\t{{\goth t}}
\def\h{{\goth h}}
\def\sl{{\goth s}{\goth l}}
\def\Ad{{\rm Ad}\,}
\def\ad{{\rm ad}\,}
\def\Ext{{\rm Ext}}
\def\SL{{\rm SL}}
\def\GL{{\rm GL}}

\def\Bun{{\rm Bun}}
\def\eqto{\stackrel{\sim}{\longrightarrow}}
\def\eqge{\stackrel{\sim}{\longleftarrow}}

\def\Gr{{\rm Gr}}
\def\j{j}
\def\prim{prim}
\def\ch{ch}
\def\CP{{\Bbb C}{\Bbb P}^1}
\def\CPN{{\Bbb C}{\Bbb P}^n}
\def\Aut{\mathop{\fam0 Aut}}
\def\rad{\mathop{\fam0 rad}}
\def\wh{\widetilde h}
\def\M{{\cal M}}
\def\MI{\goth m}
\def\longlleftarrow{\mathrel{\raise.3ex\hbox{$\longleftarrow$}
\llap{\lower.3ex\hbox{$\longleftarrow$}}}}
\def\longllleftarrow{\mathrel{\hbox{$
\mathord{\longleftarrow}$}\llap{\raise.5ex\hbox{$\longleftarrow$}}
\llap{\lower.5ex\hbox{$\longleftarrow$}}}}
\def\b{{\goth b}}
\def\wh{\widetilde h}
\def\wh{\widetilde h}
\def\O{{\cal O}}
\def\proj{\mathop{\lim\limits_{\smash{\longleftarrow}}}\limits}
\def\inj{\mathop{\lim\limits_{\longrightarrow}}\limits}

\def\XT{\mbox{\bf X}_*(T)}                  
\def\scriptXT{\mbox{\bf\scriptsize X}_*(T)} 
\def\XTu{\mbox{\bf X}^*(T)}                 
\def\XTup{\mbox{\bf X}^*(T^\vee)}           

\def\RRR{{\Bbb R}}                           
\def\Ct{\C[\t]}                         

\def\R{\mbox{\rm R}}
\def\D{\mbox{\rm D}}                        
\def\RHom{\mbox{\rm RHom}}
\def\supp{\mbox{\rm supp}}
\def\Hom{\mbox{\rm Hom}}
\def\mod{\mbox{\,\rm mod}\,}
\def\rk{\mbox{\rm rk}\,}
\def\Vect{\mbox{\Bbb Vect}}
\def\Spec{\mbox{\rm Spec}}
\def\DegTrans{\mbox{\rm Deg. transc. }}
\def\Rep{\mbox{\rm Rep}}
\def\Cat{\mbox{\rm Cat}}
\def\Lie{\mbox{\rm Lie}\,}
\def\Im{\mbox{\rm Im}}
\def\Ob{\mbox{\rm Ob}}
\def\pt{\mbox{\rm pt}}

\def\gotht{\t}

\newcommand{\NN}{{{\cal N}ilp}}
\newcommand{\pp}{{\cal P}}
\newcommand{\K}{{\Bbb K}}
\newcommand{\oo}{{\cal O}}
\newcommand{\B}{{\cal B}}
\newcommand{\mm}{{\cal M}}
\newcommand{\dm}{{D^b_L(M/L)_c}}
\newcommand{\aaa}{{\cal A}}
\begin{document}

\title{\bf PERVERSE SHEAVES ON A LOOP GROUP AND LANGLANDS' DUALITY}
\author{Victor GINZBURG \\
The University of Chicago, \\
Department of Mathematics \\
Chicago, IL 60637, USA \\
E-mail: ginzburg@math.uchicago.edu}
\date{}
\maketitle
\hfill{\small alg-geom/9511007}\break

\section*{Introduction}

\renewcommand{\theequation}{0.\arabic{equation}}
\renewcommand{\theprop}{0.\arabic{prop}}

Grothendieck associated to any complex of sheaves ${\cal F}$ on a
variety $X$ over a finite field $\F$ the function
\[
x\mapsto\chi_{{}_{\cal F}}(x)=\sum(-1)^i\Tr({\sf F}{\sf r};{\cal H}_x^i{\cal F}),
\qquad x\in X^{{\sf F}{\sf r}}
\]
on the set of $\F$-rational points of $X$ given by the alternating
sum of traces of ${\sf F}{\sf r}$, the Frobenius action on stalks of the cohomology
sheaves
${\cal H}^i{\cal F}$.
He then went on to initiate an ambitious program of giving
{\em geometric\/}
(= sheaf theoretic) meaning to various classical
{\em algebraic\/}
formulas via the above ``func\-tions-faisceaux'' correspondence
${\cal F}\mapsto\chi_{{}_{\cal F}}$.
This program got a new impetus with the discovery of perverse sheaves
\cite{BBD}, for it happens for certain mysterious  reasons that most
of the functions encountered `in nature' arise via the
``functions-faisceaux'' correspondence from perverse sheaves. In
the present paper Grothendieck's philosophy is applied to what may be
called the Geometric Langlands duality.

The relevance of the intersection cohomology technique to our problem
was first pointed out by Drinfeld~\cite{D} and Lusztig~\cite{Lu 1}.
Later, in the remarkable  paper~\cite{Lu}, Lusztig established
algebraically a surprising connection between finite-dimensional
representations of a semisimple Lie algebra and the Kazhdan-Lusztig
polynomials for an affine Weyl group. It is one of the  purposes
of the present paper to provide a geometric interpretation of~\cite{Lu}.
It should be mentioned however that the results of~\cite{Lu}
are used in a very essential way in the proof of our main theorem
in section 2.

To begin with, recall the classical notion of a
{\em spherical function}.
Fix a semisimple group $G$. Given a ring $R$, let $G(R)$ denote
the corresponding group of $R$-rational points. Let $p$ be a prime,
$\Q_p$ the $p$-adic field, and $\Z_p$ its ring of integers. Let
$\C[G(\Z_p)\backslash G(\Q_p)/G(\Z_p)]$,
be  the Hecke algebra with respect to convolution of complex valued
$G(\Z_p)$-biinvariant, compactly supported functions on $G(\Q_p)$.
Let further
$T(\C)\subset G(\C)$
be a maximal torus in the corresponding complex group and
$X_*(T)=Hom_{alg}(\C^*,T(\C))$
the lattice of one-parameter subgroups of
$T(\C)$. Let  $\C[X_*(T)]$
be the group algebra of this lattice acted upon naturally by $W$,
the Weyl group of $G$.
Then, one has the following classical result due to Satake:

\noindent
{\em There is an algebra isomorphism {\rm(}explicitly constructed
by Macdonald})
\begin{equation}
\C[G(\Z_p)\backslash G(\Q_p)/G(\Z_p)]\simeq\C[X_*(T)]^W
\end{equation}
Our aim is to produce a sheaf-theoretic counterpart of isomorphism
(0.1). This is done in several steps by reinterpreting both the
left hand and the right hand sides of (0.1). First, we introduce
the dual complex torus $T^\lor $ such that
$X_*(T)=X^*(T^\lor )$
is the weight lattice of $T^\lor $. Let $G^\lor $ be the Langlands dual of
$G$, the complex semisimple group having $T^\lor $ as a maximal torus
and having the root system dual to that of $G$. Then, the RHS of
(0.1) can be rewritten as
\begin{equation}
\C[X^*(T^\lor )]^W\simeq\C[G^\lor ]^{G^\lor }
\end{equation}
where
$\C[G^\lor ]^{G^\lor }$
stands for the algebra of polynomial class
functions on $G^\lor $. Let $\Rep_{_{G^\lor}}$ be the abelian category of finite
dimensional rational complex representations of $G^\lor $, and
$K(\Rep_{_{G^\lor}})$
its Grothendieck group. Assign to a representation
$V\in\Rep_{_{G^\lor}}$
its character, a class function on $G^\lor $. This yields a natural
algebra isomorphism
$\C\otimes_\Z K(\Rep_{_{G^\lor}})\eqto\C[G^\lor ]^{G^\lor }$.
By (0.2), the isomorphism (0.1) can thus be written as
\begin{equation}
\C[G(\Z_p)\backslash G(\Q_p)/G(\Z_p)]\simeq\C
\otimes_\Z K(\Rep_{_{G^\lor}})
\end{equation}
We next turn to the LHS of (0.1). Let
$\F=\Z_p/p\cdot\Z_p$
be the residue class field and
$\overline\F$
an algebraic closure
of $\F$. First we use the well known analogy between the fields
$\Q_p$ and
$\F((z))$
(=the Laurent formal power series field). The algebra
$\C[G(\Z_p)\backslash G(\Q_p)/G(\Z_p)]$
is unaffected by the substitution
\[
\Q_p\longleftrightarrow\F((z)),\quad
\Z_p\longleftrightarrow\F[[z]],\quad
p\longleftrightarrow z
\]
Thus, the LHS of (0.3) can be rewritten as
\[
\mbox{{\em$
G(\F[[z]])$-invariant functions on $G(\F((z)))/G(\F[z]])$}}
\]
Now, following Grothendieck, we view the discrete set
$G(\F((z)))/G(\F[[z]])$
as the set of $\F$-rational points of
$\Gr=G(\overline\F((z)))/G(\overline\F[[z]])$,
an infinite dimensional algebraic variety. We introduce a category
$P(\Gr)$ whose objects are semisimple
$G(\overline\F[[z]])$-equivariant perverse sheaves on $\Gr$. There
is a natural tensor product structure on the category $P(\Gr)$ given
by a sort of convolution. Moreover, the ``function-faisceaux''
correspondence assigns to
${\cal F}\in P(\Gr)$
a function
$\chi_{{}_{\cal F}}$
on
$G(\F((z)))/G(\F[[z]])$.
That gives an isomorphism of the Grothendieck group,
$\C\otimes_\Z K(P(\Gr))$,
with the algebra of $G(\F[[z]])$-invariant functions on
$G(\F((z)))/G(\F[[z]])$.
Thus, the algebra isomorphism (0.1) takes the following final form
\begin{equation}
\C\otimes_\Z K(P(\Gr))\simeq\C
\otimes_\Z K(\Rep_{_{G^\lor}})
\end{equation}
The main result of this paper (theorem 1.4) can now be formulated
as follows:

\smallskip

{\em There is an equivalence of tensor categories
$P(\Gr)\simeq\Rep_{_{G^\lor}}$
that induces isomorphism (0.4) on the corresponding Grothendieck
rings. Moreover, the underlying vector space of a representation
$V\in\Rep_{_{G^\lor}}$
may be identified canonically with the hyper-cohomology of the
corresponding perverse sheaf
${\cal P}(V)\in P(\Gr)$.}

\smallskip

\noindent
In the main body of the paper we replace the field $\overline\F$
by the field $\C$ of complex numbers. This replacement does not
affect the formulation of the theorem above, except that $\Gr$ is
now regarded as an (infinite-dimensional) complex algebraic variety
equipped with the usual Hausdorff topology. Thus, finite and
$p$-adic fields serve only as a motivation and will never appear
in the rest of the paper.

\smallskip

A few words on the organization of the paper are in order.

In the first chapter we formulate most of the results of the paper
without going into any technical details and sometimes without even
giving proper
definitions. This is made to simplify the reading of the paper,
for working with perverse sheaves on an infinite-dimensional
variety like $\Gr$ requires a lot of extracare in technical details
to make our approach rigorous. Thus, the first chapter contains
the statement of the main theorem and its basic applications.

Chapters 2 and 3 are devoted to the proof of the main theorem.
The key ingredient of the proof is a construction of a
{\em fiber functor\/} on the tensor category $P(\Gr)$ in terms of
the equivariant cohomology.

There
are two kinds of applications  of the main theorem
that we consider. The first one is of `global' nature,
concerning a smooth compact complex curve $X$ of genus $>1$. We
set up a Langlands-type correspondence between
`{\em automorphic perverse sheaves}' on the moduli space of principal
holomorphic $G$-bundles on $X$ and homomorphisms $\pi_1(X)\to G^\lor
$, respectively.
The detailed treatment of this subject including definition of the
{\em modular stack\/} of principal $G$-bundles on a smooth complex
curve is given in chapter 6. We include here new proofs of
 two (known) basic geometric results underlying the theory.
The first one says that any algebraic $G$-bundle on an affine
curve is trivial. The second one says that the {\it global Nilpotent variety}
is a Lagrangian substack in  the cotangent bundle to the moduli space.
We then show how these results can be combined with the main theorem
to obtain the Langlands correspondence (one way).

Applications of the second type related to the topology of the
Grassmannian $\Gr$ are discussed in chapters 4--5. It is shown there
that the cohomology,
$H^*(\Gr,\C)$,
is isomorphic to the Symmetric algebra of the
centralizer of the
{\em principal nilpotent\/} in $\Lie (G^\lor) $. The principal nilpotent
itself turns out to be nothing but the first Chern class of the
Determinant bundle on $\Gr$. Among other things, that enables us
to give natural proofs of the results of Lusztig~\cite{Lu} and
Brylinski~\cite{Br} concerning the $q$-analogue of the weight
multiplicity, and also
to establish connections with various interesting
questions in represention theory.

\smallskip

This paper is a slightly revised TeX-version of the 
typewritten manuscript with
the same title distributed in 1989 (without chapter 6). 

I am greatly indebted to V.~Drinfeld without whom this work would
have never been carried out. He initiated the whole project by raising the question
whether the category $P(\Gr)$ studied below is a tensor category.
I am grateful to A.~Beilinson who explained me a lot about the geometry
of modular stacks which is crucial for chapter 6. The present paper
may be seen, in fact, as a preparatory step towards 
the Geometric Langlands correspondence \cite{BeDr}.
 
Thanks are also
due to R.~Bezrukavnikov,\,R.~Brylinski,\,B.~Kostant\, and G.~Lusztig
for some very useful comments
and suggestions.
\vskip 10mm
\newpage
\centerline{\bf Table of Contents}
\vskip 5mm
$\hspace{5mm}$ {\footnotesize \parbox[t]{115mm}{{\bf
1.{ $\;$} {\sc Main results} \newline
2.{ $\;$} {\sc Tensor category of perverse sheaves} \newline
3.{ $\;$} {\sc Application of the Equivariant cohomology} \newline
4.{ $\;$} {\sc The loop group cohomology and the principal nilpotent} \newline
5.{ $\;$} {\sc Equivariant filtration and $q$-analogue of the weight multiplicity} \newline
6.{ $\;$} {\sc Moduli spaces and Hecke operators} \newline
7.{ $\;$} {\sc  Appendix A: Tensor functors and $G$-bundles}\newline
8.{ $\;$} {\sc Appendix B: Equivariant hyper-cohomology}}}}
\vskip 10mm

\section{Main results}

\renewcommand{\theequation}{\thesubsection}
\renewcommand{\theprop}{\thesubsection}

\subsection{Duality} Given a complex torus $T$, let $T^\lor $ denote
the group of 1-dimensional local systems on $T$ with the tensor product
group structure. In ground to earth terms, $T^\lor $ is
the group of all homomorphisms:
$\pi_1(T)\to\C^*$.
This is again a complex torus which is called the dual of $T$.
More explicitly, let $\goth t$ be the Lie algebra of $T$,
${\goth t}^*$ the dual space, $X_*(T)$ the lattice of holomorphic
homomorphisms: $\C^*\to T$, and $X^*(T)$ the (dual) lattice of
holomorphic homomorphisms: $T\to\C^*$. Taking the differential
of such a homomorphism at the identity defines embeddings:
$X_*(T)\hookrightarrow\Hom_\C(\C,{\goth t})=\goth t$
and
$X^*(T)\hookrightarrow\Hom_\C({\goth t},\C)={\goth t}^*$.
Then, we have canonically:
\[
T\cong{\goth t}/X_*(T),\quad T^\lor \cong{\goth t}^*/X^*(T)
\]
so that there are
canonical identifications:
\begin{equation}
({\goth t}^\lor )^*\simeq{\goth t},\quad
X^*(T^\lor )\simeq X_*(T)
\end{equation}

Now, let $G$ be a split connected semisimple group over $\Z$, and $T$
a maximal torus in $G$. We let $G^\lor $ denote the connected complex
semisimple group having the dual torus $T^\lor $ as a maximal torus
and having the root system dual to that of $G$. The group $G^\lor $
is said to be the Langlands dual to $G$~\cite{Lan}. Our immediate
goal is to give an intrinsic new construction of $G^\lor $ which
does not appeal to root systems, maximal tori, etc.

\renewcommand{\theequation}{\thesubsection.\arabic{equation}}
\renewcommand{\theprop}{\thesubsection.\arabic{prop}}

\subsection{Infinite Grassmannian} 
Given the semisimple group $G$ and
a $\C$-algebra  $R$ write $G(R)$ for the group (over $\C$) of
$R$-rational
points of $G$. Our basic choices of the ring $R$ are
$\K=\C((z))$ and $\oo=\C[[z]]$. 
Below we will more often use
polynomial rings $\C[z,z^{-1}]$ and $\C[z]$ instead of the rings
$\K$ and $\oo$, respectively.
We will use shorthand notation $LG$ and
$L^+G$ instead of $G(\C[z,z^{-1}])$ and $G(\C[z])$ respectively.
Choosing an imbedding
$G\to GL_n$ presents
$G(\K)$ as a subgroup of $\K$-valued  $n\times n$ matrices and the
group
$LG$ as the group of maps $f: \C^* \to G(\C)\subset {\rm Mat}_n(\C)$
that can be written
as a finite Laurent polynomial:
\begin{equation}
\label{(1.2.1)}
f(z)=\sum_{i=-m}^mA_i\cdot z^i,\,\mbox{
$A_i$ are $n\times n$ matrices, and} \,f(z)\in G(\C)
\,,\,\forall z\in \C^*\label{matrix}
\end{equation}

Neither topological nor algebraic structure of the groups $G(\K)$ and
$LG$ themselves will be of any importance for us. What
{\it is} important is that the corresponding `integral' groups
$G(\oo)$ and $L^+G$ each  has a decreasing chain of normal subgroups:
\[L^0 \supset L^1\supset L^2\supset\ldots\]
In either case, the subgroup $L^m , m>0$, is defined to consist of
the following loops $f$
regular at the origin and such that
$$ f(0)=1 \;\mbox{and $f$ has vanishing derivatives at $0$ up to order
$m$}\eqno(1.2.1^m)$$
(if $m=0$ we put $L^0 =\{f\,|\, f(0)=1\}$). For each $m\geq 0$,
the quotient $G(\oo)/L^m$, resp. $L^+G/L^m$, clearly has the natural
structure of a finite-dimensional algebraic group.
Observe that, for $m>0$, the quotient $L^m/ L^{m+1}$ is 
an abelian group isomorphic to the additive group $\g := Lie\,G$.

We are interested in this paper in the coset space
$\Gr:=G(\K)/G(\oo)$. 
It is called the Infinite Grassmannian for, if $G=SL_n$,
the set $\Gr$ gets identified naturally with the variety of unimodular
$\oo$-lattices of maximal rank in $\K^{^n}$.
The `finiteness' properties of the Infinite
Grassmannian listed in the following proposition allow us to apply to $\Gr$
standard algebro-geometric constructions as if it was a
finite-dimensional
algebraic variety.

\begin{prop}
\label{exaust}
\begin{enumerate}
\item The set $\Gr$ is the union of an infinite
 sequence of $G(\oo)$-stable subsets:
$\Gr_1\subset\Gr_2\subset\ldots$.
\item These subsets have a compatible structure of finite-dimensional
projective varieties of increasing dimension, i.e. $Gr_i \hookrightarrow Gr_j$
is a projective imbedding for any $i<j$.
\item For any $i$ the $G(\oo)$-action on $\Gr_i$ factors through an
algebraic action of the group $G(\oo)/L^m$, where $m=m(i)$ is large
enough.
\item Each piece $\Gr_i$ consists of finitely many $G(\oo)$-orbits.
\end{enumerate}
\end{prop}

{\bf Comments on Proof:} Part (i) of the proposition is due to 
Lusztig~\cite[n.~11]{Lu}). Let $Grass$ be the
set of all $\O$-submodules in the vector space
$\g(\K)$, viewed as an $\O$-module.
 Lusztig considered a
map $G(\K) \to Grass$ given by the assignment
$f \mapsto (Ad\,f)\,\g(\oo)$. He showed that this map yields
a bijection of the coset space
$G(\K)/G(\oo)$ with a set $\B$ of Lie subalgebras
in $\g(\K)$ subject to certain algebraic conditions. Further,
for each $i \geq 1$, following Lusztig
put 
\begin{equation}
\label{(1.2.3)}
\Gr_i := \{L \in \B\,|\, z^i\cdot\g(\oo) \subset L \subset
z^{-i}\cdot\g(\oo)\}
\end{equation}
With this definition, parts (i) and (ii) of the proposition are
immediate.
Part (iv) follows by comparing the known (cf.~\cite[ch.~8]{PS} and
section 1.4 below)
parametrization of $G(\oo)$-orbits with the above definition of the sets
$\Gr_i$. It remains to prove (iii). This can be done for each
$G(\oo)$-orbit
separately, due to (iv). Furthermore, it suffices to check the claim
for a single point of each $G(\oo)$-orbit, for the groups $L^m$ are
normal in $G(\oo)$. But every orbit contains
a point represented by a group homomorphism $\lambda:\C^* \to G$, see n. 1.4.
For such a point claim (iii) is clear. Indeed, the image of $\lambda$
is contained in a maximal torus $T\subset G$. Then, $\lambda$ may be
viewed, cf. $n^\circ 1.1$, as an element of the lattice
$X_*(T)\subset\t=Lie T$. Then, in property
(iii), one may take $m$ to be an
integer, such that $m > \alpha(\lambda)$, for any
root $\alpha \in \t^*$.  $\quad\square$

One may introduce a similar Grassmannian in the polynomial setup, putting
$\Gr' = LG/L^+G$. It turns out that the object one gets
in this way is not only {\it similar} but is in fact {\it identical}
to the previous one. Specifically, the natural imbeddings
\[LG \hookrightarrow G(\K) \quad,\quad L^+G \hookrightarrow G(\oo)\]
give a map $j: \Gr' \to \Gr$. This map is injective, for it is clear
that
$LG\,\cap\, G(\oo) = L^+G$. Moreover,
the following comparison result shows that  $j(\Gr')=\Gr$ and that
we may (and will) make
no distinction between $\Gr$ and $\Gr'$.

\begin{prop}
\label{1.2.4.}
 The map $j$
is an $LG$-equivariant isomorphism $j: \Gr'
\stackrel{\sim}{\to} \Gr$. Any $G(\oo)$-orbit in $\Gr$ is the
image of a  single $L^+G$-orbit in $Gr'$.
\end{prop}

{\bf Proof:} Define an exhaustion
$\Gr'_1\subset\Gr'_2\subset\ldots = \Gr'$ similar to the one
introduced in the
proposition \ref{exaust}. Equation (1.2.3) shows that $\Gr_i$ may be viewed as a
variety of vector subspaces in the
finite dimensional $\C$-vector space 
$z^{-i}\cdot\g(\oo)/z^i\cdot\g(\oo)$,
subject to certain algebraic conditions. Observe now that replacing
power series by polynomials affects neither the space
\[z^{-i}\cdot\g(\oo)/z^i\cdot\g(\oo)
=z^{-i}\cdot\g(\C[z])/z^i\cdot\g(\C[z])\]
nor the algebraic conditions (the `infinite tail' of all the powers of
$z$ greater than $i$ disappears in the quotient). It follows that $j(\Gr'_i)=\Gr_i$ for
any $i$. Moreover, it follows that $G(\oo)$-orbits in $\Gr_i$ are the
same
things as $L^+G$-orbits in $\Gr'_i$. That completes the proof.
$\square$

There is another realization of the Infinite Grassmannian
often used in topology.
Let $K$ be a maximal compact subgroup of $G$ and 
$S^1$ the unit circle. Write $\Omega$ for the group
of based polynomial loops
$f:S^1\to K$, i.e. expressions like (1.2.1) with $f(z)\in K$ for all
$z\in S^1$, and $f(1)=1$.
Such a loop extends uniquely to a polynomial map $f: \C^* \to G$ such
that
$f({\overline z})={\overline {f(z)}}$ where `bar' stands for the complex
conjugation (in $\C$) on the LHS and for the involution on $G$
corresponding
to the real form $K$ on the RHS. This way one gets a group imbedding
$\Omega \hookrightarrow  LG$.
One has an ``Iwasawa'' decomposition, see \cite{PS}:
$LG=\Omega\cdot L^+G$, $\Omega\cap L^+G=\{1\}$.
Thus, the Grassmannian
$\Gr=LG/L^+G$
can be identified with the group $\Omega$. We endow both $\Gr$ and
$\Omega$  with  the  topology  of  direct   limit   of   closed
finite-dimensional
subvarieties
$\Gr_1\subset\Gr_2\subset\ldots$
and let
$m:\Gr\times\Gr\to\Gr$
denote the map corresponding to the multiplication in $\Omega$.
We'll show later that, for any $i,j\ge1$, there exists an integer
$k=k(i,j)\gg0$
such that the product of $\Gr_i$ and $\Gr_j$ is contained in
$\Gr_k$ and, moreover, the multiplication map
$m:\Gr_i\times\Gr_j\to\Gr_k$
is a morphism of real-analytic sets. This makes $\Gr\simeq\Omega$
 a topological group (with an additional sub-analytic
structure).

Let $D^b(\Gr)$ be the bounded derived category of constructible complexes
on $\Gr$ whose support is contained in a big enough subset $\Gr_i$,
i.e., the direct limit of the categories $D^b(\Gr_i)$. Given
$L,M\in D^b(\Gr)$,
define
$L* M\in D^b(\Gr)$
by the formula:
$L* M=m_*(L\boxtimes M)$.
The assignment
$L,M\mapsto L* M$
is called
{\em convolution}.
The convolution is associative, i.e., for any
$M_1,M_2,M_3\in D^b(\Gr)$,
there is a functorial isomorphism:
$M_1*(M_2* M_3)=(M_1* M_2)* M_3$.

\subsection{Tensor category $P(Gr)$}
There is an obvious $L^+G$-action on $\Gr$ on the left. Each
$L^+G$-orbit is contained in a certain $\Gr_i$, and it is a smooth
locally closed algebraic subvariety of $\Gr_i$
(cf.~\cite[ch.~8]{PS}).
Each orbit is known \cite{PS} to be isomorphic to a vector bundle
over a partial flag variety (for the group $G$). In particular, the
orbits are simply-connected. By proposition 1.2.2, all the $\Gr_i$ are
$L^+G$-stable subsets of $\Gr$, and the orbits form a finite stratification
of each of the $\Gr_i$.

Let
$ IC({O})$
denote the Intersection cohomology complex of the closure of a
$L^+G$-orbit $O$, extended by $0$ to all $\Gr$. As explained,
$ IC({O})$
is a well-defined object of $D^b(\Gr)$. Let $P(Gr)$ be the abelian
full subcategory of $D^b(\Gr)$ whose objects are perverse sheaves
on $\Gr$ isomorphic to finite  direct sums of complexes
$ IC({O})$
(for various orbits $O$). Any object $L\in P(Gr)$ has
finite-dimensional hyper-cohomology $H^*(L)$.

Let $\Vect$ denote the abelian tensor category of finite-dimensional
vector spaces over $\C$ with the standard tensor product.
The following result will be proved in the next chapter.

\begin{Th}{\it
\begin{enumerate}
\item If $L,M\in P(Gr)$ then $L* M\in P(Gr)$;
\item The pair $(P(Gr),*)$ is a semisimple rigid tensor category
(cf.~\cite[def.~1.7]{DM});
\item The assignment:
$M\leadsto H^\bullet(M)$
yields an exact
fully-faithful tensor functor:
$P(Gr)\to\Vect$.
\end{enumerate}
}\end{Th}

\subsection{Equivalence of categories}
We remind the reader one of the main results of~\cite{DM}.
It says~\cite[thm.\ 2.1]{DM} that any tensor category having the
properties as in theorem~1.3 is equivalent to the category
$\Rep(G^*)$ of finite-dimensional rational representations of a
reductive group $G^*$. This applies, in particular, to the category
$P(Gr)$. Thus, starting from a semisimple group $G$ we have
constructed another reductive group $G^*$. The group $G^*$ turns out to
be isomorphic to $G^\lor $, the dual group of $G$. Thus, we get a new
construction of the dual group.

We turn to the formulation of the main result of the paper.
Let $T\subset G$ be a maximal torus.
Any homomorphism $\lambda:\C^*\to T$, viewed as a map:
$S^1\to G$, determines the coset
$\lambda\cdot L^+G\subset LG$
and, hence, a point in $\Gr$. Let $O_\lambda$ be the
$L^+G$-orbit of that point. Any $L^+G$-orbit in $\Gr$ is known to be
an orbit of the form $O_\lambda$ for some
$\lambda\in X_*(T)$,
and $O_\lambda=O_\mu$ iff $\lambda$ and $\mu$ are conjugate by
$W$, the Weyl group of the pair $(G,T)$. For
$\lambda\in X_*(T)$
we set
$ IC_\lambda:= IC(O_\lambda)$,
but it should be understood that $\lambda$ is only determined here
up to the action of $W$.


It is well known that irreducible finite-dimensional rational
representations
of the group $G^\lor $ are labeled by their extreme weights. Let $V_\lambda$
denote the irreducible representation with extreme weight $\lambda$.
Clearly, it depends only on the $W$-orbit of $\lambda$. Thus, we see
from (1.1) that simple objects of the categories $P(Gr)$ and
$\Rep_{_{G^\lor}}$ are both indexed by the same set:
$X_*(T)/W=X^*(T^\lor )/W$
(note that pairs $(G,T)$ and $(G^\lor ,T^\lor )$ have identical Weyl groups).

\begin{Th}{\it
There is an equivalence of the tensor categories $P(Gr)$ and $\Rep_{_{G^\lor}}$
which sends $ IC_\lambda$ to $V_\lambda$. Furthermore, the hyper-cohomology
functor $H^\bullet$ on $P(Gr)$ goes, under the equivalence, to the forgetful
functor:
$\Rep_{_{G^\lor}}\to\Vect$.
}\end{Th}
\subsection{A conjecture on automorphic sheaves}

The purpose of this section is to explain at a heuristic level
the main idea of the connection between theorem 1.4.1 above and
Langlands'
paramet-
rization of so-called `automorphic sheaves' on the moduli space of
principal  $G$-bundles. Rigorous approach to the moduli spaces in
question requires some special techniques, which will be worked out 
later in chapter 6. It seems instructive at this point, however, 
to outline the perspective without waiting till chapter 6.

Let $X$ be a smooth compact complex (algebraic) curve.
Fix a point $x\in X$, and let $\O_x$ be the local ring at $x$, $\MI_x\subset
\O_x$ the maximal ideal of functions vanishing at $x$, and $\K_{\,x}$ the
fields of fractions of $\O_x$, that is the field of germs of rational
functions on a punctured neighborhood of $x$.

Let $G$ be a connected  semisimple algebraic group, and
$G(\K_{\,x})$ and $G(\O_x)$ the corresponding
groups of rational points, viewed as infinite dimensional complex
groups. Set $Gr_{_X} := G(\K_{\,x})/G(\O_x)$, and write $P(Gr_{_X})$ for the
category of $G(\O_x)$-equivariant perverse sheaves on
$Gr_{_X}$ with compact support. Any such perverse sheaf turns out to be
automatically semisimple (for there are no extensions between
simple objects of $P(Gr_{_X})$, due to the fact that all
$G(\O_x)$-orbits on $Gr_{_X}$ are even dimensional).
 Thus, the setup is very similar to that of
section 1.2. Moreover, we claim that it is not only similar,
but in fact identical:

\begin{lemma}
\label{grx}
There is a natural equivalence of tensor categories
\[P(Gr_{_X}) \simeq P(Gr).\]
\end{lemma}

To prove the lemma, one first replaces the rings $\K_{\,x}$ and $\O_x$
by the corresponding $\MI_x$-adic completions. Repeating the argument of
proposition 1.2.4, one shows that this doesn't change the grassmannian
$Gr_{_X}$. It remains to note, that the curve $X$ being smooth,
choosing a local parameter $z$ at the point $x\in X$ yields
an isomorphism of rings $\hat  \O_x \simeq \C[[z]]$.
Whence, an isomorphism of the grassmannians 
$Gr_{_X} \simeq Gr$, and the induced equivalence
$P(Gr_{_X}) \simeq P(Gr)$. Observe finally, that 
this equivalence is in effect independent of the
choice of a local coordinate. Indeed, choosing
another local coordinate, $z'$, at the point $x\in X$
amounts to choosing an isomorphism  $\C[[z']] \simeq \C[[z]]$;
the latter gives an automorphism of $Gr$ that preserves
$L^+G$-orbits.$\quad\square$.

Further, let $\O_{\text{out}}$
denote the ring of regular functions on the punctured curve $X\backslash
\{x\}$. Restricting functions to a neighborhood of the puncture
yields an algebra imbedding $\O_{\text{out}}\hookrightarrow
\K_{x}$. Whence, the group imbedding 
$G(\O_{\text{out}})\hookrightarrow
G(\K_{x})$.
In chapter 6 we will endow the coset space
$\M=G(\K_{\,x})/G(\O_{\text{out}})$ with a certain structure of
an (infinite dimensional) algebraic variety. This variety
has a natural  left $G(\O_x)$-action. Moreover, we will show
that
the orbit space $G(\O_x)\backslash\M$  has a reasonable
algebraic structure (though not of a `variety' but rather
of a `stack', see e.g. \cite{Sorger}).
The well-known `double-coset construction' provides the following natural
isomorphism of stacks,
see \cite{Sorger} or theorem 6.3.1 and proposition 6.3.8 of chapter~6:
\begin{equation}
G(\O_x)\backslash\M \,\simeq\, G(\O_x)\backslash
G(\K_{\,x})/G(\O_{\text{out}})\, \simeq\, Bun_G 
\label{dcoset}
\end{equation}
where $\Bun_G$ is
the modular stack of the (isomorphism classes of)
algebraic principal $G$-bundles on $X$.

To 
any perverse sheaf $A\in P(Gr_{_X})$ we associate a 
functor $D^b(\Bun_G)\to D^b(\Bun_G)$, called
the {\it local Hecke functor} at the point $x$.
To that end, we identify (or rather {\it define}, cf. $\S 6.1.4$)
the derived category, $D^b(G(\cal O_x)\backslash\M)$, with
$D^b_{G(\O_x)}(\M)$, the $G(\O_x)$-equivariant derived category
on $\M$, see \cite{BL} or chapter 8 below for the definition
of an equivariant derived category.
Mimicing the construction of sheaf-theoretic convolution
given in section 3.2, one introduces a convolution pairing
\[ \ast : P(Gr_{_X}) \times D^b_{G(\O_x)}(\M)\to
D^b_{G(\O_x)}(\M).\]
We now use double-coset isomorphism (1.5.2) and lemma
\ref{grx}, to reinterpret this convolution as a bi-functor
\begin{equation}\label{ast1}
\ast : P(Gr) \times D^b(\Bun_G)\to
D^b(\Bun_G)
\end{equation}
Given $M\in  P(Gr)$, we define the  corresponding local (at 
$x$)
Hecke functor
$D^b(\Bun_G)\to D^b(\Bun_G)$ by the
formula $A \mapsto M \ast A$.

A complex
$A\in D^b(\Bun_G)$
is called a
{\it local Hecke eigen-sheaf}
if, for each $M\in P(Gr)$,
there exists a finite dimensional vector space $L_M$ and an isomorphism
\begin{equation}
M \ast A =L_M\otimes A.
\label{eigen}
\end{equation}

\begin{rem}
If $G=\SL_n$, then the orbit $O_\lambda \subset \Gr$
associated to a fundamental weight $\lambda$ is closed in $Gr$.
Hence the corresponding simple object $IC_\lambda \in P(Gr)$
reduces (up to shift) to the constant sheaf (this
is false for $G$ simple, $G\ne\SL_n$). Furthermore, we will see later that if
equation (1.5.4) holds for $M=IC_\lambda$, for
each fundamental weight $\lambda$, then
it holds for all $M \in P(Gr)$. Thus, if $G=\SL_n$, we have
only to check $(n-1)$ equations out of the infinite family of
 equations (1.5.4) and
these involve no intersection cohomology whatsoever
(in agreement with the
classical theory).$\quad\square$
\end{rem}

The constructions that we have made so far depended on the choice of
an arbitrary point $x\in X$. One can modify the construction slightly
by letting the point $x$ to vary within the curve $X$. In the special
case of $G=SL_n$ and of the complexes  $IC_\lambda$ where $\lambda$ is a
fundamental weight, this was done by Laumon \cite{La1} (cf. remark
above).
In the general case, one consideres the stack of quadruples:
\[\rm Heck_X=\{(P_1, P_2, x, \phi)\,|\,P_1, P_2\in \Bun_G, 
x\in X\,,\,\phi : P_1
|_{X\setminus\{x\}}
\stackrel{\sim}{\to} P_2 |_{X\setminus\{x\}} \}\]
where $P_1$ and $P_2$ are $G$-bundles on $X$ and $\phi$ is an
isomorphism
of their restrictions to $X\setminus\{x\}$. The stack $\rm Heck_X$ is
not algebraic.
It has, however, a natural stratification by the  algebraic substacks
$\rm Heck_{X,\lambda}\,,\, \lambda\in X_*(T)/W$, 
labeled by $L^+G$-orbits
in the grassmannian $Gr$. For each stratum, the
projection $pr_2 :\rm Heck_{X,\lambda} \to \Bun_G\times X\;,\;
(P_1, P_2, x, \phi) \mapsto (P_2, x)$ is a smooth morphism with
fiber isomorphic to $O_\lambda$, the correponding $L^+G$-orbit
in $Gr$. Thus, there is a natural functor $M \mapsto M_X$ assigning
to an object $M\in P(Gr)$ the corresponding perverse sheaf $M_X$
on  $\rm Heck_X$.

Write $pr_1 : \rm Heck_X \to \Bun_G$ for the projection 
$(P_1, P_2, x, \phi) \mapsto P_1$. We
define a convolution type functor 
\begin{equation}
\star : D^b(\rm Heck_X) \times D^b(\Bun_G)\to
D^b(X\times \Bun_G)\enspace,\enspace
N, A \mapsto N\star A = (pr_2)_*(N \otimes pr_1^!A)
\label{ast}
\end{equation}
Setting in particular,
$N=M_X$, we thus associate with any $M\in  P(Gr)$, the  corresponding 
{\it global}
Hecke functor given by the formula
\begin{equation}
\label{(1.5.7)}
T_M: D^b(\Bun_G)\to D^b(X\times \Bun_G)\quad,\quad A\mapsto
T_M(A)=M_X\star A
\end{equation}
We will now show, using Theorem 1.4.1, that to some Hecke eigen-sheaves
on $Bun_G$
one can associate naturally a homomorphism from, $\pi_1(X)$, the fundamental
group of the curve, into the Langlands dual group $G^\lor$. This way one gets
the geometric Langlands correspondence.

Given a principal $G^\lor$-bundle on $X$ with flat connection, and
a representation
$V\in\Rep_{_{G^\lor}}$,
write $V_{_P}$ for the associated vector bundle, i.e. a local system
on $X$.
On the other hand, let
$\pp : Rep_{_{G^\lor}} \stackrel{\sim}{\to} P(Gr)$ be the
equivalence inverse to the equivalence of theorem 1.4.1.
 We write $\pp(V)$ for the
perverse sheaf on $Gr$ arising, via the equivalence, from
$V\in\Rep_{_{G^\lor}}$. Associated to $\pp(V)$, is the
corresponding  {\it global} Hecke functor 
$T_{\pp(V)}: D^b(\Bun_G)\to D^b(X\times \Bun_G)$.

Here is our main result on the geometric Langlands correspondence.

\begin{Th}{\it
\label{main}
Let $A\in D^b(\Bun_G)$ be 
a local Hecke
eigen-sheaf at a certain point $x\in X$. Assume the following two
properties hold:
\begin{enumerate}
\item $A$ is an irreducible perverse  sheaf of geometric origin (see
\cite{BBD});
\item The characteristic variety of $A$ (see e.g. \cite{KS}
or\cite{Gi1}) is contained in the global nilpotent
cone $\NN \,\subset\, T^*\Bun_G$, defined by Laumon~\cite{La2}.
\end{enumerate}
 Then, $A$ is a global Hecke
eigen-sheaf. Moreover, there exists a unique (up to isomorphism)
flat $G^\lor$-bundle $P=P(A)$ on $X$ such  that,
 for any $V\in\Rep_{_{G^\lor}}$, we have
\begin{equation}
\label{act}
T_{\pp(V)}(A)= V_{_P}\boxtimes A\,.
\end{equation}
}\end{Th}

Proof of the theorem will be outlined in $n^\circ 6.6$.

The correspondence $A\mapsto P(A)$ established in the theorem
may be thought of as the (one way)
geometric Langlands correspondence. To explain this, fix a point
$x\in X$. Recall that taking the monodromy representation in the fiber
at $x$ of a flat $G^\lor$-bundle $P$ on $X$ sets-up a bijection
between
the isomorphism classes of flat $G^\lor$-bundles  on $X$,
and the conjugacy classes of group homomorphisms
$\phi:\pi_1(X,x)\to G^\lor$. Writing $\phi(A)$ for the
homomorphism associated to the flat $G^\lor$-bundle $P(A)$,
via the theorem, we get the geometric Langlands correspondence:
{\it Hecke eigen-sheaves satisfying (i)-(ii) of 1.5.7}$\mapsto$
{\it Homomorphisms $\pi_1(X,x)\to G^\lor$ up to conjugacy}.

 Constructing a correspondence
in the opposite direction turns out to be a much more complicated
task. Let $\Lambda$ be the set of isomorphism
classes of pairs $(P,\chi)$, where $P$ is a  flat $G^\lor$-bundle on
$X$
 such that,
$Aut(P)$, the group of automorphisms of $P$, is finite, and
$\chi$ is an irreducible representation of $Aut(P)$.
The 
conjecture below  is an extension to general semisimple
groups of the geometric Langlands reciprocity conjecture due to
Laumon \cite[conjecture~2.1.1]{La1} for $G=\GL_n$.
A proof of a $\cal D$-module analogue of this conjecture 
(at least of an essential part of (i) below) is going to appear in the series
of papers of Beilinson and Drinfeld [BeDr].

\newtheorem{env}{Conjecture}
\label{conjecture}
\renewcommand{\theenv}{}
\begin{env}{\it
(i)To any $(P,\chi)\in\Lambda$ one can associate a finite
set, called an`$L$-packet',
consisting of $\dim \chi$
perverse shaves $A$ on $Bun_G$ whose characteristic variety is
contained
in $\NN$, and such that equation \ref{act} holds.

(ii) If the flat bundle $P$ carries a variation of mixed Hodge
structure
in the sense of Deligne, then the corresponding perverse sheaves
on $Bun_G$ have an additional structure of {\it mixed Hodge modules}
in the sense of Saito \cite{Sa}.
}\end{env}

\subsection{``Topological'' gradation}
By the last sentence of theorem 1.4.1, there is a canonical
identification of the underlying vector space of an irreducible
representation
$V_\lambda\in\Rep_{_{G^\lor}}$
with the hyper-cohomology
$H^\bullet(IC_\lambda)$.
This  raises some new questions:

{\bf Question 1.} What is the representation-theoretic meaning
(in terms of $V_\lambda$) of the natural gradation  on
$H^\bullet(IC_\lambda)$ by cohomology degree?

{\bf Question 2.} What is the representation-theoretic meaning of
the natural $H^\bullet(\Gr,\C)$-action on $H^\bullet(IC_\lambda)$?

To answer these questions we fix a principal nilpotent $n$ in
${\goth g}^\lor $, the Lie algebra of the group $G^\lor $. By the
Jacobson-Morozov theorem, we can (and will) choose a Lie algebra
homomorphism:
\begin{equation}
\j:\sl_2(\C)\hookrightarrow{\goth g}^\lor\quad
\mbox{such that}\quad
\j\left(\begin{array}{ccc}0&1\\0&0\end{array}\right)=n
\end{equation}
We set:
\begin{equation}
h=
\j\left(\begin{array}{ccc}1&0\\0&-1\end{array}\right)
\end{equation}
Clearly, $h$ is a semisimple regular element of ${\goth g}^\lor $ that
has integral eigenvalues in any finite-dimensional ${\goth g}^\lor $-module.
(Moreover, it can be assumed without loss of generality, that $h$ is
the sum of positive coroots in a Cartan subalgebra of ${\goth g}^\lor $
and the nilpotent 
$n$ is the sum of root vectors corresponding to all simple roots.)

The answer to Question 1 above is given by the following

\begin{Th}{\it
The natural gradation on
$H^\bullet(IC_\lambda)$
corresponds, via theorem {\rm 1.4.1}, to the gradation on $V_\lambda$
by the eigenvalues of $h$, i.e. to the gradation:
\begin{equation}
V=\mathop\oplus\limits_{i\in\Z}V^h(i)\quad
\mbox{
where}\quad  V^h(i)=\{v\in V\,|\,h\cdot v=i\cdot v\}
\end{equation}
}\end{Th}

\setcounter{prop}{3}
\begin{remark}
Theorem 1.6.3 looks ambiguous as stated, for gradation (1.6.4)
clearly depends on the choice of $h$, that cannot be made canonical.
 To resolve this ``contradiction'', observe that the correspondence
between the two gradations in Theorem 1.6.3 depends on the equivalence
of Theorem 1.4.1\vvv. Such an equivalence is not unique and can only
be fixed up to automorphism of the category $\Rep_{_{G^\lor}}$. Any group
automorphism of $G^\lor $ gives rise to an automorphism of $\Rep_{_{G^\lor}}$.
Now, the element $h$ in (1.6.2) is uniquely determined up to
conjugacy. The effect of conjugating $h$ by $y\in G^\lor $ amounts to
applying the automorphism of $\Rep_{_{G^\lor}}$ induced by the conjugation
by $y$. That resolves the `contradiction'.
Similar meaning is assigned to the term `correspond' in
many other results stated in this paper, the objects on the `topological'
side are usually defined in a canonical way while their representation
theoretic counterparts are only determined up to {\em simultaneous\/}
conjugation. It should be emphasized however that if more than one
object are involved, then the `relative position' of the objects
may still make an `absolute' sense.
\end{remark}

\subsection{The $H^\bullet(\Gr,\C)$-action}
To answer Question 2, we have to describe the cohomology of the
Grassmannian fist. Let
$\Gr_1\hookrightarrow\Gr_2\hookrightarrow\ldots$
be the exhaustion of $\Gr$ mentioned in n.~1.2, and
$H^k(\Gr_1)\leftarrow H^k(\Gr_2)\leftarrow\ldots$
the corresponding projective system of cohomology. For each $k\ge0$
this system stabilizes, for all $\Gr_i$ are known to have compatible
Bruhat
cell decompositions so that $\Gr_{i+1}$ is obtained from $\Gr_i$
by attaching cells of dimensions $>k$ provided $i=i(k)$ is big enough.
We set by definition
$H^k(\Gr)=\mathop{\lim}\limits_{\stackrel{\leftarrow}{i}} H^k(\Gr_i)$.

To compute $H^\bullet(\Gr)$ we identify the Grassmannian with the Loop
group $\Omega$ viewed as a topological group (n.~1.2). Hence,
$H^*(\Omega)$ is a graded-commuta\-ti\-ve and cocommutative Hopf algebra.
By a well-known theorem $H^*(\Omega)$ is freely generated by the subspace
$\prim$ of its primitive elements. These elements can be obtained
\cite{B} by transgressing primitive cohomology classes of
the compact group $K$, i.e. by pulling them back to
$S^1\times\Omega$
via the evaluation map
\begin{equation}
{\rm ev}:S^1\times\Omega\to K,
\end{equation}
and then integrating over $S^1$. We will give
another description
of the primitive subspace
$\prim\subset H^\bullet(\Omega)$ after the following remark.

\begin{rem}
\label{can.bundle}
Let $O$ be the {\em regular \/} coadjoint orbit in
$({\goth g}^\lor )^*$, i.e., an orbit of maximal dimension in
the dual of the Lie algebra ${\goth g}^\lor $. Let $G^\lor (x)$ denote
the isotropy group of $x\in O$ and ${\goth g}^\lor (x):=\Lie G^\lor (x)$.
Then ${\goth g}^\lor (x)$ is an abelian Lie subalgebra of ${\goth g}^\lor $
whose dimension equals $\rk G$. (The easiest way to prove these facts
is to identify ${\goth g}^\lor $ with $({\goth g}^\lor )^*$ via an invariant
bilinear form on ${\goth g}^\lor $, and to view $x$ as a regular element
in ${\goth g}^\lor $ so that the algebra ${\goth g}^\lor (x)$ becomes the
centralizer of $x$ in ${\goth g}^\lor $. This centralizer can be clearly
viewed as the limit of a sequence of Cartan subalgebras of ${\goth g}^\lor $.)
Let $x'$ be another point of the orbit $O$ and let $u\in G^\lor $ be
any element such that $\Ad u(x)=x'$. Then the operator $\Ad u$ gives
an isomorphism ${\goth g}^\lor(x)\stackrel{sim}{\to}
{\goth g}^\lor(x')$.
Moreover, this isomorphism does not depend on the choice of $u$,
for $u$ is determined up to an element of the group $G^\lor (x)$ and the
latter, being commutative, acts trivially on ${\goth g}^\lor(x)$. Thus, all the algebras
${\goth g}^\lor (x)$, $x\in O$, can be canonically identified with each
other. 

Further, it is well-known, see \cite{Ko2}, that there is a canonical
 bijection between regular and semisimple
coadjoint orbits in $({\goth g}^\lor )^*$, respectively.
The latter are parametrized naturally by the orbit space
$\t/W$. Thus, associated with each point $t\in \t/W$, is
a regular orbit $O$ as above, hence, a {\it canonically defined}
 abelian Lie algebra
${\goth g}^\lor (x)$, $x\in O$.
We let $\goth a_t$ denote this ``universal'' abelian Lie
algebra, ${\goth g}^\lor (x)$, associated to $t\in \t/W$.
The family $\{\goth a_t\,,\,t\in \t/W\}$ 
may be thought  of  as
the family of fibers of a vector bundle on $\t/W$.
$\quad\square$
\end{rem}
\vskip 1mm

For $t=0$, i.e. in the case of a regular nilpotent $x$, we write
$\goth a$ instead of $\goth a_t$, for short. 
In particular, take $x=n$, cf.~(1.6.4), and endow the vector
space $\goth a$ with a grading induced
by the grading on ${\goth g}^\lor (x)$ by the eigenvalues of the adjoint
$h$-action on ${\goth g}^\lor (n)$.

\begin{prop}
\label{1.7.2}
There is a canonical graded space isomorphism{\rm:}
${\prim=~\goth a}$.
\end{prop}

We sketch a construction of the isomorphism: $\prim\cong\goth a$.
Let
$\C[{\goth g}]^G$
and
$\C[({\goth g}^\lor )^*]^{G^\lor }$
denote the graded algebras of invariant polynomials on $\goth g$
and $({\goth g}^\lor )^*$ respectively. These algebras are canonically
isomorphic via the following chain of isomorphisms:
\begin{equation}
\label{dualiso}
\C[{\goth g}]^G\eqto\C[\h]^W\simeq\C[(\h^\lor )^*]^W\eqge
\C[({\goth g}^\lor )^*]^{G^\lor }
\end{equation}
where the first and the third isomorphisms are due to Chevalley's
restriction theorem. We shall now construct the following diagram
\begin{equation}
\vcenter{\hbox{\diagram
\C[\g]^G_+\rrdouble^{(1.7.3)}\ddto_c&&\C[(\g^\lor)^*]^{G^\lor}_+\ddto^d\\
\\
prim\rrdashed^\sim|>\tip&&\a\\
\enddiagram}}
\label{seconds}
\end{equation}
where $\C[{\goth g}]^G_+$ stands for the augmentation ideal in
$\C[{\goth g}]^G$, etc.

To define the map $c$ decompose the $2$-dimensional sphere as the
union of two disks $D_+$ and $D_-$ with the common boundary $S^1$.
Take trivial principal $K$-bundles on
$D_+\times\Omega$
and
$D_-\times\Omega$
and identify them over the boundary
$S^1\times\Omega$
by means of the map (1.7.1). This way we obtain a principal
$K$-bundle on
$S^2\times\Omega$.
Given an invariant polynomial
$P\in\C[{\goth g}]^G$,
let
$\ch(P)\in H^\bullet(S^2\times\Omega)$
be the corresponding characteristic class of our $K$-bundle, and let
$c(P)$ be the integral of $\ch(P)$ along the factor $S^2$. One can
show that $c(P)$ is a primitive class in $H^*(\Omega)$. The assignment
$P\mapsto c(P)$
defines the map $c$ in (1.7.4). It is clear from the construction
that:
$\deg c(P)=\deg\ch(P)-2=2\cdot\deg P-2$.

To define the map $d$ in (1.7.4), notice first that the differential
of any polynomial $P^\lor $ on $({\goth g}^\lor )^*$ may be viewed as a
${\goth g}^\lor $-valued function on
$({\goth g}^\lor )^*$.
If $P^\lor $ is an invariant polynomial then,
for any
$x\in({\goth g}^\lor )^*$, we have:
$dP^\lor (x)\in{\goth g}^\lor (x)$.
Now, take $x$ to be a regular nilpotent and identity
${\goth g}^\lor (x)$
with $\goth a$. We define the map $d$ in (1.7.4) by the formula
$\C[({\goth g}^\lor )^*]^{G^\lor }\ni P^\lor \mapsto dP^\lor (x)\in\goth a$.
\vskip 1mm

Let
$I\subset\C[\g]^G$
and
$I^\lor \subset\C[(\g^\lor )^*]^{G^\lor }$
denote the augmentation ideals. Proposition 1.7.2 clearly follows
from diagram (1.7.4) and the following result which will be proved
in chapter~4.

\begin{lemma}
The maps $c$ and $d$ in (1.7.4) both vanish on the squares of
the augmentation ideals. The resulting maps{\rm:}
\[
I/I\cdot I\stackrel{c}{\longrightarrow}\prim\quad \mbox{
and}\quad
I^\lor /I^\lor \cdot I^\lor \stackrel{d}{\longrightarrow}\goth a
\]
are isomorphisms.
\end{lemma}

In the simply laced case a result similar to Proposition~1.7.2
has been also obtained by D.~Peterson (unpublished).

Let $U(\a)$ denote the enveloping algebra of the (commutative)
Lie algebra $\a$, and let $\check u$ denote the element of
$U(\a)$ corresponding to a cohomology class
$u\in H^\bullet(\Omega)$ via the isomorphism
$H^\bullet(\Omega)\cong U(\a)$
induced by that of Proposition~1.7.2.

We are now ready to give an answer to Question 2 in n.~1.6.

\begin{Th}
\label{1.7.6}
{\it
For any
$u\in H^\bullet(\Gr)$,
the natural action of $u$ on the hyper-cohomology of a perverse sheaf
from the category $P(Gr)$ corresponds (cf.~Remark~1.6.4) to the natural
action of $\check u\in\a$ in the $G^\lor $-module that corresponds to
the perverse sheaf via Theorem~1.4.1.
}\end{Th}

\subsection{Kostant theorem and the generalized exponents}
The imbedding (1.6.1) makes $\g^\lor $ an $\sl_2(\C)$-module
with respect to the adjoint action. One can decompose $\g^\lor $ into
irreducible $\sl_2(\C)$-submodules $\g^\lor _i$ so that the corresponding
highest weight vectors
$\check a_i\in\check\g^\lor _i$
form a base of the subalgebra $\a$
(in particular, the decomposition is known to have $r=\rk(G^\lor )$
direct summands).

Furthermore, all the eigenvalues of the operator $\ad h$ in
$\g^\lor _i$ are {\em even} integers for otherwise we would have:
$\dim\a>\rk G$. The structure theory of $\sl_2$-modules
then yields an equality: 
$\dim\g^\lor _i=\deg \check a_i+1$. On the other hand, proposition
\ref{1.7.2} combined with theorem \ref{1.7.6} implies:
 $\deg \check a_i=\deg\a_i=2\cdot\deg P_i-2$.
Thus, we have proved the following numerical identity
first discovered by Kostant~\cite{Ko1}:
\[
\deg\g^\lor _i=2\cdot\deg P_i-2
\]
\smallskip

Let $\t=Z_{\g^\lor }(h)$ be a Cartan subalgebra of $\g^\lor $ and $T^\lor $
the corresponding Cartan subgroup. Let
$C\subset X^*(T^\lor )$
be a coset with respect to the root lattice of $G^\lor $ and $C^{++}=
\{\mu\in C\,|\,\mu(h)\ge 0\}$, a
 dominant Weyl chamber. There is a unique
weight $\mu_{_C}\in C^{++}$ which is the minimal element in $C^{++}$,
i.e. such that for any other $\lambda\in C^{++}$,
$\lambda-\mu_{_C}$ is a sum of positive
roots. The weight $\mu_{_C}$ is
{\em minuscule}, that is all the weights of the irreducible
representation $V_{\mu_{_C}}$ form a single $W$-orbit. Any minuscule
weight is known to be the minimal element in some coset $C$.

Let $v_\mu$ denote a lowest weight vector in the representation
$V_\mu$. The following result was suggested to me by B.~Kostant.

\begin{prop}
The module $V_{\mu_{_C}}$ is cyclically generated by the action of the
algebra $\a=Z_{\g^\lor }(n)$, that is $V_{\mu_{_C}}=U(\a)\cdot v_{\mu_{_C}}$.
\end{prop}

\begin{proof}{Proof}
The orbit $O_{\mu_{_C}}$ is the unique closed orbit in the connected
component of the Grassmannian $\Gr$ that corresponds to the coset
$C$ (connected components are parameterized by elements of the finite
group $X_*(T)/\mbox{root}$ lattice). Hence, $\overline O_\mu=O_\mu$ is
a smooth variety so that $IC_{\mu_{_C}}$ is a constant sheaf on $O_{\mu_{_C}}$
(up to shift). By Theorem~1.7.2, proving the Proposition amounts
to showing that the restriction map
$H^*(\Gr)\to H^*(O_{\mu_{_C}})$
is surjective. This is equivalent, by duality, to the injectivity
of the map $H_*(O_{\mu_{_C}})\to H_*(\Gr)$. But the orbit
$O_{\mu_{_C}}$ has an even-dimensional cell decomposition which is part
of a similar decomposition of $\Gr$ (Bruhat decomposition). The injectivity
follows.
\end{proof}

For any $V\in\Rep_{_{G^\lor}}$, let $V(\mu)$ denote the weight subspace
corresponding to a weight $\mu\in X^*(T^\lor )$, and $V^\a$ the subspace
annihilated by the subalgebra~$\a$. In chapter 5 we will prove the
following
result.

\begin{prop}
Let $\lambda\in C^{++}$. Then,
$\dim V^\a_\lambda=\dim V_\lambda(\mu_{_C})${\rm;} moreover, there
exists a basis $v_1,\ldots,v_m$ of the weight space $V_\lambda(\mu_{_C})$
and non-negative integers {\rm(}called the generalized exponents{\rm):}
$k_1,\ldots,k_m$ such that the elements{\rm:}
$n^{k_1}\cdot v_1,\ldots,n^{k_m}\cdot v_m$
form a base of $V_\lambda^\a$.
\end{prop}

If $C$ is the root lattice then $\mu_{_C}=0$. In that case the dimension
equality of the first part of the Proposition is proved in~\cite{Ko2}
and the second part is proved in~\cite{Br}. Proposition 1.8.2 is
closely related to some ideas of R.~K.~Brylinski who stated it as
a conjecture.

\renewcommand{\theequation}{\thesubsection}
\renewcommand{\theprop}{\thesubsection}

\subsection{Integrals over Schubert cycles}
We now  describe (partially) the pairing between homology and
cohomology of the Grassmannian. Namely, for any
$u\in H^{2k}(\Gr)$
and any  $L^+G$-orbit $O_\lambda$ (see n.~1.4) such that
$\dim_\C O_\lambda=k$,
we shall give a representation-theoretic formula for
$\langle u,[O_\lambda]\rangle$,
the value of $u$ on the fundamental cycle of the closure of $O_\lambda$.
To that end, let $v_\lambda$ denote a lowest weight vector in
the irreducible representation $V_\lambda$ and $v^\lambda$ a lowest
weight vector of the contragredient representation. Further, let
$\check u$ be the element of $U(\a)$ corresponding to the cohomology
class $u$ via the isomorphism of Corollary~1.7.3\vvv. The element
$\check u$ acts naturally in the representation $V_\lambda$, and we
have
\begin{prop}
$\langle u,[O_\lambda]\rangle=\langle v^\lambda,\check u\cdot v_\lambda
\rangle$.
\end{prop}

\begin{proof}{Proof}
Let
$S:V_\lambda\eqto H^\bullet(IC_\lambda)$
denote the canonical isomorphism provided by theorem~1.4.1, and
$V_\lambda^\pm$
the highest (resp.\ the lowest) weight subspaces in $V_\lambda$.
We'll see later that the isomorphism $S$ sends
$V_\lambda^\pm$
into
$H^{\pm k}(IC_\lambda)$,
the top (resp.\ the lowest) intersection cohomology. We use the
map $S$ and the natural map from intersection cohomology to ordinary
homology to obtain the following commutative digram:
\[
\vcenter{\hbox{\diagram
V^+_\lambda\rrto^{\stackrel{S}{\sim}}&&H^k(IC_\lambda)\rrto^{\sim}&&H_0
(\overline O_\lambda)\\
\\
V^-_\lambda\uuto^{\check u}\rrto^{\stackrel{S}{\sim}})&&H^{-k}(IC_\lambda)
\uuto^{u\cdot}\rrto^\sim
&&H_{2k}(\overline O_\lambda)\uuto^{u\cap}
\enddiagram}}
\label{nexts}
\]

The choice of linear function
$v^\lambda:V^+_\lambda\to\C$
corresponds via the above isomorphism, to a map
$\varphi: H_0(\overline O_\lambda)\to\C$.
Hence, the commutative diagram yields:
\[
\langle u,[\overline O_\lambda]\rangle=
\varphi(u\cap[\overline O_\lambda])=
\langle v^\lambda,\check u\cdot v_\lambda\rangle
\]
\end{proof}

\renewcommand{\theequation}{\thesubsection.\arabic{equation}}
\renewcommand{\theprop}{\thesubsection.\arabic{prop}}

\subsection{Computation of Ext-groups}
We are interested in the groups
$\Ext^\bullet(IC_\lambda,IC_\mu)$,
computed in $D^b(\Gr)$ (but not in the category $P(Gr)$ where
$\Ext$'s are trivial). To find these groups observe first that,
for any
$M,N\in D^b(\Gr)$,
one has, by functoriality, a natural morphism:
\begin{equation}
\Ext^i(M,N)\to\Hom^i_{H^\bullet(\Gr)}(H^\bullet(M),H^\bullet(N))
\end{equation}
where $\Hom$'s on the right are taken in the category of graded
$H^\bullet(\Gr)$-modules and $\Hom^i$ denotes the space of homomorphisms
shifting gradation by~$i$.

Now, let
$M=IC_\lambda$
and
$N=IC_\mu$, $\lambda,\mu\in X_*(T)$.
In that case the right-hand side of (1.10.1) turns, by theorem~1.7.2,
into
$\Hom^i_\a(V_\lambda,V_\mu)$,
so that we obtain a map:
\begin{equation}
\Ext^\bullet(IC_\lambda,IC_\mu)\to\Hom^\bullet_\a(V_\lambda,V_\mu)
\end{equation}

\begin{Th}{\it
If the strata $O_\lambda$ and $O_\mu$ are contained in the same connected
component of the Grassmannian, then the morphism (1.10.2) is an
isomorphism.
}\end{Th}

This theorem follows from the main result of~\cite{Gi3} (see Remark
at the end of~\cite{Gi3}).

We shall now give a reformulation of Theorem~1.10.3 in slightly more
invariant terms, not involving a choice of principal nilpotent.

Let ${\cal N}$ be the cone of nilpotent elements of the Lie algebra $\g^\lor $,
and $\C[{\cal N}]$ the algebra of regular functions on ${\cal N}$. There is a natural
action on ${\cal N}$ of the group
$\C^*\times G^\lor $
(the factor $\C^*$ acts by multiplication and $G^\lor $ by conjugation). Let
${{\cal C}oh}({\cal N})$ denote the abelian category of $\C^*\times G^\lor $-equivariant
coherent $\O_{\cal N}$-sheaves. Any such sheaf is characterized by the
$\C[{\cal N}]$-module of its global sections, for ${\cal N}$ is an affine variety.

To each representation
$V\in\Rep_{_{G^\lor}}$
we attach the free $\O_{\cal N}$-sheaf
$V\otimes_\C\O_{\cal N}$
whose global sections from the free $\C[{\cal N}]$-module
of $V$-valued regular functions on ${\cal N}$.
The $\C^*$-action on $V\otimes_\C\O_{\cal N}$ is induced from that on ${\cal N}$, and the
$G^\lor $-action on $V\otimes_\C\O_{\cal N}$
 arises from the simultaneous action of
$G^\lor $ both on $V$ and on ${\cal N}$. In particular, for
$\lambda\in X_*(T)$
we set
${\cal V}_\lambda=V_\lambda\otimes_\C\O_{\cal N}\in{{\cal C}oh}({\cal N})$.

Theorem 1.10.3 turns out to be equivalent to the following
\begin{prop}
For any strata $O_\lambda$ and $O_\mu$ contained in the same connected
component of $\Gr$ one has a canonical isomorphism{\rm:}
\[
\Ext^\bullet_{D^b(\Gr)}(IC_\lambda,IC_\mu)\cong
\Hom^\bullet_{{{\cal C}oh}({\cal N})}({\cal V}_\lambda,{\cal V}_\mu).
\]
\end{prop}

Recently, Kashiwara-Tanisaki~\cite{KT} proved a version of the Kazhdan-Lusztig
conjecture for affine Lie algebras. The complexes $IC_\lambda$ correspond,
by that conjecture, to irreducible highest weight representations
$L_\lambda$ (of negative level) of the Lie algebra $\hat\g$
(= central extension of the algebra
$\g\otimes\C[z,z^{-1}]$). Correspondingly, the left-hand side
of~(1.10.2) turns out to be isomorphic to
$\Ext^\bullet(L_\lambda,L_\mu)$,
the $\Ext$-group in an appropriate category of $\hat\g$-modules.
Thus, Theorem~1.10.3 gives an expression for
$\Ext^\bullet(L_\lambda,L_\mu)$
in terms of finite-dimensional representations of the dual Lie algebra.
 This is quite surprising,
for there is no apparent connection between representations of $\hat\g$
and those of $\g^\lor $. To make a hint on what an
 explanation of such a connection
might be, notice that the sheaves ${\cal V}_\lambda$, ${\cal V}_\mu$  appearing
in Proposition~1.10.4 are {\it projective} objects of the category ${{\cal C}oh}({\cal N})$.
So, the Proposition manifests some instances of the Koszul duality
(see~\cite{BGS}). Moreover, it is likely that Proposition~1.10.4
is a special case of an extension to affine Lie algebras \cite{BGSo} of the
Koszul duality conjecture for semisimple Lie algebras.
The latter one has been proved by W.~Soergel~\cite{So}.
Our Theorem~1.10.3 is similar in spirit to Soergel's result
(cf.\ also \cite{BGSo}). 

A result somewhat related to  Proposition~1.10.4 was  in effect obtained much
earlier in~\cite{FP} in the context of
modular representations of a semisimple
group $G^\lor $ in positive characteristic. The corresponding analogue
for quantum groups at roots of unity was later proved in [GK].
In the last section of [GK] we proposed a conjectural
interpretation of the intersection cohomology $H^*(IC_\lambda)$
in terms of quantum groups. The latter gives way to translating
theorem~1.10.3 into a purely algebraic claim concerning
simple modules over the finite-dimensional quantum algebra
introduced by Lusztig. Finding an algebraic proof of that
claim, independent of the intersection cohomology
methods,  presents a very interesting problem.

We would like to close this discussion with a 
conjecture that would generalize theorems 1.3.1 and 1.4.1.
Fix an Iwahori(= affine Borel) subgroup $I \subset L^+G$ 
and write 
$P_I(Gr)$ for the category of $I$-equivariant perverse
sheaves on $Gr$ with compact support. This is a non-semisimple
abelian category containig $P(Gr)$ as a subcategory. 
Furthermore,  category
$P_I(Gr)$ is known, through the combination of the works of
Kazhdan-Lusztig \cite{KL2} and Kashiwara-Tanisaki
\cite{KT}, to be equivalent to a (regular block of the) category
of finite dimensional representations of a quantum group at a root
of unity. The following conjecture may be
approached, we believe, using combinatorial results of 
\cite{Lu 5} (I am grateful to Lusztig for pointing out this 
reference).

\renewcommand{\theenv}{}
\begin{env}{\it
\label{conjecture2} (i)
For any $M\in P_I(Gr)$ and $A\in P(Gr)$ we have $M~\ast~A$
$\in~P_I(Gr)$,
so that we get a bi-functor $\ast : P_I(Gr) \times P(Gr) \to P_I(Gr)$.
Moreover,

(ii) This bi-functor corresponds, via the
Kazhdan-Lusztig-Kashiwara-Tanisaki
equivalence, to the standard tensor product of representations of the
quantum group.

(iii) In particular, for any fixed $A\in P(Gr)$, the functor
$M\leadsto M\ast A$
on the (non-semisimple) category $P_I(Gr)$
is {\it exact}.}
\end{env}

Observe that it is essential in the conjecture to convolve
with an object $A\in P(Gr)$
"on the right", i.e to take 
$M\ast A$ and not $A\ast M$.

\subsection{}
The action of the principal nilpotent $n$ in a representation
$V\in\Rep_{_{G^\lor}}$
yields the filtration on $V$ by the kernels of powers of $n$. We
let
$V_i(\mu):=V(\mu)\cap\ker(n^{i+1})$, $i=0,1,2,\ldots$,
denote the induced filtration on a weight subspace $V(\mu)$
(notation of~1.8). Write the Poincar\'e polynomial:
\begin{equation}
P_\mu(V,q)=\sum_{i\ge0}q^{2i}\cdot\dim(V_i(\mu)/V_{i-1}(\mu))
\end{equation}

On the other hand, for any pair of dominant weights $\lambda$ and
$\mu$, Lusztig considered the affine Kazhdan-Lusztig polynomial
$P_{\mu,\lambda}(q)$
and proved (see~\cite{Lu}) the weight multiplicity formula:
$\dim V_\lambda(\mu)=P_{\mu,\lambda}(1)$.
We'll show in section~5 that the Localization theorem for equivariant
cohomology yields the following $q$-analogue of Lusztig's formula,
involving polynomials~(1.11.1).

\begin{Th}{\it
For any dominant $\lambda$ and $\mu$ we have{\rm:}
\[
P_\mu(V_\lambda,q)=q^{\lambda(h)-\mu(h)}\cdot P_{\mu,\lambda}(q^2).
\]
}\end{Th}

A similar result was proved (under certain restrictions) in~\cite{Br}
by totally different means.

\section{Tensor category of Perverse Sheaves}

This section is entirely devoted to the proof of Theorem~1.3.

\subsection{}
There is a natural (1-1)-correspondence between orbits of the
$L^+G$-action on
$\Gr=LG/L^+G$
and orbits of the diagonal $LG$-action on $\Gr\times\Gr$. If these
latter orbits were finite-dimensional we would have been able to
define the Intersection cohomology complexes of their closures.
The abelian category of perverse sheaves on $\Gr\times\Gr$ isomorphic
to direct sums of such complexes would have been clearly equivalent
to the category $P(Gr)$ and would have had a natural convolution
structure defined by the formula (cf.~\cite{Spr}):
\[
M* N=(p_{13})_*(p_{12}^*  M\otimes p^*_{23}N)
\]
where
$p_{ij}:\Gr\times\Gr\times\Gr\to\Gr\times\Gr$
denotes the projection along the factor not named in the subscript.

Unfortunately, all $LG$-orbits in $\Gr\times\Gr$ have infinite dimension.
That is why we were forced in n.~1.2 to use the loop group $\Omega$ in
order to define the convolution $*$. We will now give a purely
algebraic construction of the convolution. Another, more canonical,
algebraic construction based on equivariant technique will be given in
the next chapter.

Let $O_\mu^2$ denote the $LG$-orbit in $\Gr\times\Gr$ corresponding
to an $L^+G$-orbit
$O_\mu\subset\Gr$, $\mu\in X_*(T)$,
and let
$p_{1,2}:\Gr\times\Gr\to\Gr$
denote the first and the second projections. Given
$\lambda,\mu\in X_*(T)$,
set
$O_{\lambda,\mu}=p_1^{-1}(O_\lambda)\cap O_\mu^2$.
Clearly, $O_{\lambda,\mu}$ is a finite-dimensional subvariety of
$\Gr\times\Gr$. Furthermore, the first projection makes
$O_{\lambda,\mu}$
a locally-trivial fibration over $O_\lambda$ with fibre $O_\mu$.
Let $IC_{\lambda,\mu}$ denote the Intersection cohomology complex
on $\overline O_{\lambda,\mu}$, the closure of $O_{\lambda,\mu}$.

\begin{prop}
There is natural isomorphism{\rm:}
$IC_\lambda* IC_\mu\cong(p_2)_* IC_{\lambda,\mu}$.
\end{prop}

\begin{proof}{Proof}
Identify $\Gr$ with the Loop group $\Omega$ and define an automorphism
$f$ of variety
$\Omega\times\Omega$
by the formula
$f:(x_1,x_2)\mapsto(x_1,x_1^{-1}\cdot x_2)$, $x_1,x_2\in\Omega$.
The map $f$ fits into the commutative triangle:
\[
\vcenter{\hbox{\diagram
\Omega\times\Omega\rrto^f_\sim\ddrto_{p_2}&&\Omega\times\Omega\ddlto^m\\
\\
&\Omega&\\
\enddiagram}}
\]
The triangle yields:
$IC_\lambda* IC_\mu=m_*(IC_\lambda\times IC_\mu)=(p_2)_*
\circ f^*(IC_\lambda\times IC_\mu)$.
Hence, proving the Proposition amounts to showing that
$f^*(IC_\lambda\times IC_\mu)\cong IC_{\lambda,\mu}$,
i.e. that
$\overline O_{\lambda,\mu}=f^{-1}(\overline O_\lambda\times
\overline O_\mu)$.
It suffices to check that
$O_{\lambda,\mu}=f^{-1}(O_\lambda\times O_\mu)$.

To prove the last equality we introduce the map
$q:\Omega\times\Omega\to\Omega$
defined by
$q(x_1,x_2):=x_1^{-1}\cdot x_2$.
We claim that
$O_\mu^2=q^{-1}(O_\mu)$.
Indeed, the sets
$O_\mu^2$ and
$q^{-1}(O_\mu)$ are, clearly, stable under the diagonal
$\Omega$-action on $\Omega\times\Omega$. Hence, each of these sets
is completely determined by its intersection with $\Omega\times\{1\}$.
But we have:
\[
O_\mu^2\cap(\Omega\times\{1\})=O_\mu\times\{1\}=q^{-1}(O_\mu)\cap
(\Omega\times\{1\}),
\]
and the claim follows. Thus, one finds:

\[
O_{\lambda,\mu}=
p_1^{-1}(O_\lambda)\cap O^2_\mu=
p_1^{-1}(O_\lambda)\cap q^{-1}(O_\mu)=
f^{-1}(O_\lambda\times O_\mu).
\]

This completes the proof of the Proposition.
\end{proof}

\begin{rem}
We have shown in the course of the proof of Proposition 2.1.1
that the product
$m_*(O_\lambda\times O_\mu)$
of any two orbits
$O_\lambda$ and  $O_\mu$
is a finite union of orbits. It is clear, on the other hand,
that each of the sets
$ \Gr_i$ forming an exhaustion of $ \Gr$ (see n.~1.2)
is a finite union of orbits.
Hence, for any $i, j$, the set
$m_*( \Gr_i\times \Gr_j)$
is contained in a big enough $ \Gr_k$,
as claimed in n.~1.2.
\end{rem}

\subsection{}
The following result implies the first part of Theorem 1.3.
Its proof is based in a very essential way
on a result of Lusztig~\cite{Lu}.
The geometric construction of the Hecke algebra
used in the proof of Proposition 2.2.1 below
is nowadays well known, cf.~[Spr].
It was independently discovered in the course of the proof of
the Kazhdan-Lusztig conjecture by a number of people
(including Beilinson-Bernstein, MacPherson, Lusztig and others).

\begin{prop}
For any $\lambda,\mu\in\XT$,
the complex $ IC_\lambda* IC_\mu$
is isomorphic to a finite direct sum of the complexes
$ IC_\nu$, $\nu\in\XT$.
\end{prop}

\begin{proof}{Proof}
By Proposition 2.1.1 we have:
$ IC_\lambda* IC_\mu=(p_2)_* IC_{\lambda,\mu}$.
The latter complex is isomorphic to a finite direct sum
of shifted Intersection cohomology complexes,
by the Decomposition theorem \cite{BBD}.
Note further, that
the variety $O_{\lambda,\mu}$ is stable under
the diagonal action of the group $ L^+G$.
So, the complex $(p_2)_* IC_{\lambda,\mu}$,
and hence all of its direct summands,
are locally-constant along $ L^+G$-orbits in $ \Gr$.
Hence, these direct summands are of the form
$ IC_\nu[n_\nu]$ for some $\nu\in\XT$,
where $[n_\nu]$ denotes the ''shift'' in the derived category.
Thus, we obtain:

\begin{equation}
 IC_\lambda* IC_\mu=\sum_\nu IC_\nu[n_\nu],
\quad\nu\in\XT
\end{equation}

Recall next that associated to the semisimple group $G$ is
the  {\it affine Hecke algebra} $H$, cf.~[Spr]. This is a
$\Z[t, t^{-1}]$-algebra
which has a base formed by
certain distinguished elements $c_w$, $w\in W_a$
(= affine Weyl group),
the so-called Kazhdan-Lusztig basis.
The algebra $H$ is known to have
a geometric realization as the Grothendieck group of
a semisimple category whose objects are direct sums of
certain (shifted) Intersection cohomology complexes
$ I_w[n]$, $w\in W_a$, $n\in\Z$,
on the Flag variety associated to an affine Lie algebra.
The complex $ I_w[n]$ corresponds to the element
$t^n\cdot c_w\in H$.
The ring structure on $H$ corresponds
to a convolution of complexes.

There is a natural inclusion $\lambda\mapsto w(\lambda)$
of the set $\XT/W$ into $W_a$.
Furthermore, the assignment
$ IC_\lambda\mapsto I_{w(\lambda)}$ is compatible
with convolutions so that one has the formula:

\begin{equation}
 I_{w(\lambda)}* I_{w(\mu)}=\sum_\nu I_{w(\nu)}[n_\nu],
\end{equation}
corresponding to (2.2.2) term by term.

Now, Lusztig has proved \cite[Corollary 8.7]{Lu}  that
for any $\lambda,\mu\in\XT$
one has an equality in $H$:

\begin{equation}
c_{w(\lambda)}\cdot c_{w(\mu)}=
\sum_\nu m^\nu_{\lambda, \mu}\cdot c_{w(\nu)},
\quad\nu\in\XT
\end{equation}
where
$m^\nu_{\lambda, \mu}$ are some non-negative integers
(they are equal to the multiplicity of
the irreducible representation $V_\nu$ in $V_\lambda\otimes V_\mu$).
What is important for us in (2.2.4) is that
the coefficients $m^\nu_{\lambda, \mu}$ are independent of $t$.
This implies that there are no ''shifts'' in the
right-hand side of (2.2.3), i.e.~$n_\nu=0$ for all $\nu$.
Hence, the same is true for the right-hand side of (2.2.2)
and the Proposition follows.\endpr

\subsection{}

We will show now that the convolution on $ P(Gr)$ is commutative.
More precisely, one has

\begin{prop}
There is a natural functor isomorphism:
\[
A*B\cong B*A,\quad A, B\in\Ob\  P(Gr).
\]
\end{prop}

\medskip\noindent Proof of the Proposition is essentially standard.
Let $\theta$ be a Cartan anti-involution on $ G$
such that $\theta( K)= K$ where $ K$ is a maximal compact subgroup of $ G$
(choosing $\theta$ amounts to choosing a maximal torus $ T\subset G$
such that $ T\cap K$ is a maximal torus in $ K$).
The induced (pointwise) anti-involution on
$ LG$ will be also denoted by $\theta$.
By definition, the anti-involution $\theta$
preserves the subgroup $\Omega$ and acts identically on the set $\XT$,
viewed as a subset of $ LG$.
It is clear also that $\theta( L^+G)= L^+G$.
Hence, for any $\lambda\in\XT$, the map
$\theta$ preserves the set
\[
( L^+G\cdot\lambda\cdot L^+G)\cap\Omega.
\]
But this set may be identified with the stratum $O_\lambda$ in $ \Gr$.
So, we have an isomorphism
$\theta^* IC_\lambda\cong IC_\lambda$
for any $\lambda\in\XT$ and, hence, an isomorphism
$\theta^*A\cong A$ for any $A\in P(Gr)$.

The isomorphism $\theta^*A\cong A$ above is not canonical, in general.
To choose a canonical one we argue as follows.
For each $\lambda\in\XT$, viewed as a point in $ \Gr$,
let $i_\lambda: \{\lambda\}\hookrightarrow \Gr$
denote the inclusion. Since $\lambda$ is fixed by $\theta$,
for any complex $A\in P(Gr)$, there is
a {\em canonical\/} isomorphism
$ I_\lambda: i_\lambda^*A\stackrel{\sim}{\to} i_\lambda^*(\theta^*A)$.
We now define a
{\em canonical isomorphism\/}
$ I: \theta^*A\stackrel{\sim}{\to} A$
to be the unique isomorphism
such that for any $\lambda\in\XT$ the composite:
\[
i_\lambda^*A\buildrel I_\lambda\over\longrightarrow i_\lambda^*(\theta^*A)
\buildrel I\over\longrightarrow i_\lambda^*A
\]
is the identity morphism on $i_\lambda^*A$.

Finally, we define a functor isomorphism $A*B\cong B*A$ to be
the composite of isomorphisms:
\[
A*B\buildrel I^{-1}\over\longrightarrow
(\theta^*A)*(\theta^*B)\stackrel{\sim}{\to}
\theta^*(B*A)
\buildrel I\over\longrightarrow B*A,
\]
where the isomorphism in the middle arises from the fact that $\theta$
is an anti-involution on the group $\Omega$.
\endpr\medskip

It is clear that the complex $ IC_e$ supported on $e$,
the unit of the group $\Omega$, is the unit
in the category $ P(Gr)$, i.e., for any
$M\in P(Gr)$, there are functorial isomorphisms:
$M* IC_e\cong M\cong IC_e*M$.
Furthermore, we will show in the next chapter
that the associativity constraint and
the above defined commutativity constraint satisfy
the hexagon axiom \cite[(1.0.2)]{DM}.
That will prove that $( P(Gr), *)$ is a semisimple abelian
tensor category.

\subsection{}
Let $\sigma$ be an involutive automorphism of the Grassmannian
arising from the involution: $x\mapsto x^{-1}$ on the group $\Omega$.
Let $\D$ denote the Verdier duality functor on $\D^b( \Gr)$.
Define the {\em transposition} functor on $\D^b( \Gr)$ by
$A^t:=\sigma^*(\D A)$, $A\in\D^b( \Gr)$.

For any $M, N, L\in\D^b( \Gr)$ one has natural isomorphisms:

\begin{equation}
 \begin{array}{l}
  \mbox{\rm (i)}\quad (M^t)^t\cong M;\\
  \mbox{\rm (ii)}\quad (M*N)^t\cong N^t*M^t;\\
  \mbox{\rm (iii)}\quad \RHom(M*N, L)\cong\RHom(M, L*N^t).
 \end{array}\hbox{\hfill}
\end{equation}

Extending $\sigma$ to the anti-involution $x\mapsto x^{-1}$
on the group $ LG$, we observe that:
\[
\sigma(\Omega\cap( L^+G\cdot\lambda\cdot L^+G))=
\Omega\cap( L^+G\cdot\lambda^{-1}\cdot L^+G)
\]
for any $\lambda\in\XT$.
It follows that one has an isomorphism:
$( IC_\lambda)^t\cong IC_{\lambda^{-1}}$.
Hence, the category $ P(Gr)$ is stable under the transposition $M\mapsto M^t$.

Further, given $M, N\in P(Gr)$ define an object
${\cal H}om(M, N)\in P(Gr)$ by the formula:
${\cal H}om(M, N):=M^t*N$.

\begin{prop}
For any $M, N, L, T\in P(Gr)$ we have:
\begin{flushleft}
{\em (i)} $M^t\cong {\cal H}om(M,  IC_e)$;\\
{\em (ii)} The functor: $T\mapsto\Hom(T*N, L)$ is represented
by the object ${\cal H}om(N, L)$, that is:
\[
\Hom_{ P(Gr)}(T*N, L)\cong\Hom_{ P(Gr)}(T, {\cal H}om(N, L)).
\]
\end{flushleft}
\end{prop}

\begin{proof}{Proof}
(i) follows from definitions. To prove (ii) note that for any objects
$X, Y\in P(Gr)$ we have:
$\Hom_{ P(Gr)}(X, Y)\cong\Hom_{\mbox{\rm\scriptsize{\D}}^b( \Gr)}(X, Y)$.
Hence, the isomorphism (2.4.1 (iii)) yields:
$\Hom_{ P(Gr)}(T*N, L)\cong\Hom_{ P(Gr)}(T, L*N^t)$. But
$L*N^t\cong N^t*L$, by the commutativity of the convolution.
Thus, $L*N^t= {\cal H}om(N, L)$ and statement (ii) follows.
\end{proof}

\begin{lemma}
There is natural functor isomorphism:
\[
{\cal H}om(M_1, N_1)*{\cal H}om(M_2, N_2)\cong {\cal H}om(M_1*M_2, N_1*N_2)
\]
for $M_i, N_i\in P(Gr)$.
\endpr
\end{lemma}

Proposition 2.4.2, Lemma 2.4.3 and the results of n.~3.7 below
show that $ P(Gr)$ is an abelian rigid tensor category.

\subsection{}

Let us make a general remark. Let $f: X\rightarrow Y$ be a morphism
of topological spaces and $J$ a constructible complex on $X$.
Then, one has a natural isomorphism of hyper-cohomology:

\begin{equation}
H^\bullet(Y, f_*J)\cong H^\bullet(X, J).
\end{equation}

This isomorphism is obtained as the composite of the following isomorphisms:
\[
H^\bullet(Y, f_*J)
\cong\R^\bullet p_*(f_*J)
\cong\R^\bullet(p\circ f)_* J
\cong H^\bullet(X, J)
\]
where $p$ denotes the constant map: $Y\rightarrow\pt$.

We shall now prove the third claim of theorem 1.3.1.

The exactness of the functor $H$ is obvious since $ P(Gr)$
is a semisimple category.
Next, let $M, N\in P_\bullet(G)$. We have:
\begin{tabbing}
$H^*(\Omega, M)\otimes H^*(\Omega, N)$\quad\=$\cong$\quad\=\kill
$H^*(M*N)$\>$\cong$\>   (definition of the convolution)\\
$H^*(\Omega, m_*(M\boxtimes N))$\>
$\cong$\>
((2.5.1) applied to $J=M\boxtimes N)$\\
$H^*(\Omega\times\Omega, M\boxtimes N)$\>$\cong$\>(Kunneth formula)\\
$H^*(\Omega, M)\otimes H^*(\Omega, N).\quad\square$ \> \> 
\end{tabbing}

\section{Application of the equivariant cohomology.}

Throughout this chapter the reader is assumed to be familiar
with definitions and results on equivariant hyper-cohomology,
proved in [BL] (see also [Lu 3], [Lu 4]) and summarized in chapter 8 below.

\subsection{Convolution construction}
Let $M$ be a Lie group, $L\subset M$ a closed subgroup and
$X$ an $M$-variety. Let $D^b_L(X)$ be the $L$-equivariant derived
category on $X$, cf. \cite{BL}, and 
write $\dm$ for the $L$-equivariant derived
category of complexes on $M/L$ with compactly supported
cohomology sheaves. In this setup, we are going to define a 
convolution pairing
\begin{equation}
\label{con1}
\ast : \dm \times D^b_L(X) \to D^b_L(X)
\end{equation}
The construction of this pairing
 is essentially due to Lusztig, who considered the
special
case $X=M/L$.

First, define a `{\it right-left}' $L$-action on $M\times X$ by the formula:
$l : (m,x) \mapsto (m\cdot l^{-1}, l\cdot x)$. The right-left 
action is
clearly a free left action, and we write $M\times_L X$ for the corresponding
orbit space. Observe that $M$-{\it action-map} $M\times X \to X$
descends to a well-defined map
$a: M\times_L X \to X$. 

Second, let $p : M \to M/L$ be the projection. Then the pull-back
functor $p^*$ gives, by equivariant descent property, cf. $\S 8.2$,
 an equivalence between $D^b_L(M/L)$ and the
equivariant derived category on $M$ with respect to the
$L$-action by right translation. In particular, we get an
equivalence of $\dm$ with the equivariant derived category of
complexes
on $M$ with an additional support condition. Observe also that, for
a subset  $S\subset M/L$, the restriction map
\begin{equation}
\label{conv2}
a: p^{-1}(S)\times_L X \to X\quad\mbox{\it is proper whenever $S$ is compact.}
\end{equation}
Now, let $\mm\in \dm$ and $\aaa\in
D^b_L(X)$. By the second remark above, we view $p^*\mm$ as
an $L$-equivariant complex on $M$ with respect to right
translation. Then,
we may regard $(p^*\mm)\boxtimes\aaa$ as an element
of $D^b_{L_{rl}}(M\times X)$, the equivariant derived
category on $M\times X$ with respect to the right-left 
$L$-action. The latter is equivalent, by the equivariant
descent, to $D^b(M\times_L X)$. Thus, there is a uniquely
determined complex ${\overline {\aaa\boxtimes\mm}}$ on
$M\times_L X$ such that its pull-back to $M\times X$
is isomorphic to $(p^*\mm)\boxtimes\aaa$. The support
condition on $\mm\in \dm$ combined with \ref{conv2}
ensure that the restriction of the map $a$ to the support of
any cohomology sheaf $\cal H^\bullet({\overline {\aaa\boxtimes\mm}})$
is proper. We put 
\begin{equation}
\label{conv3}
\mm\ast\aaa=a_*\left({\overline {\aaa\boxtimes\mm}}\right)=
a_!\left({\overline {\aaa\boxtimes\mm}}\right)
\end{equation}
It is easy to verify that $\mm\ast\aaa \in D^b_L(X)$, and that one has
functorial isomorphisms
\[(\mm_1\ast\mm_2)\ast\aaa = \mm_1\ast (\mm_2 \ast\aaa)\quad,
\;\;\forall \mm_1\,,\,\mm_2 \in \dm\,,\,\aaa\in D^b_L(X)\]

In the next subsection we will adapt the above construction
to the infinite dimensional
set-up: $M=LG\,,\, L=L^+G$ and $X=LG/L^+G$
to get an algebraic definition of the convolution on the
category $P(Gr)$. With similar modifications, 
 the above
construction
may be applied,
 using the formalism of $\S 8.2$,
in the case: $M=LG\,,\, L=L^+G$ and
$X=G(\K_{\,x})/G(\oo_{\text{out}})$ to give a rigorous
definition of the Hecke operators exploited in $n^\circ 1.5$.
The condition allowing reduction to finite dimensions in
 either of those cases is essentially the following:
  there exists a normal subgoup
$L'\subset L$ such 
that~: (i) {\it The quotient $L/L'$ is a finite-dimensional Lie
group};
 (ii) {\it  The group $L'$ acts trivially on $X$}.


\subsection{}

Recall the subgroups $L^i\,\subset\,L^+G$, $i=0,1,\ldots$
introduced in $(1.2.1^m)$.
It follows from proposition 1.2.2(iii) that, for each $ \Gr_i$,
the $ L^+G$-action on $ \Gr_i$ factors through the action
of a finite-dimensional algebraic group  $L_k:= L^+G/L^k$, $k=k(i)$.
 Thus, we may speak of
$ L^+G$-equivariant sheaves on $ \Gr_i$,
meaning the sheaves that are equivariant
with respect to the group $L_k$ where $k$ is chosen so that the $L^k$-action
on $ \Gr_i$ is trivial. That does not depend on a particular
choice of~$k$.

We can now construct the convolution $\ast$ on $P(Gr)$ 
using the approach of the previous $n^\circ$ as follows
Let $\mm\,,\, \aaa\in P(Gr)$. Choose positive integers $i$ and $j$
so that the complexes $\mm$ and $\aaa$ are supported on $ \Gr_i$ and $ \Gr_j$
respectively. Let $ LG_i$ denote the inverse image of $ \Gr_i$
under the projection $ LG\rightarrow LG/ L^+G=\, \Gr$.
Clearly $ L^+G\subset LG_i$. Moreover, by n.~1.2,
there exists an integer $m=m(i,j)>> 0$
 such that $ LG_i\cdot \Gr_j\subset \Gr_m$; furthermore, there exists
another integer, $k=k(i,j)>> 0$, such that the group $L^k$ acts trivially
on the subsets $\Gr_j$ and $ \Gr_m$. 
We can therefore replace in the setup of $n^\circ 3.1$ the action
map $M \times X \to X$ by the following well defined morphism of
{\it finite-dimensional} varieties induced by the action
\begin{equation}
\label{act1}
LG_i/L^k\,\times \Gr_j\rightarrow \Gr_m
\end{equation}
Recall the group $L^k$ is normal in $L^+G$,
and set $L_k:=L^+G/L^k$, a finite dimensional Lie group.
There is a natural free
$L_k$-action on $LG_i/L^k$ on the right making 
the projection
$ LG_i/L^k\rightarrow LG_i/ L^+G=\, \Gr_i$
a principal $L_k$-bundle over $
\Gr_i$. Thus, the `right-left' $L_k$-action
on the product $LG_i/L^k\,\times \Gr_j$
is free, and the orbit-space
\[LG_i/L^k\,\times_{L_k} \Gr_j = LG_i\,\times_{_{L^+G}} \Gr_j\]
is a well-defined finite dimensional variety. Moreover, the
map \ref{act1} induces the following morphism
\[a: LG_i/L^k\,\times_{L_k} \Gr_j\,\to\, Gr_m\]
which is an
analogue of the map $a: M\times_L X \to X$ from $n^\circ 3.1$.
With that understood, the convolution-construction
of $n^\circ 3.1$ goes through verbatim.
One can show that this way we obtain the same convolution on the category $ 
P(Gr)$
as the one defined in n.~1.2 (for $\mm= IC_\lambda$ and $\aaa= IC_\mu$
our present construction is identical to that of n.~2.1 (Proposition 2.1.1)).

\subsection{}

Now let $K$ be any {\em connected} subgroup of $ L^+G$,
that is a subgroup that projects to a connected Lie subgroup $K_i\subset L_i$,
for all large enough $i\gg0$.
Then, by (8.2.4), the
constant sheaf on each $ L^+G$-orbit $O_\lambda$ has
a unique structure of $K$-equivariant complex.
Hence, the Intersection cohomology extension of such a complex
has a unique $K$-equivariant structure.
Observe next that
any object of the category $ P(Gr)$ is a direct sum of such complexes and,
moreover, the direct sum can be chosen so that there are no
non-trivial morphisms between different summands.
Thus, we have proved the following

\begin{prop}
\label{connectgrass}
For any {\em connected} subgroup $K\subset L^+G$ as above,
every object of the category $ P(Gr)$ has a unique structure of
$K$-equivariant perverse sheaf.
\endpr
\end{prop}

We apply the proposition  to a maximal torus $ T\subset G$,
viewed as a subgroup of constant loops. We get


\begin{corol}
Any object of the category $ P(Gr)$ has a canonical (unique)
structure of $ T$-equivariant complex on $ \Gr$.
\endpr
\end{corol}

Recall that any algebraic homomorphism $\lambda:\C^*\rightarrow T$
gives a point in $ \Gr$ denoted by the same symbol $\lambda$.
These points form a countable  descrete  subset  $\XT\subset  \Gr$ 
that turns out
to be exactly the set of $ T$-fixed points.

\begin{lemma}
$ \Gr^{ T}=\XT$.
\end{lemma}

To prove the Lemma, choose a maximal compact subgroup $ K\subset G$
such that $ T_c:= K\cap T$ is a maximal (compact) torus in $ K$.
Identify the Grassmannian with the Loop group $\Omega$ as in n.~1.2.
The action of the compact torus $T_c$ on $ \Gr$ corresponds
to the natural $T_c$-action on $\Omega$ by conjugation.

\begin{proof}
{Proof of the Lemma}
It is clear that $ \Gr^{ T}= \Gr^{ T_c}\cong\Omega^{ T_c}$. Further, a loop
$f\in\Omega$ is fixed by the $ T_c$-action iff the image of $f$ is contained
in $ T_c\subset K$. But any {\em polynomial\/} loop: $S^1\rightarrow T_c$
must be a group homomorphism. Hence, $f\in\XT$. The inverse inclusion:
$\XT\subset \Gr^{ T}$ is obvious.
\end{proof}

\subsection{}
We now apply the machinary of equivariant cohomology.
To any perverse sheaf $M\in P(Gr)$ we associate the space $H^*_{ T}(M)$,
the $ T$-equivariant cohomology group with respect to the
$ T$-equivariant structure on $M$ provided by Corollary 3.3.1.
The assignment: $M\mapsto H^*_{ T}(M)$ gives rise to a functor:
\[
H^*_{ T}: P(Gr)\rightarrow\mod^\bullet\mbox{\rm-}\Ct,
\]
where $\mod^\bullet\mbox{\rm-}\Ct$ denotes the tensor category
of finitely-generated graded $\Ct$-modules.

\begin{prop}
The functor $H^*_{ T}$ is a fibre functor.
\end{prop}

\noindent
(Proposition says that $H^*_{ T}$ is an exact fully-faithful tensor
functor such that $H^*_{ T}(M)$ is a projective $\Ct$-module
for any $M\in P(Gr)$).

\begin{proof}{Proof}
The functor $H^*_{ T}$ is obviously exact. Theorem 8.4.1 implies that it is
fully-faithful and that the module $H^*_{ T}(M)$ is free for any $M\in P(Gr)$.
It remains to show that $H^*_{ T}$ is a tensor functor. This can be
done by repeating the argument of n.~2.4 as follows.

Identify $ \Gr$ with the Loop group $\Omega$ and replace
$ T$-equivariant cohomology by $ T_c$-equivariant cohomology (see 8.3.4).
By abuse of notation, we shall drop the subscript ''c'' and will always
write ''$ T$'' instead of ''$ T_c$''.
Given $M, N\in P(Gr)$, set $J=M\boxtimes N$.
This is a $ T$-equivariant complex on $\Omega\times\Omega$.
To such a complex $J$ we have associated in n.~3.1 a complex $J_{ T}$
on the variety
\mbox{$(\Omega\times\Omega)_{ T}=\mbox{\rm ET}\times_{ 
T}(\Omega\times\Omega)$}.

Further, note that the compact torus acts on the Loop group $\Omega$ by group 
automorphisms,
so that the multiplication map $m:\Omega\times\Omega\rightarrow\Omega$
commutes with the $ T$-action.
Hence, the functor $J\mapsto J_{ T}$ commutes with the direct
image functor $m_*$. Hence, we obtain:
\begin{equation}
 \begin{array}{rcl}
  H_{ T}^*(M*N)=H_{ T}^*(m_*(M\boxtimes N))&=& \\
  H_{ T}^*(m_*J)=H^*(\Omega_{ T}, (m_*J)_{ T})&=& \\
  H^*(\Omega_{ T}, (m_{ T})_*(J_{ T}))&=&\mbox{\rm (2.4.1)} \\
  H^*((\Omega\times\Omega)_{ T}, J_{ T})=H_{ T}^*(M\boxtimes N)&=&\mbox{\rm 
(Corollary 8.4.3)} \\
  H_{ T}^*(M)\otimes_{\Ct}H_{ T}^*(N).&&
 \end{array}
\end{equation}
\end{proof}

We can specialize the equivariant cohomology functor $H^*_{ T}$
at various points $t\in\gotht$. By Proposition 3.4.1, we obtain
a family of exact fully-faithful tensor functors $H_t: P(Gr)\rightarrow\Vect$.
All these functors are non-canonically isomorphic to each other
(corollary 7.2.2) and, in  particular, isomorphic to
the ordinary cohomology functor $H^*(\cdot)\cong H^*_0(\cdot)$ (cf.~n.~8.5).

\subsection{}
\setcounter{subsection}{6}

Recall that $ \Gr^{ T}$, the fixed point subvariety, is the discrete
set consisting of isolated points $\lambda\in\XT$.
Let $i_\lambda:\{\lambda\}\hookrightarrow \Gr$ denote the inclusion.
Fix a regular $t\in\gotht$.
For any complex $M\in P(Gr)$, the Localization theorem 8.6 yields an 
isomorphism:
\begin{equation}
\label{fpd}
H_t(M)\cong\mathop\oplus\limits_{\lambda\in\scriptXT}H_t(i_\lambda^!M)
\end{equation}
This isomorphism will be referred to as the {\em fixed point decomposition}
and will be viewed as a gradation on $H_t(M)$ by the lattice $\XT$.

\begin{prop}
The fixed point decomposition is compatible (in the sense of (7.1.2))
with the convolution $*$ on $ P(Gr)$.
\end{prop}

\begin{proof}{Proof}
We use the notation $\Omega^2:=\Omega\times\Omega$. Given $\lambda\in\XT$,
set:
\begin{equation}
(\Omega_\lambda^2)^{ T}:=\{(\mu, \nu)\in\Omega^2\,|\,\mu, \nu\in\XT,\ 
\mu+\nu=\lambda\}
\end{equation}
and let $i^2_\lambda:(\Omega^2_\lambda)^{ T}\hookrightarrow\Omega^2$
denote the natural inclusion.
Proving the Proposition amounts to showing that,
for each $\lambda\in\XT$, the image of the push--forward morphism
\begin{equation}
(i_\lambda)_!: H_t(i_\lambda^!(M*N))\rightarrow H_t(\Omega, M*N)
\end{equation}
is identified -- via the isomorphisms (3.6.5) -- with
the image of the similar morphism:
\begin{equation}
(i^2_\lambda)_!: H_t(i^2_\lambda)^!(M\boxtimes N)\rightarrow H_t(\Omega^2, 
M\boxtimes N)
\end{equation}

In order to compare (3.6.4) with (3.6.5) we introduce the subvariety
$\Omega^2_\lambda:=m^{-1}(\lambda)\subset\Omega^2$ where $m$ is the
multiplication
map. Clearly, $\Omega^2_\lambda$ is a $ T_c$-stable subvariety of $\Omega^2$
and we have an equivariant Cartesian square:
\[
\vcenter{\hbox{\diagram
\Omega_\lambda^2\rrto^m\ddto_j|<{\!\!\!\cap}&&\{\lambda\}
\ddto^{i_\lambda}|<{\!\!\!\cap}   \\
\\
\Omega^2\rrto_m&&\Omega\\
\enddiagram}}
\label{seconds2}
\]
The Base change theorem for this square yields an isomorphism:
\begin{equation}
\mbox{\rm Image}(i_\lambda)_!\simeq
\mbox{\rm Image}
\big[
H_t(\Omega^2_\lambda, j^!(M\boxtimes N))\mathop{\longrightarrow}
\limits^{j_!}
H_t(\Omega^2, M\boxtimes N)
\big]
\end{equation}

Next, note that the set of $ T_c$-fixed points in
$\Omega^2_\lambda$ coincides with the set (3.6.3)
(this agrees with the notation $(\Omega^2_\lambda)^{ T}$).
We have the following commutative triangles of maps
\begin{equation}
\label{1}
\diagram
(\Omega_\lambda^2)^T\rrto^\varepsilon\ddrto_{i_\lambda^2}&&\Omega_\lambda^2
\ddlto^j\\
\\&\Omega^2&
\enddiagram
\end{equation}
That induces the corresponding
morphisms
on cohomology:
\begin{equation}
\label{18b1}
\diagram
H_t(i_\lambda^2)^!(M\times
N)\rrto^{\varepsilon_!}_\sim\drto_{(i_\lambda^2)_!}
&&H_t(\Omega_\lambda^2,j^!(M\boxtimes N))\dlto^{j_!}\\
&H_t(\Omega^2,M\boxtimes N)&
\enddiagram
\end{equation}
The horisontal arrow $\varepsilon_!$ in the latter triangle is an isomorphism
by localization theorem 3.3. So, we obtain from (3.6.6) and (3.6.7) that
\begin{equation}
\Im(i_\lambda)_!\cong\Im(j_!)\cong\Im(j_!\circ\varepsilon_!)\cong\Im(i_\lambda
^2)_!
\end{equation}
The composition of the isomorphisms above coincides,
in the case under consideration, with
isomorphism (7.1.2). That completes the proof of Proposition 3.6.2.
\end{proof}

\subsection{}

We will now use Proposition 3.6.2 to show that the associativity
and commutativity constraints (see n.~2.3) satisfy the hexagon axiom,
i.e., for any $X,Y,Z\in P(Gr)$ the following natural diagram
commutes (cf.~\cite{DM}):
\begin{equation}
\vcenter{\hbox{\diagram
X*(Y*Z)\dto_\sim\rto^\sim&(X*Y)*Z\rto^\sim&Z*(X*Y)\dto_\sim\\
X*(Z*Y)\rto^\sim&(X*Z)*Y\rto^\sim&(Z*X)*Y\\
\enddiagram}}\label{18b2}
\end{equation}

Assume the notation and assumptions of n.~3.6, in particular,
$t\in\gotht$ is a regular element.

\begin{lemma}
Let $f$ be an endomorphism of an object $M\in P(Gr)$,
such that the induced morphism on cohomology
$H_t(f): H_t(M)\rightarrow H_t(M)$ is the identity morphism.
Then $f$ is the identity morphism.
\end{lemma}

\begin{proof}{Proof}
Decompose $M$ into the direct sum of disjoint isotypical components
$M=M_1+\ldots+M_r$. Clearly, $f$ maps each $M_i$ into itself
and the induced morphism on cohomology is the identity.
Hence, we may assume without loss of generality that $M$ is
isotypical. Then, proving that $f$ is the identity, suffices it to show that 
for some
point $\lambda$ in the nonsingular locus of
$\supp(M)$ the induced morphism
$H_ti^!_\lambda(f):  H_t(i^!_\lambda M)\rightarrow H_t(i^!_\lambda M)$
is the identity. But the latter assertion follows from
the assumptions of the Lemma and the Fixed--point decomposition (3.6.1).
\end{proof}

\begin{corol}
Let $f, g: M\rightarrow N$ be two isomorphisms in $ P(Gr)$ such that
$H_t(f)=H_t(g)$. Then $f=g$.
\endpr
\end{corol}

Let $u: A*B\stackrel{\sim}{\to} B*A$ be the commutativity (functor) 
isomorphism
defined in Proposition 2.3.1 and let $H_t(u):H_t(A*B)\stackrel{\sim}{\to} 
H_t(B*A)$
be the corresponding isomorphism on cohomology.
Further, let
$v: H(A*B)\stackrel{\sim}{\to} H(A)\otimes H(B)$ denote the functor 
isomorphism arising from Proposition 3.5.4.
Finally, for any vector spaces $V$ and $W$, let $s:V\otimes 
W\stackrel{\sim}{\to} W\otimes V$
denote the standard isomorphism.

\begin{lemma}
For any $A, B\in P(Gr)$ the following natural diagram commutes
\[
\vcenter{\hbox{\diagram
H_t(A*B)\ddto_{H_t(u)}\rrto^v&&H_t(A)\otimes H_t(B)\ddto^s\\
\\
H_t(B*A)\rrto^v&&H_t(B)\otimes H_t(A)\\
\enddiagram}}
\label{3.7.4}
\]
\end{lemma}

\begin{proof}{Proof}
By the fixed-point decomposition and Proposition 3.6.2, proving the Lemma 
amounts
to showing that, for each $\lambda\in\XT$, the following
diagram commutes:
\[
\vcenter{\hbox{\diagram
H_t(i^!_\lambda(A*B))\ddto_{H_t(i_\lambda^!u)}\rrto^{v_\lambda}&&
\oplus_{\mu+\nu=\lambda}H_t(i_\mu^!A)\otimes H_t(i_\nu^!B)\ddto^s\\
\\
H_t(i^!_\lambda(B*A))\rrto^{v_\lambda}&&
\oplus_{\nu+\mu=\lambda}H_t(i^!_\nu B)\otimes H_t(i^!_\mu A)\\
\enddiagram}}\label{19b}
\]
where the horizontal isomorphisms $v_\lambda$ are the compositions
of the isomorphisms arising from (3.6.8). The commutativity
of this latter diagram follows immediately from the very definition of the 
morphism
$u$ and the construction of the isomorphisms (3.6.8).
\end{proof}

\begin{proof}
{Proof of the hexagon axiom}
We have to show that diagram (3.7.1) commutes.
By Corollary 3.7.3, this amounts to showing that the corresponding diagram of
the
cohomology groups  $H_t(\cdot)$ commutes. Since $H_t(\cdot)$ is a tensor functor,
this latter diagram can be rewritten as follows
(we write $H$ instead of $H_t$ for short):
\[{\small
\vcenter{\hbox{\diagram
H(X)\otimes(H(Y)\otimes H(Z))\rto^b\dto_a&
(H(X)\otimes H(Y))\otimes H(Z)\rto^c& H(Z)\otimes(H(X)\otimes H(Y))\dto^d\\
H(X)\otimes(H(Z)\otimes H(Y))\rto^e&
(H(X)\otimes H(Z))\otimes H(Y)\rto^f&(H(Z)\otimes H(X))\otimes H(Y)\\
\enddiagram}}\label{19b2}}
\]
The maps $b$, $d$ and $e$ in the diagram are the standard associativity
isomorphisms for the tensor product of vector spaces.
The maps $a$, $c$ and $f$ arising from the commutativity constraint in the
category $ P(Gr)$ coincide, by Lemma 3.7.4, with standard isomorphisms arising 
from
the commutativity constraint in the category of vector spaces.
Thus, the diagram above is nothing but the hexagon diagram
in the tensor category of vector spaces.
Hence, this diagram commutes and the commutativity of (3.7.1) follows.
\end{proof}

This completes the proof of Theorem 1.3.1.

\subsection{}
The rest of this chapter is devoted to the  proof of Theorem 1.4.1.

From now on we fix a regular element $t\in\gotht$. The argument of n.~1.4
shows that there is  a reductive group $G^*$ such that the tensor category $ 
P(Gr)$
is equivalent to the tensor category $\Rep(G^*)$ in such a way that the
functor
$H$ goes into the forgetful functor: $\Rep(G^*)\to\Vect$.

The main result of this section is the following

\begin{Th}{\it
The group $ G^*$ is isomorphic to $ G^\vee$, the dual group.
}\end{Th}

Theorem 1.4.1 is an immediate consequence of Theorem 3.8.1 since
the ordinary cohomology functor $H^*$ is isomorphic to $H_t$.

\begin{lemma}
$ G^*$ is a connected semisimple Lie group.
\end{lemma}

\begin{proof}{Proof}
For any $\lambda\in\XT$ the sequence of perverse sheaves:
$IC_\lambda,$ $IC_\lambda* IC_\lambda,$ 
$IC_\lambda* IC_\lambda* IC_\lambda,\ldots$
obviously has the infinite set $\{ IC_{\lambda^n},\ n=1,2,\ldots\}$
among their irreducible constituents. The general criterium
\cite[Corollary 2.22]{DM}
shows now that the group $ G^*$ is connected. This group must be
semisimple, for it has the only one $1$-dimensional representation.
\end{proof}

Let $ T^\vee$ be the torus dual to $ T$, so that we have
a canonical identification: $\XT=\XTup$.
The fixed point gradation \ref{fpd} yields,
by corollary 8.1.3 and proposition 3.6.2, an algebraic homomorphism:
\begin{equation}
f:  T^\vee\longrightarrow G^*.
\end{equation}
\begin{lemma}
The homomorphism $f$ is injective and its image is a maximal torus in $ G^*$.
\end{lemma}

\begin{proof}{Proof}
To prove injectivity, it suffices to show that for any
$\lambda\in\XT$ there is a perverse sheaf $M\in P(Gr)$ such that
$i^!_\lambda M\ne0$. But one obviously has:
$i^!_\lambda( IC_\lambda)\ne0$.

The first assertion of the Lemma being proved,
to prove the second suffices it to show that
$\rk( G)=\rk( G^*)$.
To that end, observe that, for any abelian tensor category $ \Cat$,
its Grothendieck group $K( \Cat)$ has a natural ring structure
induced by the tensor product. Let $Q( \Cat)$ denote the field
of fractions of the ring $K( \Cat)$. Set:
\[
\rk( \Cat):=\DegTrans(Q( \Cat)/\Bbb Q)
\]

It is clear that $\rk( P(Gr))=\rk( G)$. On the other hand,
for any semisimple group, and hence for $ G^*$, we have:
$\rk(\Rep( G^*))=\rk( G^*)$.
Since the category $\Rep( G^*)$ is equivalent to $ P(Gr)$ we find that:
$\rk( G^*)=\rk(\Rep( G^*))=\rk( P(Gr))=\rk(G)$.
\end{proof}

\subsection{}
We introduce now the following subsets in $\XTu$:
\begin{description}
\item[$ ^G\!{R_{\phantom{+}}}${ }={ }]
the root system of the pair $( G,  T)$;
\item[$ ^G\!{R_+}${ }={ }]
the set of positive roots corresponding to a choice of simple roots;
\end{description}
To this set of data one associates the dual data:

\begin{description}
\item[$ ^G\!{R^\vee}${ }={ }]
the root system dual to $ ^G\!R$,
which is viewed as the set in $\XT$ of coroots of $ G$;
\item[$ ^G\!{R_+^\vee}${ }={ }]
the set of positive roots in
$ ^G\!{R^\vee}$ corresponding to the subset $^G\!{R_+}$ in~$^G\!R$.
\end{description}
For $\lambda, \mu\in\XT$, we write: $\mu\le\lambda$ iff $\lambda-\mu$
 \ = \ a sum of positive coroots of $ ^G\!{R^\vee}$. Also, set:
$$ ^G\!{D^\vee}\,=\,\{\lambda\in\XT\,|\,\langle\lambda, x\rangle\ge0\quad\forall x\in  
^G\!{R_+}\}\;,\;\,\mbox{
the dominant Weyl chamber.}$$
Let $\gotht$ denote the Lie algebra of the torus $ T$ and
$\gotht^*$ the dual space.
We shall often view
$\XTu$, $ ^G\!R$, $ ^G\!{R_+}$
as subsets of
$\gotht^*$ and
$\XT$, $ ^G\!{R^\vee}$, $ ^G\!{R_+^\vee}$
as subsets of $\gotht$. Let:
\begin{description}
\item[$h${ }={ }]
the sum of positive roots in $ ^G\!R$.
\end{description}
Thus, $h$ is an element of $\gotht^*$ having the following properties:
\begin{equation}
 \begin{array}{ll}
  \mbox{\rm (i)}&\ h\ \mbox{is a regular element in}\ \gotht^*;\\
  \mbox{\rm (ii)}&\ \langle\lambda,h\rangle\ \mbox{\rm is an integer for any}\ 
\lambda\in\XT;\\
  \mbox{\rm (iii)}&\ \langle\mu,h\rangle<\langle\lambda,h\rangle
  \ \mbox{if}\ \mu\le\lambda,\mu\ne\lambda,\quad\mu,\lambda\in\XT.
 \end{array}
\end{equation}

We have a canonical identification of $\gotht^\vee$, the Lie algebra
of the dual torus $ T^\vee$, with $\gotht^*$. Further, we identify the torus
$ T^\vee$ with its image in $ G^*$. Thus, $ T^\vee$ is a maximal torus in $
G^*$,
$\gotht^\vee$ is a Cartan subalgebra in $\Lie( G^*)$,
and $h$ is an element of $\gotht^\vee$.

\begin{lemma}
$h$ is a regular element in $\Lie( G^*)$.
\end{lemma}

\begin{remark}
The element $h$ is regular, of course, as an element
of $\gotht^*$ (see 3.9.1 (i)), but this
does not necessarily imply that it is regular as an element of $\gotht^\vee$
since we know nothing about
the root system of the pair $( G^*,  T^\vee)$, which we are going to
 introduce now.
\endpr
\end{remark}

Let
$ ^{G^*}\!\!R\subset\XTup$ be the root system of the pair
$( G^*,  T^\vee)$.
We view the set $ ^{G^*}\!\!R$ as a subset of $\XT$ via the identification:
$\XTup\cong\XT$.
Assuming that Lemma 3.9.2 is true and taking (3.9.1 (ii)) into account,
we define the system of positive roots in $ ^{G^*}\!\!R$ by
\begin{equation}
 ^{G^*}\!{R_+}=\{\alpha\in ^{G^*}\!\!R\,|\,\langle\alpha,h\rangle>0\}.
\end{equation}
The choice of positive roots being made, we can define the set
$ ^{G^*}\!{R_+^\vee}$ of positive roots in $ ^{G^*}\!{R^\vee}$
(= \ the root system dual to $ ^{G^*}\!\!R$) and, hence,
the dominant Weyl chamber $ ^{G^*}\!D\subset\XTup$.
Note that the sets $ ^G\!{D^\vee}$ and $ ^{G^*}\!D$ are contained in the same 
lattice $\XT$.

\begin{prop}
Let $\lambda\in {^G\!D^\vee}$ and $V_\lambda$ the irreducible
representation of $ G^*$ corresponding to the complex $ IC_\lambda$
via the equivalence of the categories $\Rep( G^*)$ and $ P(Gr)$.
Then $\lambda\in {^{G^*}\!\!D}$ and $V_\lambda$ is representation with
highest weight $\lambda$.
\end{prop}

\begin{rem}
The Proposition ensures consistency of the notation $V_\lambda$
which a priori means either the representation of $ G^*$
with  highest weight $\lambda$, or the simple object of $\Rep( G^*)$
corresponding to $ IC_\lambda$.
\end{rem}

{\bf Notation.}
Given a representation $V\in\Rep( G^*)$, we let $V(\nu)$
denote the weight subspace of $V$ corresponding to a weight $\nu\in\XTup$
and set:
\[
\Spec(V)\,=\,\{\nu\in\XTup\,|\,V(\nu)\ne0\}.
\]
Similarly, for a perverse sheaf $M\in P(Gr)$ we set:
\[
\Spec(M)\,=\,\{\nu\in\XT\,|\,i^!_\nu M\ne0\}.
\]

\begin{proof}
{Proof of Lemma 3.9.2 and Proposition 3.9.5}
We know that $\langle\lambda,h\rangle\in\Z$ for any
$\lambda\in\XTup$ (by 3.9.1 (ii)).
Hence, forgetting (3.9.4) for a moment,
we can (and will) make a choice of the set
$ ^{G^*}\!\!{R_+}$
of positive roots in such a way that $\langle\alpha, h\rangle\ge0$
for any $\alpha\in ^{G^*}\!\!{R_+}$.
Let $ ^{G^*}\!\!D$ stand for the corresponding dominant Weyl chamber
(we'll show later that these $ ^{G^*}\!\!{R_+}$ and $ ^{G^*}\!\!D$
coincide with those defined by (3.9.4)).

We have the following two key observations:
\begin{equation}
\vbox{
\parshape=1 0pt 250pt
\noindent
For an irreducible representation
$V\in\Rep( G^*)$,
the function
$\Spec(V)\ni\nu\longmapsto\langle\nu,h\rangle$
attains its maximum at the highest
weight of $V$.
}
\end{equation}
And similarly:
\begin{equation}
\hphantom{x}\qquad\vbox{
\parshape=1 0pt 250pt
\noindent
(a) For $\lambda\in ^G\!{D^\vee}$,
the function $\Spec(IC_\lambda)\ni\nu\mapsto\langle\nu,h\rangle$
attains its maximum at the point
$\nu=\lambda$; and\\

(b) $\langle\nu,h\rangle<\langle\lambda,h\rangle$
for any $\nu\ne\lambda$; $\nu\in\Spec( IC_\lambda)$.}
\end{equation}
Claim (3.9.7) follows immediately from (3.9.2 (iii)) and
the implication:
$\nu\in\overline O_\lambda\ \Rightarrow\ \nu\le\lambda$.

Now, let
$\lambda\in {^G\!D^\vee}$,
let $V_\lambda$ be the irreducible representation of $ G^*$
corresponding to the complex $ IC_\lambda$ via
the equivalence of the categories $ P(Gr)$ and $\Rep( G^*)$,
and let $\mu\in ^{G^*}\!\!D$ be the highest weight of $V_\lambda$.
We know that the fixed-point gradation on $H_t( IC_\lambda)$
corresponds to the weight gradation on $V_\lambda$.
Hence, comparison of (3.9.6) with (3.9.7(a))
shows that
$\langle\mu,h\rangle=\langle\lambda,h\rangle$.
Then, the statement (3.9.7(b)) yields
an equality: $\mu=\lambda$. This completes proof of Proposition 3.9.5.

Finally, suppose that $h$ is not a regular element of 
$\gotht^\vee$. 
Then, one can find an irreducible representation $V_\lambda$
(with highest weight
$\lambda\in ^{G^*}\!\!D$)
and a simple root $\alpha\in ^{G^*}\!\!{R_+}$ such that
\[
\langle\alpha,h\rangle=0\quad\mbox{\rm 
and}\quad\lambda-\alpha\in\Spec(V_\lambda).
\]
Translating this into topological language, and setting $\nu=\lambda-\alpha$,
we get
\[
\langle\nu,h\rangle=\langle\lambda,h\rangle\quad\mbox{\rm 
and}\quad\nu\in\Spec( IC_\lambda).
\]
Again, this contradicts to (3.9.7(b)), completing the proof of Lemma 3.9.2.
\end{proof}

\begin{corol}
We have: $ {^G\!D^\vee}= ^{G^*}\!\!D$.
\end{corol}

\begin{proof}{Proof}
It follows from Proposition 3.9.5 that 
$ {^G\!D^\vee}\subset ^{G^*}\!\!D$.
Let $\lambda\in ^{G^*}\!\!D$
and let $ IC_\mu$, $\mu\in {^G\!D^\vee}$ be the object of $ P(Gr)$
corresponding to the irreducible
representation $V_\lambda$. Then, $\lambda=\mu$ by the Proposition,
hence $\lambda\in {^G\!D^\vee}$.
\end{proof}

\subsection{}
We remind the reader that, given a semisimple group $G$ with a maximal torus $ 
T$,
one associates to the pair $( G,  T)$ three natural lattices in $\gotht=\Lie( 
T)$
and three natural lattices in $\gotht^*$:
\[
 ^G\!Q\subset\XTu\subset ^G\!\!\!P\subset\gotht^*\quad\mbox{and}\quad
 ^G\!{Q^\vee}\subset\XT\subset {^G\!P^{\vee}}\subset\gotht.
\]
Here $^G\!Q$ is the lattice generated by the roots of $( G, T)$
(the root lattice),
$ ^G\!{Q^\vee}$ is the lattice generated by coroots (the coroot lattice),
$^G\!\!P$ is the lattice dual to $ ^G\!{Q^\vee}$ (the weight lattice),
and ${^G\!P^\vee}$ is the lattice dual to $ ^G\!Q$
(the coweight lattice).
Note that replacing $ G$ by $ G^\vee$
(=\  the dual group) exchanges the roles of $\gotht$ and $\gotht^*$
and reverses the order
of lattices. Further, the center of a semisimple group
can be expressed in terms of the lattices as follows:
\begin{equation}
Z( G)\cong {^G\!P^\vee}/\XT\quad\mbox{and}\quad
Z( G^\vee)\cong ^G\!\!P/\XTu.
\end{equation}

\begin{proof}
{Proof of Theorem 3.8.1}
We shall view the groups $ G^*$ and $ G^\vee$ as having the {\em same\/}
vector space $\gotht^\vee$ (=\ $\gotht^*$) as a Cartan subalgebra. This space
contains various lattices associated to $ G$, $ G^*$ and $ G^\vee$.

In n.~3.8 we have fixed the set $ ^G\!{R_+}$ of positive roots of the pair
$( G,  T)$ and showed that it determines the choice of positive roots
for the pair $( G^*,  T)$.
On the other hand, the choice of positive roots in $ ^G\!R$
leads to the choice of positive roots in the dual root system. Hence, we may 
speak
freely about positive roots, simple roots, dominant Weyl chambers, etc.~for 
all
three groups: $ G$, $ G^*$, $ G^\vee $.
We claim that the following sets in $(\gotht^\vee)^*=\gotht$ coincide:
\begin{equation}
\mbox{\it simple roots of}\  G^*\ =\
\mbox{\it simple roots of}\  G^\vee.
\end{equation}

To prove this claim we let $ ^G\!{Q_+}$ denote the semigroup in the root lattice
$ ^G\!Q$ generated by positive roots and let
$ ^G\!{Q^\vee_+}$, $ ^{G^*}\!\!{Q_+}$, $\ldots$
denote the objects defined in a similar way.
We can characterize these objects in terms of the categories
$ P(Gr)$ and $\Rep( G^*)$ and the corresponding dominant Weyl chambers.
as follows:
\begin{equation}
 ^G\!{Q^\vee_+}=\{\alpha\in\XT\,|\,\exists\lambda\in {^G\!\!D^\vee}\,:\,
\lambda-\alpha\in\Spec( IC_\lambda)\}
\end{equation}
and
$$
 ^{G^*}\!\!{Q_+}=\{\alpha\in\XTup\,|\,\exists\lambda\in {^{G^*}\!\!D}\,:\,
\lambda-\alpha\in\Spec(V_\lambda)\}.
\eqno{(3.10.3^*)}
$$
\ignorespaces
We know that the sets $\Spec( IC_\lambda)$ and
$\Spec(V_\lambda)$ correspond to each
other via the equivalence of the categories $ P(Gr)$ and $\Rep(G^*)$,
and we have shown that
$ {^G\!D^\vee}= {^{G^*}\!\!D}$ (Corollary 3.9.8).
Hence, the right-hand sides of (3.10.3) and (3.10.3${}^*$) coincide,
and we obtain that
\begin{equation}
 ^{G^*}\!\!{Q_+}\,=\, ^G\!{Q^\vee_+}.
\end{equation}

We  now derive some consequences of (3.10.4).
First, recall that each side in (3.10.4)
is a free semigroup generated by the corresponding set of simple roots.
Therefore, simple roots can be characterized as indecomposable
elements of the semigroup in question
(an element $\lambda$ is said to be indecomposable
if $\lambda=\mu+\nu\,\Longrightarrow\,\mu=0
\ \mbox{\rm or}\ \nu=0$). Hence, (3.10.4) yields:
{\it simple roots of $ G^*$\ =\ simple roots of $ G$},
and (3.10.2) follows.

Next, note that each of the lattices $Q$ is generated by the
corresponding semigroup $Q_+\subset Q$, so that the lattice $Q$
is completely determined by its positive part $Q_+$. Passing
to the dual lattices and using (3.10.4) again, we obtain:
\[
 ^{G^*}\!\!{P^\vee}= ^G\!\!P.
\]
Hence, (3.10.1) implies that
\begin{equation}
Z( G^*)\cong Z( G^\vee).
\end{equation}

Now, we will (temporarily) make an additional assumption,
that the group $ G$ is of {\em adjoint} type, i.e.~has no center.
Then, the group $ G^\vee$ is simply-connected and we claim that
(compare with (3.10.2)):
\begin{equation}
\mbox{\it simple coroots of}\  G^*\ =\ 
\mbox{\it simple coroots of}\  G^\vee.
\end{equation}
To prove the claim note first that
$\XTup= ^{G^\vee}\!\!P$
(since $G^\vee$ is simply-connected)
and, hence, any Weyl chamber in $\XTup$ (of the group $ G^\vee$)
is the free semigroup generated by the corresponding fundamental weights. 
Therefore,
for any Weyl chambers $ G^\vee_{D_1}$ and $ G^\vee_{D_2}$
separated by a wall we have:
\begin{equation}
\begin{array}{l}
\mbox{\it The wall between}\  G^\vee_{D_1}\ \mbox{\it and}\  G^\vee_{D_2}\
\mbox{\it is the hyperplane in}\ (\gotht^\vee)^*=\gotht\\
\noalign{\vskip2pt}
\mbox{\it spanned by}\
 G^\vee_{D_1}\cap G^\vee_{D_2}.
\end{array}
\end{equation}
Furthermore, the line $L=\RRR\cdot\lambda$
generated by a fundamental weight $\lambda$
has the following characterization:
\begin{equation}
L=\ \mbox{\it intersection of the walls containing}\ L.
\end{equation}

Let us now turn to the group $ G^*$.
Notice that Corollary 3.9.8 holds, of course, for any choice of positive roots 
in
$ ^G\!R$. Hence, the Corollary implies that the groups $ G^*$ and $ G^\vee$
give rise to identical partitions of $\XTup$ into Weyl chambers.
(One should be careful here since we do not yet know that the group
$ G^*$ is simply-connected, so that
$\XTup\ne ^{G^*}\!\!P$ in general. Speaking about Weyl chambers of $ G^*$
in $\XTup$ we have in mind intersections of Weyl chambers in $ ^{G^*}\!\!P$
with $\XTup$.)
Let  $\widetilde D_1$, $\widetilde D_2$ be Weyl chambers in $ ^{G^*}\!\!P$
and let $D_i:=\widetilde D_i\cap\XTup$.
We have: $D_1\cap D_2\subset\widetilde D_1\cap\widetilde D_2$.
Suppose now that $D_1$ and $D_2$ are separated by a wall, as Weyl
chambers of $ G^\vee$. Then, (3.10.7) shows
that the set $D_1\cap D_2$ spans a hyperplane in $(\gotht^\vee)^*$.
Hence, this hyperplane must be the wall between 
$\widetilde D_1$ and $\widetilde D_2$.
It follows, that the groups $ G^*$ and $ G^\vee$ give rise
to the same collection of walls in
$(\gotht^\vee)^*$. Thus, fundamental weights of $ G^*$ are proportional
to fundamental weights of $ G^\vee$, by (3.10.8).

We fix a choice of positive roots, let
$\alpha_1,\ldots,\alpha_r$ be the corresponding simple roots
of $ G^\vee$ and $ G^*$ (they are the same by (3.10.2)),
and let 
$\alpha_1^\vee,\ldots,\alpha_r^\vee$
(resp.~$\alpha_1^*,\ldots,\alpha_r^*$)
be the corresponding simple coroots of $ G^\vee$ (resp.~$ G^*$).
It follows from the preceding paragraph
(for, simple roots are dual to the fundamental weights)
that the collection of lines generated by the
$\{\alpha_i^*\}$
is the same as that generated by the
$\{\alpha_i^\vee\}$. Whence,
$\alpha_i^*=\alpha_i^\vee$, for the coroots
$\alpha_i^*$ and $\alpha_i^\vee$ are uniquely determined
by the following conditions:
\[
\langle\alpha_i^*,\alpha_i\rangle=2=\langle\alpha_i^\vee,\alpha_i\rangle
\quad\mbox{and}\quad\langle\alpha_i^*,\alpha_j\rangle<0
\ \mbox{for}\ i\ne j.
\] 

Now, it follows from (3.10.2) and (3.10.6), that for an adjoint
group $ G$, the groups $ G^\vee$ and $ G^*$ have the same
Cartan matrices and therefore have isomorphic Lie algebras.
Hence, there exists a surjective homomorphism
$ G^\vee\rightarrow G^*$ with finite kernel. This kernel must
be trivial, by (3.10.5).

Finally, let us drop the assumption that the group $G$ has no
center and let $G'$ denote the derived group of $G$.
The projection $G\to G'$ gives rise to an inclusion of
Grassmannians: $\Gr_G\hookrightarrow\Gr_G$, and hence to a fully
faithful imbedding of tensor categories:
$P(Gr)\hookrightarrow P(G')$
which may be viewed as an imbedding
$\Rep(G^*)\hookrightarrow\Rep((G')^*)$.
The last imbedding is induced (Proposition 7.1.1) by a homomorphism
$f:(G')^*\to G^*$
which, being restricted to the maximal tori, gives a surjective
homomorphism
$(T')^\lor\to T^\lor$
with finite kernel. Hence, $f$ is also a surjective homomorphism with
finite kernel. It remains to note that the Theorem is already
proved for the group $G'$, and to apply (3.10.5).
\end{proof}

\medskip\noindent {\bf Remark 3.11.}
The proof of  Theorem 3.8.1 shows that the fixed-point
decomposition \ref{fpd} corresponds to the weight decomposition
for representations of the group $G$ with respect to the maximal
torus $T$ embedded into $G$ via homomorphism (3.8.3).
\medskip

\section{Loop group cohomology and the Principal nilpotent}

In this chapter we will answer Question 1 of n.~1.6. (Theorem 1.6.3)
and will approach  Question 2 of n.~1.6.

\subsection{}
For representation
$V\in\Rep(G^\lor)$
let
$V=\oplusl_iV^{{\rm top}}(i)$, $i\in\Z$
be the gradation on the underlying vector space of the representation,
arising --- via Theorem 1.4.1 --- from the natural gradation
on hyper-cohomology. It will be referred to as the ``topological''
gradation.

\begin{prop}
There is a unique semisimple element
$h\in\g^\lor$,
such that, for any
$V\in\Rep(G^\lor)$,
the ``topological'' gradation on $V$ coincides with gradation
$V^h(\cdot)$ by the eigenvalues of $h$
(cf.~(1.6.4)).
\end{prop}

\begin{proof}{Proof}
The Kunneth formula shows (cf.~n.~2.4) that the ``topological''
gradation is compatible with tensor product. The result now
follows from Corollary~7.1.3 applied to $T=\C^*$.
\end{proof}

\begin{lemma}
The element $h$ is a regular element of $\g^\lor$.
\end{lemma}

\begin{proof}{Proof}
of the lemma is similar to that of Lemma 3.9.2\vvv. Let $\t^\lor$
be a Cartan subalgebra of $\g^\lor$ containing $h$. First, we choose
the set $S$ of simple roots of the pair
$(\g^\lor,\t^\lor)$
in such a way that
$\langle\alpha,h\rangle\ge0$
for
$\alpha\in S$.
Then, we have assertion (3.9.6). It holds for the element $h$
we are considering now (although it is different from the one considered
in Lemma 3.9.2\vvv. Furthermore, an argument in the proof of Lemma 3.9.2
shows that if $h$ is not regular, then one can find an irreducible
representation $V_\lambda$ such that
$\dim V^h_\lambda(m)>1$
where $m$ is the maximum of the function
$\Spec(V_\lambda)\ni\nu\mapsto\langle\nu,h\rangle$.
That means, in topological terms, that the dimension of the top
Intersection cohomology group
$IH^m(\overline O_\lambda)$
is greater than $1$. This is impossible, of course,
so that $h$ must be regular.
\end{proof}

\begin{lemma}
If $\lambda$ belongs to the coroot lattice, then all
the eigenvalues of the element $h$ in the representation
$V_\lambda$ are even.
\end{lemma}

Proving the Lemma amounts to showing that
$IH^{{\rm odd}}(\overline O_\lambda)=0$,
provided $\lambda$ satisfies the assumptions of the Lemma. To prove
this cohomology vanishing we recall the following results:
\begin{enumerate}
\renewcommand{\theenumi}{(\alph{enumi})}
\item For any $\lambda$, the cohomology sheaves
${\cal H}^i(IC_\lambda)$
vanish if
$i\not\equiv\dim(O_\lambda)\mod2$
(see~\cite{KL});
\item For any $\lambda$, we have:
$H^{{\rm odd}}_c(O_\lambda)=0$
for $O_\lambda$ has a stratification by complex affine spaces,
the so-called Bruhat decomposition (see~\cite{AP});
\item If $\lambda$ belongs to the coroot lattice, then the stratum
$O_\lambda$ has even dimension (see e.g. \cite{Lu}).
\end{enumerate}

The Lemma now follows from the standard spectral sequence:
\[
\oplusl_{\mu\le\lambda}{\cal H}_\mu^i(IC_\lambda)
\otimes H^j_c(O_\mu)\Rightarrow H^{i+j}(IC_\lambda),
\]
associated to the stratification
$\overline O_\lambda=\mathop\cup\limits_{\mu\le\lambda}O_\mu$
(we write: $\mu\le\lambda$ if
$O_\mu\subset\overline O_\lambda$).
\hfill $\Box$

\subsection{}
Let us interrupt the proof of Theorem 1.5.3 and turn to Question 2
from ${\rm n}^\circ$~1.6\vvv. We first analyze the relationship between the
convolution on the category $P(Gr)$ and the action of the algebra
$H^*(\Omega)$ on hyper-cohomology of the complexes $M\in P(Gr)$.

Let $prim$ denote the space of primitive elements of the Hopf algebra
$H^*(\Omega)$. We view $prim$ as an abelian Lie algebra acting
on $H^*(M)$, $M\in P(Gr)$, by nilpotent transformations. Let
$A$ ($\cong A^r$) be the unipotent algebraic group corresponding
to the Lie algebra $prim$. The $prim$-action gives rise to an algebraic
$A$-action on $H^*(M)$ by unipotent transformations. This way we
obtain a functor $R:P(Gr)\to\Rep(A)$.

\begin{lemma} The functor $R$ is an exact tensor functor.
\end{lemma}

\begin{proof}{Proof}
First, recall the general formula (2.5.1). One knows that, for any
$a\in H^*(Y,\C)$,
the natural action of $a$ on
$H^*(Y,f_*J)$
corresponds, via the isomorphism (2.5.1), to the action of
$f^*a\in H^*(X,\C)$ on $H^*(X,J)$.

Now let $X=\Omega\times\Omega$, $Y=\Omega$
and let $m:\Omega\times\Omega\to\Omega$
be the product map. We have an isomorphism (see ${\rm n}^\circ$~2.5):
\begin{equation}
H^*(\Omega,m_*(M\boxtimes N))\cong H^*(\Omega,M)\otimes H^*(\Omega,N)
\end{equation}
The action of $a\in prim$ on the left-hand side of (4.2.2) corresponds,
by the preceding paragraph, to the action of $m^*(a)$ on the right-hand
side. But
$m^*(a)=a\otimes1+1\otimes a$,
for $a$ is a primitive element in $H^*(\Omega)$. That means that
$R$ is a tensor functor and the Lemma follows.
\end{proof}

\subsection{}
We now compose the equivalence
$\Rep(G^\lor)\stackrel{\sim}{\to}P(Gr)$ provided by
theorem 1.4.1 with the functor $R$ of ${\rm n}^\circ$~4.2, to obtain
a tensor functor:
$\Rep(G^\lor)\to\Rep(A)$.
By Proposition 7.1.1, this functor yields a Lie algebra homomorphism:
\begin{equation}
prim\to\g^\lor
\end{equation}
such that the image of $prim$ consists of nilpotent elements in
$\g^\lor$.

\begin{lemma}
The homomorphism (4.3.1) is injective.
\end{lemma}

To prove the Lemma, it suffices to show that for any $a\in prim$
one can find an irreducible representation
$V_\lambda\in\Rep(G^\lor)$
such that $a$ acts on $V_\lambda$ as a non-zero operator. In topological
terms, this amounts to proving the following more general result:

\begin{prop}
For any
$a\in H^k(\Gr)$,
there is an Intersection cohomology complex
$IC_\lambda\in P(Gr)$
such that the element $a$ acts on $H^*(IC_\lambda)$ as a non-zero
operator.
\end{prop}

\begin{proof}{Proof}
Choose a Borel subgroup $B\subset G$. Let $I$ denote the ``Iwahori''
subgroup of $LG$ consisting of those maps
$f:S^1\to G$
which extend holomorphically to the exterior of the unit disc in the
Riemann sphere, and satisfy the condition:
$f(\infty)\in B$.
Any $\lambda\in X_*(T)$ determines a point in the Grassmannian and
we let $C^\lambda$ denote the $I$-orbit of that point in $\Gr$.
All the orbits $C^\lambda$, $\lambda\in X_*(T)$
are disjoint and form a decomposition of $\Gr$ by affine cells of
finite codimension. The Poincar\'e duals of the fundamental cycles
$[\oC]$, $\lambda\in X_*(T)$ from a basis in $H^*(\Gr)$.

Now let $a\in H^k(\Gr)$. Using the basis
$[\oC]$
we can write:
\begin{equation}
a=\sum_\lambda a_\lambda\cdot[\oC],
\quad\lambda\in X_*(T),\quad a_\lambda\in\C
\end{equation}
where the sum ranges over $k$-codimensional cells $C^\lambda$.
Let $a_\mu\ne0$ for some $\mu$. We'll show that $a$ acts non-trivially
on $H^*(IC_\mu)$.

It will be convenient for us, while working with cycles, to speak
in the language of the Intersection {\em homology} rather than
cohomology. So, the action of $H^*(\Gr)$ may be written as a
cap-product map:
\begin{equation}
H^k(\Gr)\times IH_i(\oO)\stackrel{\cap}{\to}
IH_{i-k}(\oO)
\end{equation}

Let  $m=\dim_\RRR(\oO)$. By the main property of the Intersection
homology there is a non-degenerate pairing:
\begin{equation}
\langle\cdot,\cdot\rangle:IH_{m-k}(\oO)\times  IH_k(\oO)\to
H_0(\oO)=\C
\end{equation}
The maps (4.3.5) and (4.3.6) are related by the following identity
in $H_0(\oO)$:
\begin{equation}
a\cap c=\langle a\cap[\oO]\,,\,c\rangle
\end{equation}
where $a\in H^k(\Gr)$, $c\in IH_k(\oO)$.

Formula (4.3.7) shows that if $a\cap[\oO]$ is a non-zero element
in $IH_{m-k}(\oO)$, then one can always find (because of the
non-degeneracy of
$\langle\cdot,\cdot\rangle$)
a class $c\in IH_k(\oO)$ such that $a\cap c\ne0$. Thus, proving the
Proposition amounts to showing that for the element $a$ from (4.3.4)
we have: $a\cap[\oO]\ne0$.

To find $a\cap[\oO]$ we notice that all the cycles $\oC$, occurring
in (4.3.4), have the same codimension. Hence, none of them, except
$\overline C^\mu$, has non-empty intersection with $\oO$ (see~\cite{PS}).
Whence, $\oC\cap[\oO]=0$ for all $\lambda\ne\mu$ and we may therefore
assume that $a=\overline C^\mu$.

Now, the stratum $O_\mu$ is known to be isomorphic to a vector bundle
over $G/P$, a partial flag manifold for the group $G$. Let
$Z=G/P\hookrightarrow O_\mu$
denote the zero-section of that bundle. One can show \cite{AP} that
the cycle $\oC$ has no intersection with the singular locus of the
closure $\oO$, and meets $O_\mu$ transversely in a subvariety
$S\subset C^\lambda$ (no closure !). Moreover, $S$ is contained
in $Z$ and, viewed as a subvariety of
$G/P$, is a Shubert variety, i.e., the closure of a $B$-orbit. We
denote by the same symbol $S$ the corresponding homology class in
$H_\bullet(Z)$.

Let
$i:Z\hookrightarrow O_\mu$
denote the inclusion. Any cycle in $Z$ may be regarded as an intersection
hogomology cycle in $\oO$, for $Z$ is contained in the non-singular
locus of $\oO$. Hence, we get a natural morphism
\begin{equation}
i:H_\bullet(Z)\to IH_\bullet(\oO)
\end{equation}
It is clear from our construction that
$\overline C^\mu\cap[\oO]=i_!(S)$.
Thus, to prove that
$\overline C^\mu\cap[\oO]\ne0$
it suffices to show that  morphism (4.3.8) is injective (for,
obviously, $S$ is non-zero class in $H_\bullet(Z)$.

To prove injectivity of (4.3.8) set $M=IC_\mu$, and note that, in
the cohomology language,  morphism (4.3.8) is nothing but the
natural map:
\begin{equation}
H^*(Z,i^!M)\to H^*(\oO,M)
\end{equation}
Now, the complex $M$ is $\C^*$-equivariant with respect to the
natural $\C^*$-action on the fibres of the vector bundle: $O_\mu\to Z$.
For such a complex one knows (see \cite[Proposition 10.4]{Gi1})
that:
\begin{equation}
H^*(Z,i^!M)\cong H^*_c(\oO,M)
\end{equation}
Furthermore, the morphism (4.3.9) can be identified with the morphism
$q:H^*_c(O_\mu,M)\to H^*(\oO,M)$
induced by the inclusion:
$O_\mu\hookrightarrow\oO$.
To prove injectivity of the latter morphism, write the long exact
sequence of cohomology:
\begin{equation}
\ldots\to H^{k-1}(j^*M)\to H^k_c(O_\mu,M)\stackrel{q}{\to}
H^k(\oO,M)\to\ldots
\end{equation}
where
$j:\oO\setminus O_\mu\hookrightarrow\oO$
denotes the inclusion. It suffices to show that the connecting
homomorphisms in (4.3.11) vanish.

To that end, let us view $M$ as a pure Hodge module (of weight $0$, say)
in the sense of Saito. Then, the weights of
$H^{k-1}(j^*M)$
are $\le k-1$, since the functor $j^*$ decreases the weights, cf.~\cite{BBD}.
Further, the functor $i^!$ increases the weights while the functor
$H^*_c$ decreases the weights \cite{BBD}. Hence,
$H^k_c(O_\mu,M)$
is pure of weight $k$ because of the isomorphism (4.3.10). Hence,
each connecting homomorphism is a map between spaces of different
weights and therefore must be zero.
\end{proof}

Because of the injectivity of the morphism (4.3.1) we can (and will)
identify the Lie algebra $prim$ with its image in $\g^\lor$,
which will be denoted by $\a$.

\subsection{}
Until the end of chapter 5 we assume $\g$, the Lie algebra of the group $G$,
is simple. Then, there is a unique (up to constant factor)
$\Ad G$-invariant polynomial $P$ on $\g$ of degree $2$. This is the
Killing form $P(x)=(x,x)$. Let
$\alpha_1=\alpha(P)$
be the corresponding primitive class in $H^2(\Gr)$ (see \no1.7).
The class $\alpha_1$ is proportional to the first Chern class of
the Determinant (line) bundle on the Grassmannian (see e.g~\cite{PS}).

Let $a_1$ denote the image of $\alpha_1$ under the homomorphism
(4.3.1). Clearly, the action of $a_1$ in any representation
$V\in\Rep(G^\lor)$
shifts the topological gradation
$V^{\rm top}(\cdot)$
by $2$. Iterating this action, one gets a linear operator:
\begin{equation}
(a_1)^k:V^{{\rm top}}(-k)\to V^{{\rm top}}(k)
\end{equation}

\begin{lemma}
The map (4.4.1) is an isomorphism for any $k\in\Z$.
\end{lemma}

\begin{proof}{Proof}
We may assume that $V=V_\lambda$ is an irreducible representation.
The result then follows from the hard Lefschetz theorem \cite{BBD}
for the
Intersection cohomology (applied to the Chern class of the restriction
to $\overline O_\lambda$ of the Determinant bundle on $\Gr$).
\end{proof}

\subsection{Proof of Theorem 1.6.3}
Let $\lambda$ be the maximal root of $\g^\lor$ so that $V_\lambda=\g^\lor$
is the adjoint representation (recall that we assumed $\g^\lor$
to be simple). Let $h$ be the element introduced in Proposition 4.1.1\vvv.
Obviously, we have:
$[\ad h,\ad a_1]=2\ad a_1$.
Hence, $[h,a_1]=a_1$.

The space $\g^\lor$ is graded by {\em even} integers via  the
``topological'' grading (Lemma 4.1.3), and lemma 4.4.2
holds for the operator $(\ad a_1)^k$ acting on $V=\g^\lor$.
It follows easily (an exercise
in linear algebra) that:
\[
\dim\ker(\ad a_1)=\dim V^{top}_\lambda(0)
\]
By definition,
$V_\lambda(0)=Z_{\g^\lor}(h):=\mbox{{\em centralizer\/
{\rm of} $h$ in $\g^\lor$}}$.
Furthermore,
$\dim Z_{\g^\lor}(h)=\rk\g^\lor$
by Lemma 4.1.2\vvv. Hence, we find that:
$\dim Z_{\g^\lor}(a_1)=\rk\g^\lor$.
Thus $a_1$ is a regular nilpotent.

To complete the proof of Theorem 1.6.3 we choose a Lie algebra
homomorphism
$j:\sl_2(\C)\to\g^\lor$
such that
$j\left(\begin{array}{crc} 0&1\\0&0\end{array}\right)=a_1    $.
It is known that the adjoint representation of $\sl_2(\C)$ in
$\g^\lor$, arising from $j$, does not contain the trivial $1$-dimensional
representation as a component. It follows that for
$h_1:=j\left(\begin{array}{crc}1&0\\0&-1\end{array}\right)$
we have:
\begin{equation}
[Z_{\g^\lor}(a_1),h_1]=Z_{\g^\lor}(a_1)
\end{equation}
Let $A$ be the unipotent subgroup of $G^\lor$ corresponding to the
Lie algebra
$Z_{\g^\lor}(a_1)$.
Equation (4.5.1) shows that the affine space
$h_1+Z_{\g^\lor}(a_1)$ is stable under the adjoint $A$-action, and
that $A$ acts transitively on
$h_1+Z_{\g^\lor}(a_1)$.

Now, let
$h\in\g^\lor$
be the element arising from Proposition 4.1.1, which was considered
at the beginning of the proof. We have:
$[h,a_1]=2\cdot a_1=[h_1,a_1]$.
Hence,
$h\in h_1+Z_{\g^\lor}(a_1)$.
By the preceding paragraph, one can find an element $u\in A$ such
that
$h=u\cdot h_1\cdot u^{-1}$.
Conjugating the homomorphism $j$ by the element $u$, we obtain a
new homomorphism
$j':\sl_2(\C)\to\g^\lor$
such that
$j'\left(\begin{array}{crc}1&0\\0&-1\end{array}\right)=h$.
\endpr

\subsection{Proof of Proposition 1.10.4}
Let $A$ be the isotropy group of the principal nilpotent $n\in \g^\lor$
so that
$\Lie A=Z_{\g^\lor}(n)$.
Given a locally-free sheaf ${\cal V}$ on ${\cal N}$, let ${\cal V}_n$ denote its geometric
fibre at $n$. The assignment ${\cal V}\to {\cal V}_n$ gives rise to a functor from
the category of $G^\lor$-equivariant sheaves on ${\cal N}$ to the category
of $A$-modules. This way we obtain a natural morphism
($\a:=\Lie A$):
\begin{equation}
\Hom_{{{\cal C}oh}({\cal N})}({\cal V}_\lambda,{\cal V}_\mu)\to\Hom_\a(V_\lambda,V_\mu)
\end{equation}
We'll now show that (4.6.1) is an isomorphism.

Identify sheaves on ${\cal N}$ with the corresponding $\C[{\cal N}]$-modules of
global sections and note that:
\[
\Hom_{\C[{\cal N}]}(V_\lambda\otimes\C[{\cal N}],V_\mu\otimes\C[{\cal N}]\cong
V_\lambda^*\otimes V_\mu\otimes\C[{\cal N}]
\]
Hence,  morphism (4.6.1) turns into the map (= restriction to
$n$):
\[
(V_\lambda^*\otimes V_\mu\otimes\C[{\cal N}])^{G^\lor}\to
(V_\lambda^*\otimes V_\mu)^\a
\]
Note that the highest weight of any irreducible constituent
$V\subset V_\lambda^*\otimes V_\mu$
belongs to the root lattice of $G^\lor$,
for orbits $O_\lambda$ and $O_\mu$ are assumed to be in the
same connected component of $Gr$. But for any such module $V$,
Kostant proved in \cite{Ko2} that the restriction map:
$(V\otimes\C[{\cal N}])^{G^\lor}\to V^\a$
is an isomorphism.
\endpr

\section{Filtration and the $q$-analogue of the weight multiplicity}

\subsection{}
In this section we extend results of section 4 to the equivariant
setting and study their relationship with the fixed point decomposition
introduced in section 3.

\subsection{}
Fix an arbitrary element $t\in\t$. We can identify the categories
$P(Gr)$ and $\Rep(G^\lor)$ in such a way that the specialized
equivariant cohomology functor $H_t$ on $P(Gr)$ corresponds to the
forgetful functor on $\Rep(G^\lor)$. The canonical filtration on
$H_t$ (see (8.2.1)) gives rise to a certain increasing filtration
$W_\bullet(V)$ on the underlying vector space of any representation
$V\in\Rep(G^\lor)$.

Recall the element $h$ defined in (1.6.2). The following result is
an equivariant analogue of Theorem 1.6.3.

\begin{Th}{\it
The filtration $W_\bullet$ coincides with the filtration by the
eigenvalues of $h$, i.e. the space $W_i(V)$ is spanned by all those
eigenspaces of $h$ whose eigenvalues are $\le i$.
}\end{Th}

First, we prove the following weakening of the theorem:

\begin{lemma}
There exists a semisimple element $h_1$ such that the filtration
$W_\bullet$ coincides with the filtration by the eigenvalues of $h_1$.
\end{lemma}

\begin{proof}{Proof of  Lemma}
Using the equivalence of categories we can (and will) view the functor
$H^*_T$ from \no3.4 as fiber functor on $\Rep(G^\lor)$. By Proposition
7.2.1, there exists a $\C^*$-equivariant principal
$G^\lor$-bundle $P$ on $\t$ such that the functor $H^*_T$ is
isomorphic to the functor $\Gamma_P$. Hence, the functor $H_t$ is
isomorphic to the functor:
$V\mapsto\mbox{the fibre of $\,P\times_{G^\lor}V$}$ at $t$.

Restrict the bundle $P$ to the $1$-dimensional subspace
$\C\cdot t\subset\t$.
We obtain a $\C^*$-equivariant principal $G^\lor$-bundle on the line.
Now, Proposition 7.3.3 completes the proof.
\end{proof}

Let $\goth n$ denote the nilpotent subalgebra of $\g^\lor$ spanned
by the eigenspaces of the operator $\ad h_1$ corresponding to all
negative eigenvalues. Observe that 
the element $h_1$ in Lemma 5.2.2 is only determined
up to a summand in $\goth n$.

\begin{proof}{Proof of Theorem 5.2.1}
The choice of $h_1$ gives an isomorphism:
$V\cong{\rm gr}^WV$,
for each
$V\in\Rep(G^\lor)$.
On the other hand, the associated graded of the functor $H_t$ with
respect to the canonical filtration is isomorphic to $H^*$, the
ordinary cohomology functor (see (8.5.2)). Thus, the natural
gradation on $H^*$ corresponds, by Lemma 5.2.2, to the gradation
on $V$ by the eigenvalues of $h_1$. Now, Theorem 1.6.3 yields
$h_1\equiv h\mod\goth n$. Since $h$ is regular, we see that $\goth n$
is the nilradical of a Borel subalgebra of $\g^\lor$.
Furthermore, there is an element
$u\in\exp{\goth n}$
such that
$h_1=u\cdot h\cdot u^{-1}$.
But, conjugating by $u$ does not affect the filtration by the eigenvalues.
Thus, we may assume that $h_1=h$.
\end{proof}

\subsection{}
We now turn to an equivariant analogue of Theorem 1.7.6, still
assuming that $t$ is an arbitrary element of $\t$.

First, note that $H_t(\Omega)$, the equivariant cohomology of the
Loop group $\Omega$, has a natural structure of a commutative and
co-commutative Hopf algebra. Moreover, one can show by a standard
argument that this algebra is isomorphic to the symmetric algebra
over the space $prim_t$ of its primitive elements. Further,
the construction of primitive elements,
 given in
\no1.7, carries over to the equivariant setup so that to any
$\Ad K$-invariant polynomial $P$ on $\Lie K$ one can associate a
cohomology class $\alpha_T(P)\in H^*_T(\Omega)$. Furthermore,
this class specializes to a primitive class $\alpha_t(P)\in H_t(\Omega)$
and the elements $\alpha_t(P)$ from a basis of $prim_t$ when
$P$ runs over the set of primitive generators of the algebra of
invariant polynomials.

Next, the element $t$ may be viewed as a point of
$(\t^\lor)^*$.
Hence, it determines a closed coadjoint orbit
$O(t)\subset(\g^\lor)^*$.
Let $O$ be the unique {\em regular\/} coadjoint orbit in
$(\g^\lor)^*$
which contains the orbit $O(t)$ in its closure. We fix a point
$x\in O$ and let
$\a_t\subset g^\lor$
be the Lie algebra of isotropy group of $x$ (under the coadjoint
action). To any $\Ad G^\lor$-invariant polynomial $P^\lor$ on
$(\g^\lor)^*$
one can associate the element of $\a_t$ defined by
$a_t(P^\lor):=dP^\lor(x)$ (cf.~\no1.7).

It is convenient to fix a Killing form on $\g$. This form gives rise
to an invariant quadratic polynomial on $\g$, hence, to a cohomology class
$\alpha_T\in H^2_T(\Omega)$, hence to a primitive class $\alpha_t
\in prim_t$. On the other hand, the Killing form gives rise
to the linear function $(t,\cdot)$ on $\t$. We let $t^\lor$ denote
the element of $\t^\lor$ corresponding to this function via the
identification:
$\t^*\cong\t^\lor$.
Further, we pick up a principal nilpotent $n_t$ in
$Z_{\g^\lor}(t^\lor)$
and set
$a_t:=t^\lor+n_t$.
The algebra $\a_t$ introduced above is (up to conjugation) just the
centralizer of $a_t$ in $\g^\lor$. So, we assume that $a_t\in\a_t$.

Recall the bijective correspondence
$P\leftrightarrow P^\lor$
between invariant polynomials on $\g$ and $\g^\lor$, see (1.7.3).

\begin{Th}{\it
The natural action of
$\alpha_t(P)\in prim_t$
on
$H_t(IC_\lambda)$, $\lambda\in X_*(T)$,
corresponds to the action of
$a_t(P^\lor)\in\a_t$
in the representation $V_\lambda$. In particular the class $\alpha_t$
corresponds to the element $a_t$.
}\end{Th}

\begin{corol}
There is a natural Hopf algebra isomorphism: $H_t(\Gr)\cong U(\a_t)$.
\end{corol}

\begin{rem}
 We defined the algebra $\a_{t}$ as the  Lie
algebra of  $ G^{\lor}(x)$,  the isotropy group of an element
$x\in O$. Corollary 5.3.2 says that there is  a  natural
isomorphism between   the  family  $\{\a_{t}\}$  and the family
$\{prim_{t}\}$. 
\end{rem}

     To prove Theorem 5.3.1 we study first the restriction  of
$T$-equivariant cohomology  classes of the Grassmannian to the
lattice $X_{*}(T)\subset\mbox{\Gr}$.   For   any    point
$\lambda\in X_{*}(T)$                 we                 have:
$H_{T}^{*}(\lambda)=\C[\t]$,   for   $\lambda$ is  a
$T$-fixed point. Hence, the space $H_{T}^{*}(X_{*}(T))$ may be
identified with the space  of  $\C$-valued   function   on
$\t\times X_{*}(T)$, polynomial    with    respect    to    the
$\t$-factor.

\begin{prop}
Let $P$ be an invariant polynomial
on $\g$ and $\alpha_{T}(P)\in  H_{T}^{*}(\Gr)$,  the
corresponding cohomology class.  Then, the restriction of that
class to $X_{*}(T)$  is given by the following function
on $\t\times X_{*}(T)$:
\[
(t,\lambda)\longmapsto\langle dP(t),
\lambda\rangle
\]
Here $dP(t)$ stands for the differential  at  the
point t of the restriction of the polynomial $P$ to $\t$,
and $\lambda$ is viewed as a vector in $\t$,
so that the function{\rm:} $\lambda\mapsto\langle dP(t),
\lambda\rangle$ is  a  linear  function.
\end{prop}

\begin{proof}{Proof} Cohomology  classes  $\alpha_{T}(P)$  were
constructed in  n.~1.7  by  means  of   a   certain   principal
$G$-bundle on  $S^{2}\times\Omega$.  Let  $\widetilde  L$  be  its
restriction to $S^{2}\times X_{*}(T)$.  It  follows  from  the
construction that the structure group of the bundle $\widetilde L$
can be reduced from $G$ to $T$. Hence, the restrictions of the
classes $\alpha_{T}(P)$  to $X_{*}(T)$ come from a $T$-bundle $L$
on $S^{2}\times X_{*}(T)$ such that $\widetilde L  =  G\times_{T}L$.
Thus,   we   can   forget  about  the  group  $G$  and
concentrate our  attention  on  the   $T$-bundle   $L$ which
corresponds, algebraically,  to  restricting  polynomials $P$ to
the Cartan subalgebra $\t$.

     Now, pick   up   some   $\lambda\in   X_{*}(T)$  and  let
$L_{\lambda}$ be  the  restriction  of  the  bundle   $L$   to
$S^2\times\{\lambda\}$. Using  a universal $T$-bundle:  $ET\to BT$,  we form
the diagram:
\begin{equation}
BT\stackrel{\pi}{\gets}ET\times_{T}L_{\lambda}
\stackrel{p}{\to}T\setminus L_{\lambda}\cong S^{2}
\end{equation}
where the  projection  $\pi$  is  a   fibration   with   fiber
$L_{\lambda}$ and the projection $p$ is a fibration with fiber
$ET$. There   is   a   natural   $T$-action   on   the   space
$ET\times_{T}L_{\lambda}$ arising   from   the  $T$-action  on
$L_{\lambda}$ (we used here that $T$  is  an  abelian  group).
Coupling the latter $T$-action with the second copy of $ET$ we
obtain from (5.3.4) a diagram:
\begin{equation}
BT\stackrel{\widetilde  \pi}{\gets}ET\times_{T}
(ET\times_{T}L_{\lambda})\stackrel{\widetilde p}{\to}BT\times
(T\setminus L_{\lambda})\cong BT\times S^{2}
\end{equation}
The map  $\widetilde \pi$  here  forgets  about  the  new  factor
$ET$, while the map $\widetilde p$ projects this  new  factor  $ET$
onto $BT$.  The fiber of the projection $\widetilde p$ remains the
same as that  of  the  projection  $p$  in  (5.3.4),  i.e.  is
isomorphic to $ET$.  The space $ET$ being contractible, we see
that $\widetilde p$ is  a  homotopy  equivalence.  Hence,  diagram
(5.3.5) gives rise to a morphism of cohomology:
\begin{equation}
H^{*}(BT)\stackrel{\widetilde \pi^{*}}
{\longrightarrow}H^{*}(ET\times_{T}(ET\times_{T}L_{\lambda}))
\simeq H^{*}(BT\times S^{2})
\end{equation}
Upon substituting   $H_{*}(BT)=\C[\t]$ we  finally
obtain a map
$c_{\lambda}:\C[\t]    \to\C[\t]\otimes
H^{*}(S^{2})$. One  sees from definitions that the restriction
of the class $\alpha_{T}(P)$ to $\t\times\{\lambda\}$ equals the
integral of  the  element  $c_{\lambda}(P) \in\C[\t]
\otimes H^{*}(S^{2})$ over the $2$-sphere.

     To complete the proof of the Proposition,  it suffices to
prove the following formula for the map $c_{\lambda}$:
\begin{equation}
c_{\lambda}(P)=P\otimes 1    +    \langle
dP,\lambda\rangle \otimes   u,\qquad   P\in\C[\t]
\end{equation}
where $u$ denotes a generator of $H^{2}(S^{2},\Z)$.
To prove             this            formula            write:
$c_{\lambda}(P)=c_{\lambda}'(P)\otimes 1+     c_{\lambda}''(P)
\otimes u$.  We first compute $c_{\lambda}'$. To that end, pick
up a point $pt$ in $S^{2}$ and compose all the  morphisms  (5.3.6)
with the   restriction   morphism:  $H^{*}(BT\times S^{2}) \to
H^{*}(BT\times pt)=  H^{*}(BT)$.  The  resulting  morphism  is
induced by the following diagram of maps:
\[
BT\stackrel{p_{1}}{\longleftarrow}(ET\times
ET)/T\stackrel{p_{2}}{\longrightarrow}BT
\]
where $p_{1}$, $p_{2}$  denote  the  first   and   the   second
projections. One shows easily  that  this  diagram induces the
identity morphism on cohomology. Hence, $c_{\lambda}'(P)=P$.

     To proceed   further,  we  may  assume  without  loss  of
generality that  $T=S^{1}$, so that $H^{*}(BT)=\C[v]$
where $v$  is a generator of $H^{2}(BT,\Z)$.  Note
next that  both  maps  (5.3.6)   and   (5.3.7)   are
{\em algebra}
homomorphisms. Hence,  we  have  only to check that these maps
coincide on the generator $v$.  But  $c_{\lambda}''(v)$  is  a
scalar. To  find  it,  we  pick  up  a  point $pt$ in $BT$ and
compose the morphism (5.3.6) with  the  restriction  morphism:
$H^{*}(BT \times   S^{2})   \to   H^{*}(pt   \times  S^{2})  =
H^*(S^{2})$. The resulting morphism is induced by the diagram
(5.3.5). Now,  one  verifies easily that this diagram sends the
generator $v \in H^{2}(BT)$ to $m \cdot u$ where  the  integer
$m$ is the degree of the homomorphism $\lambda:S^{1} \to T =
S^{1}$. Thus  $c_{\lambda}''(v)  =  m$,  in   agreement   with
(5.3.7).
\end{proof}

\begin{proof}{Proof of theorem 5.3.1}
The  argument  of
n.~4.2 shows  that  the  action  of  the  abelian  Lie  algebra
$prim_{t}$ on the equivariant cohomology of complexes  $M  \in
P(Gr)$ gives  rise to a tensor functor from the category $P(Gr)$
to the category of $prim_{t}$-modules. That gives, via Theorem
1.4.1, a  tensor  functor on the category $\Rep(G^{\lor})$.
Proposition 7.1.1 now shows that,  for each $a  \in  prim_{t}$,
the endomorphism $\exp a$ on the underlying vector space of any
representation is  induced  by  the  action  of   a   uniquely
determined element     of    the    group    $G^{\lor}$.    By
differentiation, one gets a Lie algebra homomorphism:
\begin{equation}
\varphi_{t}:prim_{t} \longrightarrow
\g^{\lor}
\end{equation}
     The family  $\{\varphi_{t}  \}$  algebraically depends on
the parameter $t \in \t$.

     Recall the canonical filtration on $H_{t}(\Gr)$
and for any $a \in prim_{t}$ let $\overline a \in {\rm gr}^{W}_{\bullet}
H_{t}(\Gr)$ denote the ``principal symbol'' of $a$,
viewed as an ordinary cohomology class, cf.~(8.5.2).  We  can
find a  complex  $IC_{\lambda} \in P(Gr)$ such that the element
$\overline a$ acts on $H^{*}(IC_{\lambda})$ as a non-zero  operator
(Proposition 4.3.2). Then, $a$ acts on $H_{t}(IC_\lambda)$ as a
non-zero operator,   for    $H^{*}(IC_{\lambda})=    {\rm
gr}^{W}_{\bullet}H_{t}(IC_{\lambda})$. Hence,  morphism (5.3.8) is
injective.

     We will  show  later  (in  the  course  of  the  proof of
Proposition 5.4.1) that $\varphi_{t}(\alpha_{t})$, the image of
the cohomology  class  corresponding  to  the  Killing form on
$\g$, is a regular element  of $ \g^{\lor}$.
Hence, $Z_{\g^{\lor}}(\varphi_{t}(\alpha_{t}))$, the centralizer
of $\varphi_{t}(\alpha_{t})$  in  $\g^{\lor}$,   has
dimension =   $\rk\g^{\lor}$.   We   find  that  $\dim
Z_{\g^{\lor}}(\varphi_{t}(\alpha_{t})) =\dim  (prim_{t})  =
\dim\varphi_{t}(prim_{t})$. But   the  image  of  morphism
(5.3.8), being an abelian subalgebra of $\g^{\lor}$,
is contained in $Z_{\g^{\lor}}(\varphi_{t}(\alpha_{t}))$. Thus,
$\varphi_{t}(prim_{t})=Z_{\g^{\lor}}(\varphi_{t}(\alpha_{t}))$.

     Assume first that $t$ is a regular element in $\t$. Then,
for each  complex  $M\in  P(Gr)$,  we have  the   fixed-point
decomposition \ref{fpd}  on $H_{t}(M)$.  Any cohomology class in
$H_{t}(\Gr)$  acts  on  $H_{t}(M)$  as  a   diagonal
operator with  respect  to this decomposition.  Moreover,  the
class $\alpha_{t}(P)$ associated to  an  invariant  polynomial
$P$ acts on $H_{t}(M)$ as the semisimple element $dP^{\lor}(t)
\in \t^{\lor}$ (by Proposition 5.3.3). This proves the Theorem
for all   regular   $t\in   \t$.  For  an  arbitrary  $t$,  by
continuity, be     have:     $\varphi_{t}(\alpha_{t}(P))     =
dP^{\lor}(x)$ where  $x\in  (\g^{\lor})^{*}$  is an
element whose semisimple part belongs to the orbit $O(t)$.
But if      $P$      is     the     Killing     form,     then
$\varphi_{t}(\alpha_{t}(P))$ was  shown  to   be   a   regular
element. Hence,  the point $x$ must be regular and the Theorem
follows.
\end{proof}

\begin{rem}
For $t=0$ Theorem 5.3.1 yields Theorem 1.7.2.
\end{rem}

\subsection{}
We shall now analyze the relative  position (cf. remark 1.6.4)  of
the fixed  point  decomposition  \ref{fpd}  with  respect to the
canonical filtration  on  equivariant  cohomology  $H_{t}(M)$,
$M\in P(Gr)$. In fact,  we shall do this simultaneously for all
$M$ by finding the relative position inside $\g^{\lor}$ of  the
subalgebra $\a_{t}$ = {\em the  image of the homomorphism\/}
(5.3.8) with respect to the element $h$  arising  from  Theorem
5.2.1 (note  that  the  {\it absolute\/} position of either
$\a_{t}$ or $h$ inside $\g^{\lor}$ makes no  sense  since  these
objects are only determined up to conjugacy).

     Recall that  the  choice  of  the  Killing  form on  $\g$
determines a semisimple conjugacy class $O(t^{\lor})$ in
$\g^{\lor}$ and a cohomology class $\alpha_{t}$ and, hence, its
image $a_{t}\in \a_{t}$.   We   choose   a   representative
$t^{\lor}\in O(t^{\lor})$ which commutes with $h$. Note
also that  conjugating $h$ by an element of $\exp{\goth n}$ (see
$n^{\circ}$~5.2) doesn't affect the filtration $W_{\bullet}$.

The relative position of the element $a_{t}$ (and, hence,
of the subalgebra $\a_{t}$, i.e.,  the  centralizer  of
$a_{t}$) with respect to $h$ is described by the following

\begin{prop}
Conjugating  $h$  by   an
element of the group $\exp{\goth n}$, if necessary, one can get:
$a_{t}=n+t^{\lor}$.
\end{prop}

\begin{proof}{Proof}  On  the  space  $\g^{\lor}$   of   the   adjoint
representation  we  have  the canonical filtration $W_{\bullet}$
such   that   $h\in    W_0(\g^{\lor})$    and    ${\goth
n}=W_{-2}(\g^{\lor})$  (the  filtration  is  indexed  by  even
integers because of Lemma 4.1.3).  We know also that $a_{t}\in
W_{2}(\g^{\lor})$  and  Theorem  1.7.2,  combined  with Theorem
5.2.1,  yields:  $a_{t}\equiv n\mod W_0(\g^{\lor})$.
Hence, $a_{t}\in n + W_0(\g^{\lor})$. We set ${\goth b}=
W_0(\g^{\lor})$.  This is clearly a Borel subalgebra  of
$\g^{\lor}$.

  Let us compare the  elements  $a_{t}$  and  $n+t^{\lor}$.
Both of   them  belong  to  $n+\b$  and,  hence,  are  regular
(see \cite[p.108]{Ko3}). The  semisimple  part  of  $n+t^{\lor}$   is
clearly $G^{\lor}$-conjugate  to  $t^{\lor}$.  The same is true
for $a_{t}$ by Theorem 5.3.1.  Hence, $a_{t}$ and $n+t^{\lor}$
are $G^{\lor}$-conjugate. But, Kostant showed in the course of
the proof   of    \cite[Theorem~1.2]{Ko3} that     any     two
$G^{\lor}$-conjugate elements of $n+{\goth b}$ are conjugate
by an element $u\in \exp{\goth n}$. Hence, we have:  $a_{t}=
u(n+t^{\lor})u^{-1}$. Replacing $h$ by $u\cdot h\cdot u^{-1}$
completes the proof.
\end{proof}

\subsection{}
 Let  us  now  choose  the  particular element
$t\in\t$ so that $t^{\lor}=h$.  Hence,  $a_{t}=n+h$ is regular
and $\a_{t}$ is a Cartan subalgebra of $\g^{\lor}$.

     Given a representation $V\in \Rep(G^{\lor})$ and $\mu  \in
X_{*}(T^{\lor})$, let  $V(\mu)$  denote  the weight subspace of
$V$ of weight $\mu$ with respect to the Cartan subalgebra
$\a_{t}$. Further, let
\[
V_{i}(\mu):=\ker(n^{i+1})\cap V(\mu), \qquad
i=0,1,2,\ldots
\]
     be the filtration on $V(\mu)$ by the kernels of powers  of
the $n$-action  on    $V$.    On    the    other   hand,   let
$W_{\bullet}(V(\mu))$ denote the filtration on $V(\mu)$  induced by
restriction of the $W$-filtration on $V$.

     We have the following elementary result.

\begin{lemma}
The filtrations $W_{\bullet}(V(\mu))$ and
$V_{\bullet}(\mu)$ coincide up to shift; more precisely, one has:
\[
V_{i}(\mu)=W_{2i+\mu(h)}(V(\mu))=
  W_{2i+\mu(h)+1}(V(\mu)).
\]
\end{lemma}

\begin{proof}{Proof}
Set  $u=\exp  n$.  Then,  $a_{t}=h+n=u\cdot
h\cdot u^{-1}$. For $v\in V(\mu)$ and any $i\ge0$ we have:
\[
h\cdot(u^{-1}\cdot n^{i}\cdot   v)=(u^{-1}
\cdot a_{t}\cdot u)\cdot u^{-1}\cdot n^{i}\cdot v = (\mu(h)  +
2i)\cdot (u^{-1}\cdot n^{i}\cdot v),
\]
so that $u^{-1}\cdot  n^{i}\cdot  v\in  W_{2i+\mu(h)}(V(\mu))$.
Now, if $v \in V_{k}(\mu)$, then:
\[
v=u^{-1}\cdot u \cdot v = u^{-1}\cdot v  +
u^{-1}\cdot n\cdot  v  + \ldots +(1/k!)\cdot u^{-1}\cdot n^{k}
\cdot v
\]
     The last  sum  clearly belongs to $W_{2k+\mu(h)}(V(\mu))$,
hence, $V_{k}(\mu)\subset W_{2k+\mu(h)}(V(\mu))$.  The opposite
inclusion is proved in a similar way.
\end{proof}

     Now let  $M={\cal P}(V)\in   P(Gr)$   be   the   perverse   sheaf
corresponding to a representation $V\in\Rep(G^{\lor})$. Recall
that the weight decomposition  of  $V$  with  respect  to  the
Cartan subalgebra  $\a_{t}$  corresponds  to  the  fixed point
decomposition (3.7.1) of the equivariant cohomology $H_{t}(M)$
(via the identification: $V\cong H_{t}({\cal P}(V))$).  Further,  the
canonical filtration  $W_{\bullet}$  on  $H_{t}(M)$  induces   a
filtration (which   we   call   ``canonical''  and   denote  by
$W_{\bullet}$ again) on each direct summand of the decomposition
(3.7.1). Let $M_{\mu}:=H_{t}(i^{!}_{\mu}M))$ denote the summand
corresponding to a point  $\mu\in  X_{*}(T)$.  Theorem  5.2.1,
and Proposition 5.4.1, combined with Lemma 5.5.1, yield the following
result.

\begin{prop}
For  any  $\mu\in X_{*}(T)$,  the
canonical filtration  on  ${\cal P}(V)_{\mu}$ corresponds,  up to the
shift $\mu(h)$, to the  filtration  $V_\bullet(\mu)$ on  the  weight
subspace $V(\mu)$.       \hfill    $\Box$
\end{prop}

\subsection{Proof  of   Theorem   1.11.2}
Let   $V   \in
\Rep(G^{\lor})$ and  $M={\cal P}(V)  \in  P(G^{\lor})$. By Proposition
5.5.2 we have:
\begin{equation}
P_{\mu}(V,q)=\sum_{i}q^{2i}\cdot\dim(W_{2i+
\mu(h)}(M_{\mu})/W_{2i+\mu(h)-2}(M_{\mu}))
\end{equation}

     Now, let  $i_{\mu}:\{\mu\}\hookrightarrow\Gr$ denote the
embedding. We need the following result which will  be  proved
in the next $n^{\circ}$.

\begin{prop}
The  natural  morphism:
$H_{t}(i^{!}_{\mu}M)\to H_{t}(M)$  is strictly compatible with
the canonical filtration, for any $M\in P(Gr)$.
\end{prop}

     The Proposition  implies  that the filtration $W_{\bullet}$
on $M_{\mu}$  corresponds  to  the  canonical  filtration   on
$H_{t}(i^{!}_{\mu}M)$. Next,  note  that  we  have  a canonical
isomorphism:                    $H^{*}_{T}(i^{!}_{\mu}M)\cong
H^{*}(i^{!}_{\mu}M) \otimes \C[\t]$  since  $\mu$  is  a  fixed
point. Evaluation      at       $t\in\t$       yields       an
isomorphism: $H_{t}(i^{!}_{\mu}M) \cong  H^{*}(i^{!}_{\mu}M)$,
and it is clear that the canonical filtration on the left-hand
side corresponds to the filtration by degree on the right-hand
side. Thus, we obtain from (5.6.1), that:
\begin{equation}
P_{\mu}(V,q)=\sum_{i}q^{2i}\cdot\dim
H^{2i+\mu(h)} (i^{!}_{\mu}M)
\end{equation}

     Assume now   that   $V=V_\lambda$   is  an  irreducible
representation. Then,  $M=IC_\lambda$. The   dimensions   of
stalks of    the   cohomology   sheaves   of   the   complexes
$IC_\lambda$ were computed in \cite{KL}.
It was proved there that:
\begin{equation}
\sum_{j\ge0}q^j\cdot\dim
H^j(i^!_\mu IC_\lambda)=              q^{\lambda(h)\cdot}
P_{\mu,\lambda}(q^2),
\end{equation}
where $P_{\mu,\lambda}$ are  the   Kazhdan-Lusztig   polynomials.
Inserting (5.6.4)  into  (5.6.3)  completes  the  proof of the
Theorem. \endpr

\subsection{}
To start proof of  Proposition 5.6.2 we
need some preparations. First, note that the torus $T\subset G$
can be embedded into a bigger torus $\hat T$,  the maximal torus
of the semidirect product $\C^*\ltimes LG$, where $\C^*$ acts on  $LG$
 via ``rotation of the loop''
(see e.g.~\cite{AP})). We have $\hat T= \C^* \times T$. Observe that,
the subgroup $L^+G$ being $\C^*$-stable with respect to the 
rotation of the loop action, 
 the
torus $\hat T$ acts  on $\Gr=LG/L^+G$
in in a natural way. Moreover, for any  point $\mu \in X_{*}(T)$, 
one can find a one-parameter subgroup $\C^* \hookrightarrow \hat T$
that 
contracts a neighbourhood of $\mu$ to $\mu$ as the parameter
$z\in \C^*$ in the subgroup tends to zero.  Given $\mu$,
fix such a subgroup $\C^{*}$. 

Further, let $L^-G$ be the subgroup
of $LG$ formed by the loops $f: \C^* \to G$ that are regular
at $z=\infty$ and such that $f(\infty)=1$. The subgroup 
$L^-G$ should be thought of as `complementary' to $L^+G$,
cf. \S 8.4. In particular, $L^-G$-orbits form an (infinite
dimensional)
 cell
decomposition of $Gr$ which is `transverse', in a sense,
to the stratification by the $L^+G$-orbits. 

Fix a point $\mu \in X_{*}(T)$ and write $L^-G\cdot\mu$ for the
$L^-G$-orbit through $\mu$.
     Let $O_{\lambda}$ be some other $L^{+}G$-orbit in $\Gr$ such that
$\mu \in\overline O_{\lambda}$. 
We write $\{U_i\}$ for the
collection of the intersections of $L^-G\cdot\mu$ with all the $L^+G$-orbits
$O_i$ that contain $\mu$ in its closure and such that $O \subset\overline
O_{\lambda}$. Enumerating the pieces $U_i$ in such a way that
the dimensions of the $U_i$ form a non-decreasing sequence,
we
obtain 
a finite partition: $\overline
O_{\lambda}=\mathop\cup\limits_{i=1}^{n}U_{i}$ with the 
properties (i)-(v) listed below.

\begin{enumerate}
\item Each  $U_{i}$  is  a  $\hat  T$-stable locally-closed
(singular) subvariety of $\overline O_{\lambda}$;
\item $U_{i}$ contains a single $\hat T$-fixed point,  say
$\mu_{i}$, and $\mu_{n}=\mu$;
\item The  subgroup  $\C^{*}  \subset  \hat T$ contracts
$U_{i}$ to $\mu_{i}$;
\item The   subvariety  $Y_{k}:=\mathop\cup\limits_{i\le k}U_{i}$
is closed in $\overline O_{\lambda}$ (for each  $k=1,2,\ldots,n$)  and
$U_{k}$ is a Zariski-open part of $Y_{k}$.
\end{enumerate}

Let us view the complex $IC_{\lambda}$ as a pure Hodge  module
(of weight $0$, say) in the sense of Saito~\cite{Sa}.  In addition to the
properties (i)--(iv) above one has, due to transversality of $L^-G$-
and
$L^+G$-orbits:
\begin{enumerate}
\setcounter{enumi}{4}
\item The restriction  of  the complex $IC_{\lambda}$ to any
$U_{i}$ is pure.
\end{enumerate}

     We have following inclusions:
\[
Y_{k}\stackrel{j_{k}}{\hookrightarrow}\overline
O_{\lambda}\qquad,\qquad Y_{k-1}\stackrel{v}{\hookrightarrow}
Y_{k}\stackrel{u}{\hookleftarrow}U_{k}
\]
Set $L_{k}:=j^{*}_{k}IC_{\lambda}$, $k=1,2,...,n$.

\begin{lemma}
The  equivariant  cohomology:  (a)
$H^{*}_{T}(u_{!}\cdot u^{!}L_{k})$  and  (b) $H^{*}_{T}(L_{k})$
are pure and free $ \C [\t]$-modules.
\end{lemma}

\begin{proof}{Proof} (a)
Let  $i_{k}:\{\mu_k\}\to  Y_{k}$  denote  the
inclusion. We have an isomorphism:
\begin{equation}
H^{*}(u_{!}\cdot u^{!}L_{k})\cong H^{*}
(i^{!}_{k}L_{k})
\end{equation}
which is due to the fact that $L_{k}$ is a $\C^{*}$-equivariant
complex on $U_{k}$ and the $\C^{*}$-action contracts $U_{k}$
to the point $\mu_{k}$. But the functor $H^{*}u_{!}$ decreases
the weights,  the  functor $H^{*}i^{!}_{k}$ increases weights,
and the complex $u^{!}L_{k}$ is pure  by  the  property  (v)
above. Hence, both sides of (5.7.2) are pure.

     Further, the classifying space $BT$ can be represented as
direct limit of pure varieties ($\cong$ product of copies of
$\C{\Bbb P}^{d}$). Hence, the spectral sequence for the equivariant
cohomology (8.1.1)  yields  an equivariant analogue of (5.7.2)
and shows that all the equivariant cohomology  involved  there
are pure. Finally,  $H^{*}_{T}(i^{!}_{k}L_{k})$  is obviously a
free $\C[\t]$-module  since  the  complex  $i^{!}_{k}L_{k}$  is
supported on  a fixed point.  This completes the proof of part
(a) of the Lemma.

     Part (b)  follows from part (a) by induction on $k$ using
the following long  exact sequence of equivariant cohomology:
\begin{equation}
\ldots \to H^{*}_{T}(u_{!}\cdot u^{!}L_{k})
\to H^{*}_{T}(L_{k})\to H^{*}_{T}(L_{k-1})\to\ldots
\end{equation}
\end{proof}

\begin{proof}{Proof of Proposition 5.6.2}
First,  note  that  all
terms in  the  long exact sequence (5.7.3) are pure,  by Lemma
5.7.1. Hence,  all the  connecting  homomorphisms  in  (5.7.3)
vanish so  that  the  long exact sequence breaks up into short
exact sequences. Moreover, these exact sequences are {\it split
\/} as sequences of $\C[\t]$-modules, since $H^{*}_{T}(L_{k-1})$ is a
free $\C[\t]$-module,  by  Lemma  5.7.1(b).  Thus,   the   map:
$H^{*}_{T}(u_{!}\cdot u^{!}L_{k})\to  H^{*}_{T}(L_{k})$  is an
injection onto a direct summand.  Note that this injection can
be identified       with       the      natural      morphism:
$H^{*}_{T}(i^{!}_{k}L_{k})\to H^{*}_{T}(L_{k})$,  due   to   an
equivariant analogue of isomorphism (5.7.2).

     Now, put $k=n$ in  the  above  argument.  Then,  we  have
$\mu_{n}=\mu$ (property (ii)), $i_{n}=i_{\mu}$, $L_{n}=IC_{\lambda}$.
Thus, the $\C[\t]$-module $H^{*}_{T}(i^{!}_{\mu} IC_{\lambda})$
maps injectively onto a direct summand of  $H^{*}_{T}(IC_{\lambda})$.
The  Proposition follows.
\end{proof}

\section{Moduli spaces and Hecke operators}
\def\labelenumi{\thesubsection\theenumi}
The purpose of this chapter is to give some definitions and constructions
that will put informal arguments of n.~1.5 on a solid mathematical basis.
In particular, we give proof of the `well-known'
double-coset construction of the moduli space of $G$-bundles in terms
of
loop groups, and also give meaning to the Poincar\'e duality for the Grassmannian 
$\Gr$,
used in the proof of Proposition 4.3.3.

The constructions below have been clarified during my talks with A.~Bei\-linson.
I am glad to express to him my deep gratitude. I am also indebted to
R.~Bezrukavnikov
for streamlining the argument in proof of theorem 6.3.1.

\subsection{A formal construction}
Assume given the following set of data:
\begin{enumerate}
\item A partially ordered inductive set $D$ (``inductive'' means that,
for any $\alpha$, $\beta\in D$, there exists $\gamma\in D$ such that
$\alpha\leq\gamma$ and $\beta\leq\gamma$;
\item A collection $\{\cal M_{\mu/\lambda}\}$ of {\it smooth} algebraic
varieties $\cal M_{\mu/\lambda}$ indexed by all pairs $(\lambda, \mu)\in
D\times D$ such that $\lambda\leq\mu$;
\item Zariski-open imbedding $j_{\alpha,\beta}:\M_{\gamma/\alpha}
\hookrightarrow\M_{\gamma/\beta}$, for each triple $\alpha$, $\beta$,
$\gamma\in D$ such that $\alpha<\beta\leq\gamma$;
\item A smooth projection $p_{\nu,\mu}:\M_{\nu/\lambda}
\longrightarrow\M_{\mu/\lambda}$ making $\M_{\nu/\lambda}$ an affine
bundle over $\M_{\mu/\lambda}$ (i.e.\ $p_{\nu,\mu}$ is a locally-trivial 
fibration
with affine linear space as fibre), for each triple $\lambda$,
$\mu$, $\nu\in D$ such that $\lambda\leq\mu<\nu$.
\end{enumerate}

The imbeddings $j_{\alpha,\beta}$ and the projections $p_{\nu,\mu}$ should
satisfy the following two properties:

For any $\alpha<\beta<\gamma\leq\mu$, (resp. $\alpha\leq\lambda<\mu<\nu$), the
following triangles commute:

\begin{equation}
\vcenter{\hbox{\diagram
\M_{\mu/\alpha}\rrto^{j_{\alpha,\gamma}}\ddrto_{j_{\alpha,\beta}}&&\M_{\mu/\gamma}
&\M_{\lambda/\alpha}&&\M_{\nu/\lambda}\llto^{p_{\nu,\lambda}}\ddlto^{p_{\nu,
\mu}}\\
\\
&\M_{\mu/\beta}\uurto_{j_{\beta,\gamma}}&&&\M_{\mu/\alpha}\uulto^{p_{\mu,\lambda
}}&\\
\enddiagram}}\label{first}
\end{equation}

For any $\alpha<\beta\leq\lambda<\mu$ the following diagram is a cartesian
square:
\begin{equation}
\vcenter{\hbox{\diagram
\M_{\mu/\alpha}\rrto_{j_{\alpha,\beta}}\ddto_{p_{\mu,\lambda}}&&\M_{\mu/\beta}
\ddto^{p_{\mu,\lambda}}\\
\\
\M_{\lambda/\alpha}\rrto_{j_{\alpha,\beta}}&&\M_{\lambda/\beta}\\
\enddiagram}}
\label{second}
\end{equation}
Given such data, define an ``infinite-dimensional variety'' $\widehat\M$
as follows. For each $\alpha\in D$, the projections
$p_{\nu,\mu}:\M_{\nu/\alpha}
\longrightarrow\M_{\mu/\alpha}$, $(\mu\leq\nu)$, form a projective system
of affine bundles over $\M_{\alpha/\alpha}$, and we set $\widehat\M_\alpha :=
\proj_{\alpha\leq\mu}\M_{\mu/\alpha}$. Thus, $\widehat\M_\alpha$ is a
pro-algebraic variety, a bundle over $\M_{\alpha/\alpha}$ whose fibre is
a projective limit of affine linear spaces.

Further, the imbeddings $j_{\alpha,\beta}:\M_{\mu/\alpha}
\longrightarrow\M_{\mu/\beta}$, $(\alpha\leq\beta)$, give rise to
an direct system of open imbeddings $\widehat\M_\alpha\hookrightarrow
\widehat\M_\beta$. We set $\widehat\M=\inj_\alpha\widehat\M_\alpha$
and endow $\widehat\M$ with direct limit topology.

To get a better understanding of the structure of $\widehat\M$, fix
two elements $\alpha$, $\beta\in D$ such that $\alpha<\beta$. We have
the following commutative diagram:

\begin{equation}
\vcenter{\hbox{\diagram
&\widehat\M_\alpha\rrto^{j}\ddlto_{p_\alpha}\ddrto^{p}&&\widehat\M_\beta\ddto_
{p_b}\\
\\
\M_{\alpha/\alpha}&&\M_{\beta/\alpha}\llto_{p_{\beta,\alpha}}\rto_{j_{\alpha,
\beta}}&\M_{\beta/\beta}\\
\enddiagram}}
\label{third}
\end{equation}

In the diagram, the maps $p_\alpha$, $p_\beta$ and $p$ are standard
projections of a projective limit onto its components. Further, identify
$\M_{\beta/\alpha}$ with a Zariski-open part of $\M_{\beta/\beta}$ via
the imbedding $j_{\alpha,\beta}$. Then the projection $p$ in (\ref{third})
becomes the restriction to $\M_{\beta/\alpha}$
of the affine bundle $p_\beta:\widehat\M_\beta
\longrightarrow\M_{\beta/\beta}$. Thus, the
structure of the imbedding $\widehat\M_\alpha\hookrightarrow\widehat\M_\beta$
is determined, up to fibre of the projection $p_\beta$ (a non-essential
infinite-dimensional affine space), by the low row of (\ref{third}).
In other words, to construct $\widehat\M_\beta$ one starts with
$\M_{\alpha/\alpha}$, takes an affine bundle $\M_{\beta/\alpha}$ over
$\M_{\alpha/\alpha}$, then attaches to it a closed subvariety
(${}=\M_{\beta/\beta}\,\setminus\,\M_{\beta/\alpha}$), and finally takes
an infinite-dimensional affine bundle over the resulting space. It is clear from
the description above that $\widehat\M_\alpha$ is an open part of
$\widehat\M_\beta$ so that the $\widehat\M$ is the union of the family
of increasing open pieces $\{\widehat\M_\alpha\}$.

\setcounter{subsubsection}{\value{prop}}
\subsubsection{Constructible complexes on $\widehat\M$.}
For any $\lambda\in D$, define a triangulated category 
$D^b_\lambda(\widehat\M)$
as the category formed by all families  $\{F_{\alpha/\alpha}\in
D^b(\M_{\alpha/\alpha}),\ \alpha\geq\lambda\}$ that satisfy the following
compatibility condition
\begin{equation}
p^*_{\beta,\alpha}(F_{\alpha/\alpha})\simeq 
j^*_{\alpha,\beta}(F_{\beta/\beta})
\quad\mbox{for any }\beta>\alpha\geq\lambda  \label{forth}
\end{equation}
(where the maps $p_{\beta,\alpha}$ and $j_{\alpha,\beta}$ are the same as
 in the low row of (\ref{third})).

 For any pair $\lambda\leq\mu$ there is an obvious exact functor
 $D^b_\lambda(\widehat\M)\longrightarrow D^b_\mu(\widehat\M)$
 assigning to a family $\{F_{\alpha/\alpha},\ \alpha\geq\lambda\}$
 its part with indices $\alpha\geq\mu$. So the categories
 $ D^b_\lambda(\widehat\M)$ form an ``direct system'' and we set
 $ D^b(\widehat\M)=\inj_\lambda D^b_\lambda(\widehat\M)$. Similar
 construction works for perverse sheaves. So, there is a well-defined
 notion of a perverse sheaf on $\widehat\M$.

 Speaking formally, giving a constructible complex $F$ on $\widehat\M$
 requires giving a constructible complex $F_{\beta/\alpha}\in
  D^b(\M_{\beta/\alpha})$, one for each pair $\beta\geq\alpha(\geq\lambda)$,
 so that these complexes should be compatible in a natural way with
 all diagrams (\ref{first}) and (\ref{second}). However, one can see easily
 that such a collection $\{F_{\beta/\alpha}\}$ is completely determined
 by its part $\{F_{\alpha/\alpha},\ \alpha>\lambda\}$. Namely,
 given a collection $\{F_{\alpha/\alpha}\}$ satisfying (\ref{forth}),
 one can put $F_{\beta/\alpha}:=p^*_{\beta,\alpha}(F_{\alpha/\alpha})$ and
 objects so defined will be automatically compatible with diagrams
 (6.1.1) and (\ref{second}).

\subsubsection{Homology and Cohomology of $\widehat\M$.}
For any $\alpha$, the system of affine fibrations $\M_{\mu/\alpha}
\longleftarrow\M_{\nu/\alpha},(\mu\leq\nu)$, gives rise to an direct
system of cohomology isomorphisms: $H^i(\M_{\mu/\alpha})\stackrel{\sim}
{\longrightarrow}H^i(\M_{\nu/\alpha})$. We set $H^i(\widehat\M_\alpha):=
\inj H^i(\M_{\mu/\alpha})$. Now, the imbeddings $\widehat\M_\alpha
\hookrightarrow\widehat\M_\beta$ give rise to a projective system of
cohomology morphisms and we put $H^i(\widehat\M_\alpha):=
\proj_\alpha H^i(\widehat\M_\alpha)$.

Notice that there is an isomorphism $H^\bullet(\widehat\M_\alpha)\simeq
H^\bullet(\M_{\alpha/\alpha})$ induced by the projection. So, in down
to earth terms, giving a cohomology class of $\widehat\M$ amounts to giving
a collection $\{c_\alpha\in H^\bullet(\M_{\alpha/\alpha})\}$ such that,
for any pair $\beta>\alpha$, we have $p^*_{\beta,\alpha}(c_\alpha)=
j^*_{\alpha,\beta}(c_\beta)$ where $p_{\beta,\alpha}$ and $j_{\alpha,\beta}$
are as in (\ref{third}).

Similarly, one defines $H^{BM}_k(\widehat\M)$,
{\it $k$-codimensional} Borel-Moore homology group, to be
formed by collections \[\{c_\alpha\in H^{BM}_{d_\alpha-k}(\M_{\alpha/\alpha}),\
d_\alpha=\dim_{\Bbb R}\M_{\alpha/\alpha},
\ \alpha\in D\}\]
such that $p^*_{\beta,\alpha}(c_\alpha)=j^*_{\alpha,\beta}(c_\beta)$
(the pull-back morphism $p^*_{\beta,\alpha}$ is well defined for smooth
maps and shifts degree by the fibre dimension).

\subsection{Case of a group action}
Let $\M$ be an ``infinite dimensional'' space acted on by a
``pro-algebraic'' group $L$. More precisely, assume the following:
\begin{enumerate}
\item There is a (decreasing) family of normal subgroups
$\{\,L^\alpha\subset L\mid \alpha\in D\,\}$ indexed by an inductive set $D$ 
such that
 $\alpha\leq\beta\Rightarrow L^\beta\subset L^\alpha,(\alpha,\beta\in D)$,
 and such that $\cap_{\alpha\in D}L^\alpha=\{1\}$;
\item There is an exhaustion of $\M$ by open subsets 
$\{\,\M_\alpha\subset\M\mid
\alpha\in D\,\}$, indexed by the same set $D$ as in (i), such that
$\alpha\leq\beta\Rightarrow\M_\alpha\subset\M_\beta$; and moreover
$\cup_{\alpha\in D}\M_\alpha=\M$;
\item For every $\alpha\in D$ the group $L^\alpha$ acts freely on $\M_\alpha$
and the orbit-space $L^\alpha\backslash\M_\alpha$ has the structure of a
finite-dimensional smooth algebraic variety;
\item For any $\alpha<\beta$ the quotient $L^\alpha/L^\beta$ has the
stucture of a finite dimensional
 {\it unipotent} algebraic group; furthermore the induced action
of the group $L^\alpha/L^\beta$ on $L^\beta\backslash\M_\beta$ is
algebraic.
\end{enumerate}

Given a space $\M$ with $L$-action satisfying conditions 6.2(i)--(iv) above,
define a set of data 6.1(i)--6.1(iv) as follows. For each pair $\alpha\leq\beta$
put $\M_{\beta/\alpha}:=L^\alpha\backslash\M_\beta$. Conditions 6.2(iii,iv)
show that $\M_{\beta/\alpha}$ has the stucture of a smooth algebraic variety.
Furthermore, there are natural imbeddings $j_{\alpha,\beta}:L^\gamma\backslash
\M_\alpha\hookrightarrow L^\gamma\backslash\M_\beta,(\alpha<\beta\leq\gamma)$,
and natural projections $p_{\nu,\mu}:
L^\nu\backslash\M_\lambda\longrightarrow 
L^\mu\backslash\M_\lambda,(\lambda\leq
\mu<\nu)$. These maps clearly satisfy (\ref{first})--(\ref{second}).
Thus, we are in the setup of n.~6.1, so that the space $\widehat\M$ can be
defined. Now, any point $x\in \M_\alpha$ gives a point in the projective
limit $\widehat\M_\alpha=\proj_\mu\M_{\mu/\alpha}$. This way one gets an
imbedding $\M\hookrightarrow\widehat\M$ with dense image. We endow $\M$
with the topology induced from the topology on $\widehat\M$ via the imbedding.

The spaces $\M$ and $\widehat\M$ are very much alike. For each $\alpha\in D$,
for instance, there are projections 
$\M_\alpha\longrightarrow\M_{\alpha/\alpha}$
and $\widehat\M_\alpha\longrightarrow\M_{\alpha/\alpha}$. Both of them are 
affine
fibrations with infinite-dimensional fibres. The only difference between them
is that the fibre of the first fibration is isomorphic to the group $L^\alpha$
while the fibre of the second is isomorphic to
$\proj_\beta(L^\alpha/L^\beta)$, the completion 
of~$L^\alpha$.

Assume, in addition to properties 6.2(i)--(iv), that the following holds:
\begin{enumerate}
\setcounter{enumi}{4}
\item For each $\alpha\in D$, the group $L/L^\alpha$ has a
structure of a finite-dimensional algebraic group and that structure
is compatible with isomorphisms $L/L^\alpha\simeq
(L/L^\beta)/(L^\alpha/L^\beta)$, for
all pairs $\alpha<\beta$;
\item For any $\alpha\in D$, $\M_\alpha$ is an $L$-stable subset of $\M$,
and the induced $L/L^\alpha$-action on $L^\alpha\backslash\M_\alpha$ is
algebraic.
\end{enumerate}

Given an $L$-action on $\M$ satisfying properties 6.2(i)--(iv), one can make
$L\backslash\M$, the orbit-space, into a stack. We will never use the stack
approach, however, replacing it by the construction of n.~6.1 instead.
For example, we {\it define} a perverse sheaf on $L\backslash\M$ to be
an $L$-equivariant perverse sheaf on $\M$, that is a collection $\{F_\alpha,\
\alpha\in D\}$ with isomorphisms of $p^*_{\beta,\alpha}(F_\alpha)\simeq
j^*_{\alpha,\beta}(F_\beta)[\dim(L^\alpha/L^\beta)]$, for each pair
$\alpha<\beta$, (cf.~(\ref{forth})).

\subsection{Moduli space of G-bundles} 
Let $X$ be a smooth complex connected algebraic curve.
By a $G$-bundle we mean an algebraic principal $G$-bundle 
(with $G$-acting on the right) which is locally trivial
in the \`etale topology on $X$.

\begin{Th}{\it
\label{triv}
 Let $X$ be a smooth affine complex algebraic curve and $G$ a semisimple complex
connected  group. Then, any algebraic principal $G$-bundle on $X$ is
trivial.
}\end{Th}

\begin{remark}
(i) The theorem is  false for general connected
reductive groups.
For example if $G=GL_1$, a $G$-bundle is a line bundle, and there are plenty
of non-trivial line bundles on an affine curve.

(ii) 
If $G$ is semisimple and simply-connected group then a version of
Hilbert's `theorem 90' conjectured by Serre in "Cohomologie
Galoisienne"  (ch. III.14)
says that, for the field $\C(X)$ of rational
functions on a curve $X$, we have 
$$H^1(\mbox{Gal}\ (\mbox{{\small ${\Bbb C(X)}$}}), G(\C))=0\eqno(*)$$
 This cohomology vanishing
was 
proved, e.g., in \cite{BS}. It
ensures that a $G$-bundle which is locally trivial
in the \`etale topology is actually locally trivial
in the Zariski topology. The same actually holds for
a non-simply group $G$ as well (as was pointed out to me by
Telemann). Indeed, if $G$ is semisimple,
then $\pi_1(G)$ is a finite abelian group. In such a case, one has
$H^2(X,\pi_1(G))=0$. Now let
 $\tilde G$ is the
simply connected cover of $G$.
The following short exact sequence, viewed as a sequence of constant sheaves
on~$X$
\[0\to\pi_1(G) \to \tilde G \to G\to 0\]
yields 
the long exact sequence of cohomology. The latter, combined with 
the vanishing of $H^2(X,\pi_1(G))=0$, shows that equation
(*) holds for $G$ provided it holds for  $\tilde G$.

Thus, we may (and will) assume below
all $G$-bundles to be trivial at the generic point of the curve
$X$.

(iii) Originally, theorem 6.3.1 was deduced from a so-called
`strong approximation theorems' for adeles, cf. \cite{adel}. 
That approach gives a bit
stronger result but seems to be less geometric.
\end{remark}

Proof of the theorem proceeds in several steps.
 
 {\sc Step 1.} We need some simple general facts. Let $X$ be an arbitrary
 smooth algebraic variety and $G$ a linear algebraic group.
 Let $H$ be an algebraic subgroup of $G$
 and $P$ a principal $G$-bundle on $X$. We say that the structure
 group of $P$ can be reduced from $G$ to $H$ if there is a principal
 $H$-bundle $P'$ such that $P \simeq  P'\times_H G$.
 
 The following result can be easily derived  from definitions.
 
\begin{lemma}
\label{lm1} The structure group of a $G$-bundle $P$ on $X$ may be reduced
 to $H$ if and only if the associated bundle
 $ P/H= P \times_G (G/H)$ has a regular global section. $\quad\square$
\end{lemma}

 From lemma \ref{lm1} one obtains in particular the following result.
 
\begin{lemma}
\label{lm2}
 Assume that $H$ is a {\it normal} subgroup of $G$ and $P$ a $G$-bundle.
If $P/H$ is trivial as a
$G/H$-bundle then $P$ is trivial as a $G$-bundle. $\quad\square$
\end{lemma}

Now let $H={\bold G_a}$ be the additive group. Then, $H$-bundles on
$X$ are classified by the first sheaf cohomology
$H^1(X, \cal O_X)$. If $X$ is affine, then this group vanishes,
by Serre's theorem. It follows, that any $\bold G_a$-bundle on
an affine algebraic variety is trivial. More generally, we have

\begin{lemma}
\label{lm3} Let $H$ be a unipotent group. Then any $H$-bundle on an
affine algebraic variety is trivial.
\end{lemma}

Proof of lemma: We may choose a filtration $1=H_0 \subset H_1 \subset \ldots
H_n=H$ where the $H_i$ are normal subgroups of $H$ such that
$H_i/H_{i-1}\simeq {\bold G_a}$ for all $i=1,2,...$. 
We prove by descending
induction on $i$ that $P/H_i$ is trivial $G/H_i$-bundle. Assume, this is 
already proved for some $i$. To carry out the induction step observe first
that since $P/H_i$ is trivial, the structure group of $P$ reduces to
$H_i$, by lemma 1. Hence, we may assume that $P$ is an $H_i$-bundle.
Then the associated bundle $P/H_{i-1}$ has structure group $H_i/H_{i-1}=
{\bold G_a}$, hence is trivial, by the remark preceeding the lemma.
Applying lemma \ref{lm2} again, we see that the structure group of $P$
can be reduced from $H_i$ to $H_{i-1}$. That completes the proof.
$\quad\square$

\vskip 5mm
{\sc Step 2.} Now let $G$ be a connected semisimple group and $\B$ the Flag manifold
for $G$, i.e. the variety of all Borel subgroups in $G$. Choosing a Borel
subgroup $B \subset G$ yields a $G$-equivariant identification
$\B\simeq G/B$. Recall that $\B$ is a projective variety.

\begin{lemma}
\label{lm4}
 Assume that $X$ is a smooth compact algebraic curve.
Then, for any $G$-bundle $P$, the associated bundle
$P\times_G \B = P/B$ has a regular section.
\end{lemma}

Proof: By definition, the bundle $P$ is trivial on $U$, a Zariski
open subset of $X$. Hence, $(P\times_G \B) |_{_U} \simeq U\times \B$,
so that there always exists a regular section $s: U \to P\times_G \B$ 
over $U$ (e.g. a constant section). 
We claim that $s$ can be extended to a regular section all over $X$.
To that end, observe that $X \setminus U$ is a finite set of ponts,
for $X$ is 1-dimensional. Thus, we have to extend $s$ to each point
$x \in X \setminus U$ separately. The problem being local over a
neiborhood of $x$, we may assume the bundle $P$ to be trivial.
Then $s$ becomes a regular map from a punctured neighbourhood of
$x$ to $\B$. But any such map can be uniquely extended, by continuity,
to the point $x$, for $\B$ is a projective variety.$\quad\square$

\begin{lemma}
\label{lm5}(theorem in rank 1 case). Any $SL_2$- or 
$PGL_2$-bundle on an affine curve is 
trivial.
\end{lemma}

Proof: We view an $SL_2$-bundle as a rank 2 vector bundle, $V$, with
trivial determinant. By lemma \ref{lm4}, the structure group of this bundle
can be reduced to the subgroup of upper-triangular $2\times 2$-matrices.
Hence, $V$ is an extension of line bundles
$L \to V \twoheadrightarrow L'$.
This extension is controlled by an element of the first cohomology group
$H^1(X, \cal Hom(L',L))$.
The latter group vanishes, for $X$ is affine. Hence, $V\simeq L\oplus
L'$. 

Since $X$ is an affine curve, there exist regular
sections $s$ and $s'$ of the line bundles $L$ and $L'$ respectively
with disjoint
zero sets. Then $s \oplus s'$ is a nowhere vanishing section
of $V$. This section gives an imbedding of $\cal O_X$ as a trivial
subbundle in $V$. The quotient bundle, $V/\cal O_X$, has trivial
detrminant ($\det V=1$), hence is itself trivial. Therefore,
$V$ is an extension of $\cal O_X$ by $\cal O_X$. Such an extension
splits by the same cohomology vanishing as above. This proves the
lemma for $SL_2$-bundles. The $PGL_2$-case is reduced to the previous
one by remark (ii) after the statement of the theorem.
$\quad\square$

\vskip 5mm
\begin{proof}{Proof of the Theorem}
Choose a Borel subgroup $B$ in $G$ and write
$U$ for the unipotent radical of $B$ and $T=B/U$ for the
`universal maximal torus'. 
Applying lemma \ref{lm1} and lemma \ref{lm4}
to the subgroup $H=B$ we see that the structure group of $P$ may
be reduced from $G$ to $B$. Assuming this, we view $P$ as a $B$-bundle
and form the associated
$T$-bundle $P/U$. It suffices to show that this $T$-bundle is
trivial. Indeed, in that case $P/U$ has a regular section, hence
the structure group of the $B$-bundle $P$ can be reduced to $U$.
But any $U$-bundle on $X$ is trivial by lemma \ref{lm3}, and we are done.
Thus, the point is to find a reduction of $P$ from  $G$ to $B$ such that the 
associated
$T$-bundle is trivial. 

Let $\lambda : T=B/U \to \Bbb C^*$ be a dominant weight,
 $\lambda : B \to \Bbb C^*$ its pull-back to $B$, and
${\cal O}(\lambda)$ the induced $G$-equivariant line
bundle on $\B=G/B$. Given a principal
$G$-bundle $P$, form the associated bundle
$P\times_G \B$ with fiber $\B$. The fiber of this bundle is
(non-canonically) isomorphic to $\B$, and there is
a canonical line bundle on the total space of
$P\times_G \B$ whose restriction to each fiber is $\cal O_\lambda$.
Let  $L_\lambda$ denote that line bundle
on the total space.

Recall we are looking for a reduction of the structure
group of $P$ from $G$ to $B$ such that the associated
$T$-bundle is trivial. It is easy to see via lemma \ref{lm1}, that
getting such a reduction amounts to finding a regular
section $s : X \to P\times_G \B$ such that, for each dominant
weight $\lambda$, 
the pull-back $s^* (L_\lambda)$ is a trivial line bundle on $X$.

Start with any section $s': X \to P\times_G \B$, which exists
by lemma \ref{lm4}. This 
reduces the structure group to $B$. Fix some fundamental weight $\alpha$.
Let $P_\alpha\supset B$ be the minimal parabolic in $G$ corresponding to
$\alpha$. We may view the $B$-bundle $P$ as a $P_\alpha$-bundle
We want to find a new section $s''$ such that : 
(1) $(s'')^* (L_\beta)= (s')^*(L_\beta)$ whenever $ \beta$ is a fundamental weight,
$\beta  \neq \alpha$; and 
(2) $(s'') ^* (L_\alpha)$ is trivial. 
We look for $s''$ among the sections, whose image in the associated
$G/P_{\alpha}$-bundle coincides with the image of $s'$.
Then condition (1) above holds automatically. Note 
that if $U_\alpha$ denotes the unipotent radical of
the parabolic $P_\alpha$ then $P_\alpha/U_\alpha=SL_2$.
Lemma \ref{lm5} applied to the $P_\alpha/U_\alpha$-bundle
$P/U_\alpha$ yields condition (2).
Repeating this process for all fundamental weight one by one, we
find a section $s$ such that $s^*(L_\alpha)$ is trivial for
all fundamental weights $\alpha$, hence is trivial for all dominant
weights. The theorem is proved. 
\end{proof}
\vskip 3mm
Write $\goth g$ for the Lie algebra of the group $G$.
 Given a $\C$-algebra $A$,
 let $G(A)$ denote the group of $A$-rational points of $G$.
Similar notation, $\goth g(A)$, will be used in the Lie algebra case
even for a not necessarily unital algebra $A$.
 Recall the algebras $\K_{\,x}\,,\,\O_{x}\,,\,\MI_{x}$
introduced in $n^\circ 1.5$.
Thus, wee have the Lie algebra $\g(\K_{\,x})$ with
subalgebra $\g(\O_{x})$. The Lie algebra $\g(\O_{x})$ has the
chain of ideals $\g(\O_{x})\supset\g(\MI_{x})\supset\g(\MI^2_{x})
\supset\cdots$. 
Similarly there are groups $G(\K_{\,x})$ and $G(\O_{x})$.
The ring homomorphism $\O_{x}\to \O_{x}/\MI^m_{x}$
induces a group homomorphism $\pi: G(\O_{x})\to G(\O_{x}/\MI^m_{x})$.
We 
put $L^m := \pi^{-1}(1)$, the inverse image of the identity in $G(\O_{x}/\MI^m_{x})$.
The groups $L^1\supset L^2\supset\cdots$ form 
a chain of normal (congruence) subgroups in
$L^0=G(\O_{x})$. Observe that,
for any $j>i>0$, the quotient $L^i/L^j$ has an obvious structure of a
finite-dimensional unipotent group.

Recall (see $\S 1.5$), the imbedding $\O_{\text{out}}
\hookrightarrow \K_{\,x}$
 gives rise to a group imbedding $G(\O_{\text{out}}):=
G(\O_{\text{out}})\hookrightarrow G(\K_{ x})$. We set $\M=G(\K_{ x})/G(\O_{\text{out}})$.
There is a natural left $G(\K_{ x})$-action on $\M$. Further, for each
$i>0$, define a subset $\M_i\subset \M$ as follows:
\[
\M_i=\{\,f\cdot G(\O_{\text{out}})\in\M=G(\K_{ x})/G(\O_{\text{out}})\mid \g(\MI^i_{x})
\cap (\mbox{Ad}\,f)\g(\O_{\text{out}})\,=0\,\}
\]

\begin{prop}{}
{
\renewcommand{\labelenumi}{\theenumi}
\endgraf
\begin{enumerate}
\item There is a natural bijective correspondence between the set
$G(\O_{x})\backslash\M$ of $G(\O_{x})$-orbits on $\M$
and the set of isomorphism classes of algebraic principal $G$-bundles on $X$;
\item For any $i\geq1$, the group $L^i$ acts freely on $\M_i$ and the
orbit space $L^i\backslash\M_i$ has natural stucture of a smooth algebraic
variety;
\item Each piece $\M_i$ is a $G(\O_{x})$-stable subset of $\M$.
\end{enumerate}
}
\end{prop}

\begin{proof}{Sketch of Proof}
The punctured curve $X\,\setminus\,\{x\}$
being affine, any algebraic $G$-bundle on $X\,\setminus\,\{x\}$
is trivial by theorem \ref{triv}. Hence, any algebraic $G$-bundle $P$ on $X$ has
regular section, say $s_{\text{out}}$, over $X\,\setminus\,\{x\}$.
Let $s_x$ be a local section of $P$ on a neighborhood of $x$. The bundle $P$
is completely determined by the  transition function
$f=s_{\text{out}}\cdot s^{-1}_x\in G({\Bbb K}_{x})$.  Choosing other sections
$(s_{\text{out}},s_x)$ changes the element $f$ within the
double coset $G(\O_{x})\backslash G(\K_{ x})/G(\cal O_{\text{out}})$. 
That proves (i).
Now, the orbit space $L^i\backslash\M$ can be identified via the 
correspondence
of part (i) with the set of the isomorphism classes of pairs $(P,s)$, where
$P$ is a $G$-bundle on $X$ and $s$ is an $i$-th jet of local section of $P$
on a neighborhood of marked points. Part (ii) now follows from \cite{Sorger}. Part (iii)
is clear.
\end{proof}

Thus, the couple $(G(\O_{x}),\M)$ satisfies conditions 6.2(i)--(vi)
so that all the constructions of section 6.1 are applicable. In particular, 
there
is a topology on $\M$ induced from that on $\widehat\M$.

\subsection{Genus 0 case}
Take $X=\CP$ with the single marked point $x=\infty$. View $S^1$ as the
unit circle in $\C\subset\CP$. That gives a group imbedding $LG\hookrightarrow
G(\K_{\,x})=G({\Bbb K}_{\,\infty})$.
 Furthermore, we have $L^+G=G(\cal O_{\text{out}})$ so that
one gets an imbedding $\Gr=LG/L^+G\hookrightarrow G(\K_{\,x})/G(\cal
O_{\text{out}})
=\M$.
The imbedding has a dense image. Choose two opposite Borel subgroups $B^+$
and $B^-$ in $G$ so that $B^+\cap B^-=T$ is a maximal torus. Define the 
Iwahori
subgroups $I^\pm\subset G(\K_{\,x})$ by
\[
I^+=\{\,f\in G(\cal O_{\text{out}}) \mid f(0)\in B^+\,\}\mbox{ and }
I^-=\{\,f\in G(\O_\infty)\mid f(\infty)\in B^-\,\}.
\]
For each homomorphism $\lambda\in X_*(T)$, viewed as a point of $\M$, let
$C_\lambda=I^+\cdot\lambda$ and $C^\lambda=I^-\cdot\lambda$ be the
corresponding orbits of the Iwahori subgroups. These orbits are called
{\it Bruhat cells}. The orbits $C_\lambda$ are contained in $\Gr\subset\M$
and form a cell decomposition of the Grassmannian indexed by the lattice
$X_*(T)$. This cell decomposition defines the decomposition of $\Gr$ into the
$L^+G$-orbits $O_\lambda$; in fact, for any $W$-orbit $|\lambda|\subset
X_*(T)$ we have $O_{|\lambda|}=\sqcup_{\mu\in|\lambda|}C_\mu$. The closures
$\overline C_\lambda$ clearly form a base of $H_\bullet(\Gr)$.

The orbits $C^\lambda$ are not contained in $\Gr$ and have infinite 
dimensions.
However, $\overline C^\lambda$, the closure (in $\M$) of such an orbit, gives 
a
well-defined Borel-Moore homology class in $H_\bullet^{BM}(\M)$ in the
sense of n.~6.1.

We now define a bilinear cap product map
\begin{equation}
\label{6.4.1}
H^{BM}_i(\M)\times H_k(\Gr)\longrightarrow H_{k-i}(\Gr)
\end{equation}
where $\Gr$ is taken with the direct limit topology as in n.~1.7.
Let $Z\in H^{BM}_i(\M)$ and $c\in H_k(\Gr)$ The cycle $c$ clearly
has compact support. Hence, there is $n\gg0$ such that 
$\mathop{supp}c\subset\M_n$
where $\M_1\subset\M_2\subset\cdots$ is the exhaustion of $\M$
(cf.~n.~6.1). Let $Z_n$ be the representative of $Z$ in $H_\bullet(\M_{n/n})$
and $p_n:\M_n\longrightarrow\M_{n/n}$ the projection (notation of
diagram \ref{third}). We assume, replacing the cycle $c$ by a homotopy
equivalent cycle if necessary, that $p_n$ maps $c$ isomorphically onto
$p_n(c)$ and that $Z_n$ intersects $p_n(c)$
transversely. Then $c$ meets $p_n^{-1}(Z_n)$ transversely and we let
$Z\cap c$ denote the cycle in $H_{k-i}(\Gr)$ given by that intersection.
As a special case of (\ref{6.4.1}) we have a natural pairing
\[
H^{BM}_i(\M)\times H_i(\Gr)\stackrel{\cap}{\longrightarrow}H_0(\Gr)=\C
\]
which gives rise to a morphism:
\[
H^{BM}_i(\M)\longrightarrow H^i(\Gr).
\]
In particular, each cycle $\overline C^\lambda$ gives rise to a cohomology
class of $\Gr$. The following is known~\cite{PS}:
\begin{lemma}
The cell $C^\lambda$ meets $C_\lambda$ transversally in a single point
$\lambda$;
moreover the closure $\overline C^\lambda$ does not meet any other cell 
$C_\mu$
with $\dim C_\mu\leq\dim C_\lambda$.\qquad $\Box$
\end{lemma}
It follows from the Lemma that the cycles $\overline C^\lambda$ form a basis
in $H^*(\Gr)$ dual to the basis $\{\overline C_\lambda\}$ in homology.

\subsection{Global nilpotent variety}
Below we write $Bun_G$ for the moduli stack of principal
$G$-bundles on $X$, a smooth complex compact curve of genus
$>1$. Recall that the cotangent space
$T^*_PBun_G$ at a point $P\in Bun_G$ is given by the
Kodaira-Spencer formula
\[T^*_PBun_G= H^1(X,{\goth g}_{_P})= H^0(X,{\goth g}^*_{_P}\otimes
\Omega^1_X)=H^0(X,{\goth g}_{_P}\otimes
\Omega^1_X)\,,\]
where the last equality depends on the choice of an invariant bilinear
form on $\goth g$.
 
 We call a (possibly singular) algebraic subvariety $Y$
of a smooth
symplectic
algebraic variety $(X, \omega)$ {\it isotropic}, if for any {\it smooth}
locally closed subvariety $W \subset Y$, we have $\omega|_{_W}=0$.
It can be shown that this definition is equivalent to the more
conventional
one: the subvariety
$Y$ is said to be
isotropic if the tangent space, $T_yY$, at any
regular
point $y\in Y$ is an isotropic subspace in $T_yX$. The above
definitions also apply to smooth stacks that can be
locally represented as a quotient
of a smooth algebraic variety modulo an action of an algebraic group.
The stacks $Bun_G$ and $T^*_PBun_G$ are of this type, as explained
in $n^\circ 6.3$.

Recall that the cotangent space $T^*_PBun_G$ has a natural symplectic
structure.
Following Laumon \cite{La2} define 
$\NN \,\subset \, T^*Bun_G$, the {\it global nilpotent
variety}, as follows 
\[\NN=\{(P,x)\in T^*Bun_G \;|\; x\in H^0(X,{\goth g}_{_P}\otimes
\Omega^1_X)\,, \,x\;\mbox{is nilpotent section.}\}\]

\begin{Th}{\it
\label{lagrange}
Regular points of each irreducible
component of the variety $\NN$
form a Lagrangian subvariety in $T^*Bun$.
}\end{Th}

\begin{rem} This theorem was first proved, in the
special case $G=SL_n$, by Laumon~\cite{La2}. Laumon's
argument cannot be generalized to arbitrary semisimple groups.
In the general case, the theorem was proved by Faltings
\cite[theorem II.5]{Fa}. The proof below seems to be
more elementary than that of Faltings; it is based on nothing but
a few
general results of Symplectic geometry.
\end{rem} 

Let $(X_1,\omega_1)$ and $(X_2,\omega_2)$ be complex 
algebraic symplectic
manifolds, and $pr_i : X_1\times X_2 \to X_i$ the projections.

\begin{lemma}
\label{lagr1}
Let $\Lambda_1 \subset X_1$ and
$\Lambda \subset X_1\times X_2$ be isotropic algebraic
smooth subvarieties
(the latter with respect to the symplectic form
$pr_1^*\omega_1-pr_2^*\omega_2$). Then the smooth
locus of  $\Lambda_2:= pr_2\left(pr_1^{-1}(\Lambda_1)\cap
\Lambda\right)$ is an isotropic subvariety of $X_2$.
\end{lemma}

\begin{proof}{Proof}
Set $Y:=pr_1^{-1}(\Lambda_1)\cap
\Lambda$. Then simple linear algebra shows that,
 for any $y\in Y$ the image of the tangent
map $(pr_2)_{_*}: T_yY \to T_{pr_2(y)}X_2$ is isotropic.
Let $W\subset \Lambda_2:=pr_2(Y)$ be an irreducible smooth subvariety.
Observe that the map $pr_2 :pr_2^{-1}(W)\cap Y \to W$ is 
{\it surjective}. Hence, there exists a non-empty smooth
Zariski-open dense  subset $U\subset 
\left(pr_2^{-1}(W)\cap Y\right)_{red}$ such
that the restriction $pr_2: U \to W$ has surjective differential
at any point of $U$. Therefore the tangent space at the
generic point of $W$ is isotropic. Whence the tangent space at
every point of $W$ is isotropic. Thus, any smooth subvariety
of $\Lambda_2$ is isotropic, and lemma follows.
\end{proof}

Let $f: M\to N$ be a morphism of smooth algebraic varieties.
Identify $T^*(M\times N)$ with $T^*M\times T^*N$ 
via the standard map multiplied by $(-1)$ on the factor
$T^*N$.
The cotangent space, $T^*X$, to any manifold $X$ has a canonical
1-form, usually denoted `$p dq$'.
The canonical
1-form on $T^*(M\times N)$ gets identified,
under the identification above, with `$p_1 dq_1 - p_2 dq_2$'.
We endow $T^*M\times T^*N$ with the symplectic form induced from that
on $T^*(M\times N)$ via this identification.

Introduce the closed subvariety
\[Y_f= \{(m,\alpha),(n,\beta)\in T^*M\times T^*N \;|\;
 n=f(m)\,,\,f^*\beta=0\}\,.\]

\begin{lemma}
\label{lagr2}
$Y_f$ is an isotropic subvariety in $T^*M\times T^*N$.
\end{lemma}

\begin{proof}{Proof} The conormal bundle to the graph of $f$ is the subvariety
\[\Lambda = \{(m,\alpha),(n,\beta)\in T^*M\times T^*N \;|\;
 n=f(m)\,,\,
\alpha=f^*(\beta)\}\,.\]
Observe that the canonical 1-form `$p_1 dq_1 - p_2 dq_2$' on
$T^*M\times T^*N$ vanishes identically on $\Lambda$.
Hence, $\Lambda$ is a Lagrangian, in particular, an
isotropic subvariety, and we may 
apply lemma \ref{lagr1}
to $X_2= T^*M\,,\, X_1 =T^*N$, and $\Lambda_1=T^*_NN=$the
zero-section and the $\Lambda$ above. Observe now that we have by definition
$Y_f=\Lambda\cap pr_1^{-1}(T^*_NN)$.
Hence by lemma \ref{lagr1} the subvariety $p(Y_f)$ is isotropic.
\end{proof}

\begin{lemma}
\label{lagr3}
If $M$ and $N$ are smooth algebraic stacks, and
$f: M\to N$ is a representable morphism of finite type, then the
assertion of lemma \ref{lagr2} still holds.
\end{lemma}

\begin{proof}{Proof} Due to locality of the claim we may (and will)
assume $N$ is quasi-compact.
Let ${\tilde N}$ be a smooth algebraic
variety and  ${\tilde N} \to N$ be a smooth surjective equidimensional morphism.
Set ${\tilde M} := M\times_N{\tilde N}$. 
Then
we have $Y_f \subset
T^*N\times_N M$, and
\[Y_f\times_N{\tilde N} \,\subset\, T^*N\times_N{\tilde N}\times_NM
\,\subset\, T^*{\tilde N}\times_NM\]
We must show that the image of
$Y_f\times_N{\tilde N}$ is an isotropic subvariety
in $ T^*{\tilde N}$.
Let ${F} : {\tilde M} \to {\tilde N}$ be the natural morphism,
and 
$Y_{_ {F}}\subset 
T^*{\tilde N}\times_{_{\tilde N}}{\tilde
M}$
the corresponding subvariety of lemma \ref{lagr2}. Observe
that $T^*{\tilde N}\times_{_{\tilde N}}{\tilde
M}=T^*{\tilde N}\times_NM$ and $Y_f\times_N{\tilde N}=Y_{_ {\tilde
f}}$.
Hence, lemma \ref{lagr2} shows that the image of $Y_f\times_N{\tilde
N}$
in $T^*{\tilde N}$ is isotropic.
The claim follows.$\quad\square$
\end{proof}

\begin{proof}{Proof of theorem \ref{lagrange}}
Choose a Borel subgroup $B\subset G$ with Lie algebra 
${\goth b}$. Write ${\goth n}$ for the nilradical of ${\goth b}$.
In the setup of lemma \ref{lagr3}
put $M=Bun_G$, the moduli stack of principal $G$-bundles on the
curve $X$, and $N=Bun_B$, the moduli stack of principal
$B$-bundles. Let
 $f: Bun_B \to Bun_G$  be the natural morphism.
Observe that, for any $P\in Bun_B$, by Serre duality we have:
\[T^*_PBun_B= H^1(X,{\goth b}_{_P})= H^0(X,{\goth b}^*_{_P}\otimes
\Omega^1_X)=H^0(X,{\goth b}_{_P}\otimes
\Omega^1_X)\]
It follows, since any nilpotent element of $\g$ is conjugate to 
${\goth n}$, that in the notation of  lemma \ref{lagr3}
we have $\NN=p(Y_f)$. Observe further that 
$Bun_B$ is the union of a countable family of open substacks of
finite type over $Bun_G$ each.
Thus,  lemma \ref{lagr3} implies that $\NN$ is the union
of a {\it countable} family of isotropic substacks.
But any union of a countable family of isotropic substacks
is itself isotropic, for the
field of complex numbers is uncountable.
 It follows that the set $\NN$ is isotropic.

Finally, Hitchin \cite{Hi} showed that the global nilpotent variety
is the special fiber (over 0) of a certain morphism $\pi :T^*Bun_G \to H$,
where $H$ is a complex vector space of dimension $\dim Bun_G$
(at this point we use that genus $X$ is $>1$).
Furthermore, since $Bun_G$ is an equi-dimensional smooth stack,
each irreducible component of $T^*Bun_G$ has dimension $\ge 2\dim
Bun_G$.
It follows that
any irreducible component of the special fiber
$\pi^{-1}(0)$ has dimension $\geq \dim Bun_G$. But we have proved that
each component of $\NN=\pi^{-1}(0)$ is an isotropic subvariety.
Thus, $\NN$ is Lagrangian.
\end{proof}

Remarks. (i) The above argument shows that $NN$ is a complete intersection
in $T^*Bun_G$, in paricular is Cohen-Macauley. 
Further, since dimension of any fiber is $\le$ dimension of the
special fiber, it follows
that every irreducible component of any fiber of $\pi$ has dimension
$\dim Bun_G$, hence, $\pi$ is flat.

(ii) Hitchin actually worked in the setup of {\it stable} $G$-bundles
and
not in the setup of stacks. But his construction of the map $\pi$
extends to the stack set-up verbatim. 
Notice that our argument did not use any additional properties
of the map $\pi$ established in \cite{Hi}.

\subsection{Sketch of proof of theorem 1.5.7}
Let $A$ be a perverse sheaf
(not necessarily local Hecke eigen-sheaf) satisfying  conditions (i)
and
(ii) of  theorem 1.5.7.
Observe first that each simple object $IC_\lambda \in P(Gr)$ is of
geometric origin, for it may be viewed as
 the intersection-cohomology extension
of a (shifted) constant sheaf on the orbit $O_\lambda$. Thus,
decomposition theorem (see \cite{BBD}) for the convolution $IC_\lambda \ast A$
applies. The theorem yields that, for any $M\in P(Gr)$, we have
$T_M(A)= \oplus L_j[m_j]$, where each $L_j$ is a simple perverse
sheaf on $X\times Bun_G$, and $[m_j]$ stands for the shift in the
derived
category.

{\bf Notation:} Given a complex $F\in D^b(Y)$,
write $SS(F)$ for the characteristic variety of $F$,
cf. e.g. \cite{KS}.
Thus $SS(F)$ is a Lagrangian cone-subvariety in $T^*Y$.

Assume now that $SS(A) \subset \NN$. Then, an easy direct calculation,
based on a theorem of Kashiwara \cite{KS} giving an estimate on
the characteristic variety of a proper direct image, shows that
\[SS\left(T_M(A)\right)\,\subset T^*_XX\times \NN  \subset T^*(X\times
Bun_G)\quad
,\quad \forall M\in  P(Gr)\]
Using the expression for $T_M(A)$ provided by the decomposition
theorem we find
$$SS\left(L_j\right)\,\subset T^*_XX\times \NN\quad
,\quad \forall M\in  P(Gr),\quad \forall j\eqno(6.6.1)$$
By theorem \ref{lagrange} we know that $\NN$ is a Lagrangian
subvariety in $T^*Bun_G$. Hence, there exists a stratification
$Bun_G=\amalg_\nu Bun_\nu$ such that any complex on
$X\times Bun_G$ whose characteristic variety is contained
in $T^*_XX\times \NN$ is locally constant along the
strata of the stratification $X\times Bun_G=\amalg_\nu$
$(X\times Bun_\nu)$.
Hence, any perverse sheaf $L_j$ that occurs in (6.6.1) is the
intersection
complex associated with a simple local system on a certain stratum
$X\times Bun_\nu$. Observe that $\pi_1(X\times Bun_\nu)=
\pi_1(X)\times \pi_1(Bun_\nu)$, and any simple
representation of a direct product of groups
 is the tensor product of simple representations of the factors.
We find that each perverse sheaf $L_j$ is of the form
$\cal L \boxtimes IC(Bun_\nu,\cal A_\nu)$ where $\cal L$ is an irreducible
local system on $X$ and $IC(Bun_\nu,\cal A_\nu)$ is the 
intersection complex associated with a simple local system $\cal A_\nu$
on ${Bun_\nu}$. Write $IC_j(M)$ for the possible
{\it non-isomorphic} sheaves
$IC(Bun_\nu,\cal A_\nu)$, that occur above with non-zero multiplicity,
and $\cal L(M)$ for the corresponding local system on $X$.
Thus, we obtain
$$T_M(A)= \oplus_j\, \cal L_j(M) \boxtimes IC_j(M)[m_j] \quad
,\quad \forall M\in  P(Gr)\eqno(6.6.2)$$
Next, fix a point
$x\in X$, and let $i_x :\{x\}\times Bun_G \hookrightarrow X\times
Bun_G$
denote the imbedding.
One checks from definitions that a global Hecke operator, $T_M$,
is related to the correspondinding local Hecke operator at $x$
by the following
natural isomorphism
$$ i^*_xT_M(F) = M\ast F \quad
,\quad \forall M\in  P(Gr),\forall F \in D^b(Bun_G)\eqno(6.6.3)$$
where convolution on the RHS is defined via the double-coset
isomorphism (1.5.2) that involves the choice of  $x$ in an essential way.

We now use the assumption that the perverse sheaf $A$ is a local
Hecke eigen-sheaf (at $x$). In view of equations (6.6.2)
and (6.6.3) the assumption
reads
$$\bigoplus_j \, i^*_x\cal L_j(M) \boxtimes IC_j(M)[m_j] =
M\ast A\simeq L_M\otimes A
\quad
,\quad \forall M\in  P(Gr)\eqno(6.6.4)$$
Since all the sheaves $IC_j(M)$ in equation (6.6.2) were assumed to be
{\it pairwise non-isomorphic} and the eigen-sheaf $A$ is assumed simple, we conclude
that isomorphism (6.6.4) can
hold only if the direct sum on the left of (6.6.4) contains a single 
non-zero summand, say for $j=j_o$
Whence, we have $T_M(A)= \cal L_{j_o}(M) \boxtimes IC_{j_o}(M)$
and, moreover,
$IC_{j_o}
(M)\simeq A$.
Thus, $A$ is a global Hecke eigen-sheaf, and the first claim of the
theorem
is proved.

To prove the rest of the theorem, we first argue locally at the point
$x$. Let $A$ be our  local  Hecke eigen-sheaf, so that
$M\ast A \simeq L_M \otimes A$ holds for any $M\in P(Gr)$. Since
$A$ is {\it simple}, this equation yields
$L_M \simeq Hom_{P(Gr)}(A, M\ast A)$. That shows that we may assume
without loss of generality 
the assignment $M \mapsto L_M$ to be a functor $P(Gr)\to \,\mbox{{\Bbb Vect}}\,$.
Observe next that, for any $M,N\in P(Gr)$ we have
\[(M\ast N)\ast A=M\ast (N\ast A)= M\ast (L_N\otimes A)=
L_N\otimes (M\ast A)=(L_M\otimes L_N) \otimes
A\]
Applying the functor $Hom_{P(Gr)}(A,\bullet)$ to both sides we see
that the functor $M \mapsto L_M$ is a tensor functor.

If $A$ is a {\it global} Hecke eigen-sheaf one can repeat the argument
of the previous paragraph `globally', replacing equation
$M\ast A \simeq L_M \otimes A$ by the equation
$T_M(A)\simeq \cal L_M\boxtimes A$. The argument shows that
the assignment $\cal L: M \mapsto \cal L_M$ gives a tensor  functor on $P(Gr)$
with values in $\cal Loc(X)$, the tensor category of locally constant
sheaves on $X$ with the standard tensor product.

Applying theorem 1.4.1, we obtain  the following
composition of functors:
$$\cal R=\pp\circ\cal L: Rep_{_{G^\lor}}\to
P(Gr)
\to
\cal Loc(X)$$
This composition gives a tensor functor $\cal R:
Rep_{_{G^\lor}}\to\cal Loc(X)$.
We claim that any such tensor functor $\cal R$ is isomorphic to the
functor $V \mapsto V_{_P}$ , for an appropriate principal
$G^\lor$-bundle
$P$ on $X$ with flat connection. To prove the claim, view a
flat $G^\lor$-bundle as a representation of the fundamental group,
$\pi_1(X,x)$. This way the functor $\cal R$ may be regarded as 
a tensor functor
\[\cal R:Rep_{_{G^\lor}}\to Rep_{\pi_1(X,x)}\]
But any such functor is known \cite{DM} to be induced by a 
group homomorphism $\pi_1(X,x)\to G^\lor$. This proves the claim,
and the theorem follows.$\quad\square$

\section{Appendix A: Tensor functors and $G$-bundles}

     Throughout the   Appendix  we  let  $G$  denote  a  linear
algebraic group over $\C$,  and $F_{G}:\Rep(G)\to\Vect$ the
forgetful functor. The reader should be warned that the results
below involving group $G$ are  applied  to  the  ``dual  group''
$G^{\lor}$ in the main body of the paper. I hope this will not
lead to confusion.

\subsection{} Let  $T$  be  another  group,  and $F_{T}$ the
corresponding forgetful functor on $\Rep(T)$.  We  have  the
following result \cite[Corollary 2.9]{DM}:

\begin{prop}
Let $S:\Rep(G)\to \Rep(T)$ be
a tensor functor such that $F_{G}\circ S=F_{T}$.
Then, there exists a unique  algebraic  homomorphism  $s:T\to G$
such that  the  functor $S$ is induced by $s$,  i.e.,  for any
representation $r:G\to\GL(V)$, we have: $S(r)=r\circ s$.
\endpr
\end{prop}

The following special case of Proposition 7.1.1 is particularly
useful. Let $T$ be a torus, and $X^{*}(T)$ the weight lattice.
Assume further, that for each representation $V\in \Rep(G)$ we
are given  a  gradation  on its underlying vector space by the
lattice $X^{*}(T)$:
\[
V=\mathop\oplus\limits_{\lambda}V(\lambda),\quad\lambda\in
X^{*}(T),
\]
which is   compatible   with   tensor   product, i.e. for  any
$V_{1},V_{2}\in\Rep(G)$ we have:
\begin{equation}
(V_{1}\oplus V_{2})(\lambda)=\sum_{\mu+\nu=\lambda} V_{1}(\mu)
\oplus V_{2}(\nu)
\end{equation}
Then one has:

\begin{corol} There   is   a   unique   homomorphism
$s:T\to G$ such that,  for any representation $V\in \Rep(G)$,
the gradation  $V(\cdot)$ coincides with the weight-gradation,
that is the gradation:
\begin{equation}
V(\lambda)=\{v\in V\mid    s(t)\cdot    v=
\lambda (t)\cdot v,\, t\in{\LARGE T},\lambda \in X^{*}(T)\}
\end{equation}
\end{corol}

\begin{proof}{Proof}
For each $V\in \Rep(G)$ define an action of
the torus $T$ on $V$ by letting $t\in T$ act on $V(\lambda)$  as
multiplication by  $\lambda(t)$.  In  this  way  we  obtain  a
functor $S:\Rep(G)\to\Rep(T)$,  which is a tensor  functor  by
(7.1.2). The result now follows from Proposition 7.1.1.
\end{proof}

\subsection{}
 Let $A=\mathop\oplus\limits_{i\ge0}A_{i}$ be a finitely-generated
commutative graded $\C$-algebra. The gradation on $A$ yields an
algebraic $\C^{*}$-action on the affine scheme $\Spec A$. Now, let
$P$ be a $\C^{*}$-equivariant principal algebraic $G$-bundle on
Spec A (that means that the group $\C^*$ acts freely on $P$ on the
right, the  group  $\C^{*}$ acts on $P$ on the left,  these two
actions commute and the projection: $P\to\Spec A$ commutes with
the $\C^{*}$-action). For any representation $V\in \Rep(G)$, we
can form the associated bundle $P\times_{G}V$,  which is a $\C
^{*}$-equivariant algebraic vector bundle on $\Spec A$. The
space $\Gamma_{P}(V)$ of its  global  sections has  a  natural
structure of    a   finitely-generated   graded  $A$-module. Let
$\mod^{\bullet}\!\!-\!\!A$ denote  the   abelian   category   of   finitely
generated graded   $A$-modules.    The   assignment   $V\mapsto
\Gamma_{P}(V)$ clearly defines an exact fully-faithful  tensor
functor (i.e. a fibre functor):
\[
\Gamma_{P}:\Rep(G) \to{\rm mod}^{\bullet}\!\!-\!\!A
\]

The following result is  a  $\C^{*}$-equivariant  analogue  of
\cite[Theorem 3.2]{DM}.

\begin{prop}
Any  fibre  functor:
$\Rep(G)\to{\rm mod}^{\bullet}\!\!-\!\!A$ is canonically isomorphic to the functor
$\Gamma_{P}$ for a uniquely  determined  $\C^{*}$-equivariant
principal $G$-bundle $P$.   \endpr
\end{prop}

\begin{corol}
Any fibre functor $F:\Rep(G)\to\Vect$  is
(non-canonically) isomorphic to the forgetful functor $F_{G}$.
\end{corol}

\begin{proof}{Proof}  Such  a  functor  $F$  may  be  regarded  as   a
fiber-functor: $\Rep(G)\to{{\rm mod}\!\!-\!\!\C}$. Then, Proposition 7.2.1 says
that there is a principal homogeneous space $P$ such that  the
functor $F$      is     isomorphic     to     the     functor:
$V\mapsto P\times_{G}V$. A choice of a point in $P$  yields  an
isomorphism: $P\cong G$          and,         hence,         an
isomorphism: $P\times_{G}V\cong G\times_{G}V\cong V$.
\end{proof}

\subsection{}
Assume  now  that  $G$ is a connected reductive Lie
group. In this $n^{\circ}$ we shall give a classification  of
$\C^{*}$-equivariant principal $G$-bundles on the line $\C$, where
the group $\C^{*}$ acts  on  $\C$  in  the  standard  way  (by
multiplication).

     Given a $\C^{*}$-equivariant principal $G$-bundle $P$ on
$\C$, we let $P_{t}$ denote the fiber of $P$ over a point $t\in
\C$. The fiber $P_{t}$ is a principal homogeneous $G$-space  and
the group $\Aut(P_{t})$ of its automorphisms is non-canonically
isomorphic to  $G$  (the  map  $P_{t}\to P_{t}$  is  called  an
automorphism if it commutes with the $G$-action).

The fiber  $P_0$  over  zero   is   clearly   a $\C^{*}$-stable
subvariety  of  $P$, and  the  action of $\C^{*}$
gives rise to a homomorphism $\tau:\C^{*}\to {\rm Aut}(P_{\circ})$.

     Let $P_{1}$ be the  fiber of $P$ over  the  unit.  To  each
$g\in\Aut(P_{1})$ we associate a family of  automorphisms
$\{g_{t}\in\Aut(P_{t}),\,t\ne 0\}$ defined as the composition:
\[
g_{t}:P_{t}\stackrel{t^{-1}}{\to}P_{1}
\stackrel{g}{\to}P_{1}\stackrel{t}{\to}P_{t},
\]
     where the first map is induced by the action  of  $t^{-1}
\in\C^{*}$, the second one by the action of $g$,  and the last
one by the action of $t\in\C^{*}$. Let $Q$ be the set of those
$g\in\Aut(P_{1})$ that the family $\{g_{t},\,t\in\C^{*}\}$ has a
limit: $g_{t}\to g_0\in\Aut(P_0)$ as $t$ approaches
$0$. It  is clear that $Q$ is a subgroup of $\Aut(P_{1})$ and that
the assignment   $g\mapsto   g_0$   yields    a    group
homomorphism $\lim:Q\to\Aut(P_0)$.

\begin{lemma}

\begin{enumerate}
\item The group $Q$ is a parabolic  subgroup  of
${\rm Aut}(P_{1})$;
\item The kernel of the  homomorphism  '$\lim$'  is  equal  to
$\rad Q$, the unipotent radical of $Q$;
\item The image of  the  homomorphism  'lim' is  equal  to
$L_{\tau}$, the centralizer in $\Aut(P_0)$ of the above
defined homomorphism $\tau:\C^{*}\to\Aut(P_0)$.
\end{enumerate}
\end{lemma}

Thus, to any $\C^*$-equivariant principal $G$-bundle $P$ on $\C$ we can
associate a triplet $(G_1,Q,\tau)$ where $G_1:=\Aut(P_1)$ is a reductive group
isomorphic to $G$, $Q$ is a parabolic subgroup of $G_1$ and 
$\tau:\C^*\rightarrow
Q/\rad Q$ is a regular central homomorphism (here ``central'' means that the
image
of $\tau$ belongs to the center of $Q/\rad Q$, and ``regular'' means that, for 
any
lifting $\tilde\tau:\C^*\rightarrow Q$, the centralizer in $G_1$ of the image 
of
$\tilde\tau$ belongs to $Q$).

The triplets $(G_1,Q,\tau)$, as above, form a category Tripl. A morphism:
$(G_1,Q,\tau)\rightarrow (G'_1,Q',\tau')$ in that category is, by definition,
a group isomorphism $G_1\simeq G'_1$ which maps $Q$ to $Q'$ and $\tau$ to 
$\tau'$.
\begin{prop}
The functor: $P\mapsto (G_1,Q,\tau)$ gives an equivalence of the category of
$\C^*$-equivariant principal $G$-bundles on $\C$ and the category 
Tripl.\hfill$\Box$
\end{prop}

Fix a $\C^*$-equivariant $G$-bundle on $\C$. Given $V\in \mathop{\fam0 
Rep}(G)$,
let $V_1$ denote the fibre over the point 1 of the associated bundle
$P\mathbin{\times_G}V$. The gradation on the space of global sections
of $P\mathbin{\times_G}V$ induced by the $\C^*$-action gives rise to a
 filtration $W$ on $V_1$. Namely, say that $v\in W_i(V_1)$ iff there
 exists a section $s\in\Gamma_P(V)$ of degree $\leq i$ such that $v=s(1)$.

 Next, note that the vector space $V_1$ can be identified with
 $P_1\mathbin{\times_G}V$, where $P_1$ is the fibre of $P$ over 1.
 Hence, the group $G_1=\Aut(P_1)$ acts naturally on $V_1$. We shall now
 characterize the filtration $W$ on $V_1$ in terms of the $G_1$-action.
 To that end, take the homomorphism $\tau:\C^*\rightarrow Q/\rad Q$
 associated with the bundle $P$ and let $h=\tau'(1)$ be the derivative of
 $\tau$ at the identity. Thus, $h$ is a semisimple element of the Lie
 algebra of the group $Q/\rad Q$. Chose $\widetilde h\in\mathop{\fam0 Lie}Q$,
 a semisimple representative of $h$. The element $\widetilde h$ acts on the
 space $V_1$ in a natural way and we have

\begin{prop} The filtration $W_\bullet(V)$ coincides with the filtration by 
the
eigenvalues of $\widetilde h$.
\end{prop}

\begin{rem} The filtration by the eigenvalues of $\widetilde h$ will
not change if $\widetilde h$ is replaced by $\wh+x$, $x\in\mathop{\fam0 Lie}
(\rad Q)$. Hence, it depends only on $h$ and not on its representative $\wh$.
\end{rem}

\section{Appendix B: Equivariant Hyper-Cohomology}
In this appendix we recall a number of equivalent constructions of the
derived category of equivariant complexes (equivariant derived category,
for short) and define equivariant hyper-cohomology of an equivariant
complex.
The reader is referred to the foundational work of Bernstein-Lunts
[BL]
for more details. Some of the constructions we are using have been
introduced by Lusztig [Lu 3] before [BL] appeared (cf. also [Lu 4]
for some more recent additional results). 

\subsection{}
Let $T$ be a Lie group and $ET\rightarrow BT$ a universal principal 
$T$-bundle,
so that $T$ acts freely on $ET$ on the right. There are (at least) two 
approaches
to this bundle.

The first approach, used in topology, is to view the universal bundle as an
direct limit of its finite dimensional approximations, i.e.\ to use a
diagram
\begin{equation}
\vcenter{\hbox{\diagram
ET^0\rto^{i_0}\dto^{T}&ET^1\rto^{i_1}\dto^{T}&ET^2\rto^{i_2}\dto^{T}&ET^3\rto^
{i_3}\dto^{T}&\cdots\\
\mbox{\it pt}\rto&BT^1\rto^{i_1}&BT^2\rto^{i_2}&BT^3\rto^{i_3}&\cdots\\
\enddiagram}}
\label{next}
\end{equation}
The vertical arrows in this diagram are finite-dimensional principal
$T$-bundles and the horizontal embeddings: $ET^n\hookrightarrow ET^{n+1}$
are $T$-morphisms. Such a diagram is called an approximation of the universal
bundle if the following holds: {\it  For any $k\geq0$ there is an integer
$n(k)\gg0$ such that all the spaces $ET^n$, $n\geq n(k)$, are 
$k$-contractible,
i.e.\ have vanishing homotopy: $\pi_j(ET^n)=0$ for all $j\leq k$.} The 
condition
implies that the cohomology of the spaces $BT^n$ stabilize, i.e.\ the 
imbeddings
$i_n:BT\hookrightarrow BT^{n+1}$ induce isomorphisms
\[
i^*_n:H^k(BT^{n+1})\stackrel{\sim}{\longrightarrow}H^k(BT^n)\qquad\mbox{for
any }k\leq n.
\]

The second, more algebraic, approach to the universal bundle is to view $ET$
as the standard simplicial scheme
\[
\vcenter{\hbox{\diagram
pt&\lto T&\lto<.3ex>\lto<-.3ex>T\times T&\lto\lto<.5ex>\lto<-.5ex>T\times 
T\times T\ldots\\
\enddiagram}}
\]
The diagonal $T$-action on each of the spaces above is obviously free, giving
a simplicial model for the universal $T$-bundle (i.e.\ a $T$ bundle in the
category of simplicial schemes).

\subsection{}
Let $Y$ be a variety with a smooth $T$-action. Set 
$Y_T=ET\mathbin{\times_T}Y$,
viewed either as an direct limit of the finite-dimensional spaces
$Y^n_T=ET^n\mathbin{\times_T}Y$ or as a simplicial scheme. Anyway, there are
diagrams
\begin{equation}
\vcenter{\hbox{\diagram
&ET\times Y\ddlto_{\pi}\ddrto^{\rho}&&&ET^n\times 
Y\ddlto_{\pi_n}\ddrto^{\rho_n}&\\
\\
Y_T&&Y&Y^n_T&&Y\\
\enddiagram}}
\label{another}
\end{equation}
and the projection $ET\longrightarrow BT$ gives rise to a fibration
$Y_T\longrightarrow BT$ with fibre $Y$.

Following Bern\-stein-Lunts, define a
$T$-equivariant
complex on $Y$ as a collection $(M,M^n_T,\phi_n, \psi_n, n\geq0)$ where
$M\in D^b(Y)$, $M^n_T\in D^b(Y^n_T)$ and $\phi_n$, $\psi_n$ are isomorphisms
$\phi_n:M^n_T\simeq i^*_n M^{n+1}_T$, $\psi_n:\rho^*_nM\simeq \pi^*_n M^n_T$
with certain compatibility conditions. Such collections
form a triangulated category $D_G(Y)$, {\it the equivariant derived category}
of $Y$. An object of $D_G(Y)$ is called {\it an equivariant complex}. An 
equivariant
morphism $X\rightarrow Y$ of $T$-varieties induces a compatible system of
morphism $X^n_T\longrightarrow Y^n_T$, hence, gives rise to all the standard
functors on the equivariant derived categories, e.g. direct images. inverse 
images,
etc.

In the simplicial approach, one defines an equivariant complex to be an object
of $D^b(Y_T)$, the bounded derived category of constructible complexes on the
simplicial scheme  $Y_T$. If the $T$-action on $Y$ is free so 
that
the orbit space $\overline Y=T\backslash Y$ is well-defined, then the second 
projection
$ET\times Y\longrightarrow Y$ induces a simplicial morphism
$F:Y_T\longrightarrow
\overline Y$ with fiber $ET$. Since $ET$ is a contractible scheme one deduces 
the
following.

\medskip
\noindent{\bf Equivariant descent: }{\it Let  $Y\longrightarrow\overline Y$ be 
a
locally-trivial principal $T$-bundle. Then the inverse image functor $F^*$
gives an equivalence of $D^b(\overline Y)$, the ordinary derived category, 
with
$D^b(Y_T)$.}

\medskip

A connection between the two approaches to equivariant complexes can be
established as follows. Given a diagram (\ref{next}) and a $T$-variety $Y$,
view $ET^n\times Y$, as a $T$-variety with the diagonal $T$-action.
Then for each $n\geq0$ we get the following simplicial analogue of diagram
(\ref{another}):
\[
\vcenter{\hbox{\diagram
&(ET^n\times Y)_T\ddlto_{\pi_n}\ddrto^{p_n}&\\
\\
\llap{$ET^n\mathbin{\times_TY={}}$}Y^n_T&&Y_T\\
\enddiagram}}
\]
where $p_n$ stands for the simplicial map induced by the second projection
$ET^n\times Y\longrightarrow Y$, and $\pi_n$ is the ``descent map F'' for the
variety $ET^n\times Y$ (with free action). Now given an object $M\in 
D^b(Y_T)$,
for each $n\geq0$, take $p^*_nM\in D^b((ET^n\times Y)_T)$. By the equivariant
descent property, there is a unique complex $M^n_T\in D^b(Y^n_T)$ such that
$P^*_nM=\pi^*_nM^n_T$. There are natural isomorphisms $\phi_n:M^n_T\simeq
i^*_nM^{n+1}_T$ and $\psi_n:\rho^*_nM\simeq \pi^*_nM^n_T$ so that the 
collection
$\{\,M, M^n_T,\phi_n,\psi_n;\,n\geq0\,\}$ gives an object of $D_G(Y)$. This
way one gets a functor $D^b(Y_T)\longrightarrow D_G(Y)$ that turns out to be 
an
equivalence of categories.

Restricting an object of $D^b(Y_T)$ to
\raisebox{4pt}{\hbox{\diagram
Y&\lto<.3ex>\lto<-.3ex>T\times Y&\lto\lto<.5ex>\lto<-.5ex>T\times T\times Y\\
\enddiagram}},
 the beginning of the simplicial scheme $Y_T$,
we get an ordinary complex $M\in D^b(Y)$ equipped with an isomorphism
\begin{equation}
I:p^*M\simeq q^*M \label{easy}
\end{equation}
($p, q: T\times Y\longrightarrow Y$ stand for the second projection and the 
action,
respectively). The isomorphism $I$ satisfies the cocycle condition:
\[
p_1^*I\circ q_2^*\simeq m^*I,\qquad\mbox{where $p_1, q_2,m:T\times T\times Y
\longrightarrow T\times Y$}
\]
are given by the formulas
$p_1:(t_1,t_2,y)\mapsto(t_2,y)$, $q_2:(t_1,t_2,y)\mapsto(t_1,t_2,y)$ and
$m:(t_1,t_2,y)\mapsto(t_1,t_2,y)$. There is also a normalization property:
{\it the restriction of the isomorphism $I$ to $1\times Y\subset T\times Y$
is the identity morphism $M\longrightarrow M$} (notice that both
$p^*M$ and $q^*M$ being restricted to $1\times Y$ become canonically 
isomorphic
to $M$).

It should be emphasized that, given an object $M\in D^b(Y)$ together with
an isomorphism (\ref{easy}) satisfying the conditions above, it is not 
possible
in general to extend these data to an object of $D^b(Y_T)$ due to the failure
of Grothendieck's descent for derived categories. Grothendieck's descent 
holds,
however, for perverse sheaves. Hence, we obtain:
\begin{equation}
\parbox{\quotesize}{
Giving an equivariant {\it perverse sheaf} on $Y$ is the same
as giving an ordinary perverse sheaf $M$ together with an
isomorphism (\ref{easy}) satisfying the conditions above.
}
\end{equation}
Let us also mention the following:
\begin{equation}
\parbox{\quotesize}{
The constant sheaf on $Y$ always has the structure of  an
equivariant complex. Such a structure is unique, provided the group $T$ is
connected (because of the normalization $I|_{1\times Y}=\mbox{id}$ in 
(\ref{easy})).
}
\end{equation}
\begin{rem} For a smooth variety $Y$ there is yet another construction
of the equivariant derived category based on ${\cal D}$-modules 
(see~\cite{Gi2}).
\end{rem}
\subsection{Equivariant cohomology}
Let $\{\,M,M^n_T,\phi_n,\psi_n\,\}$ be an equivariant complex on $Y$ in the
sen\- se of Bern\-stein-Lunts. For each $n$, the fibration $Y^n_T\longrightarrow
BT^n$
gives rise to the Leray spectral sequence:
\begin{equation}
E^{p,q}_2=H^p(BT^n)\otimes H^q(Y,M)\Longrightarrow H^{p+q}(T^n_T,M^n_T).
\label{threeone}
\end{equation}
These spectral sequences from a projective system corresponding to the 
inductive
sequence of $T$-fibrations $Y^0_T\hookrightarrow Y^1_T\hookrightarrow\ldots$.
The spectral sequences show that, for each $k\geq0$, the projective system
$H^k(M^0_T)\leftarrow H^k(M^1_T)\leftarrow\ldots$ stabilizes, due to the 
stabilization
of the cohomology of the spaces $BT^n$, $n=0,1,\ldots$.

Define the $T$-equivariant cohomology by the formula
\begin{equation}
H^k_T(M)=\proj_n H^k(Y^n_T, M^n_T) \label{threetwo}
\end{equation}
(abusing the notation we often denote an equivariant complex
$\{\,M,M^n_T,\phi_n,\psi_n\,\}$ by a single symbol $M$). Taking the projective
limit of spectral sequences (\ref{threeone}) we obtain the spectral sequence
for equivariant cohomology:
\begin{equation}
E^{p,q}_2=H^p(BT)\otimes H^q(Y,M)\Longrightarrow H^{p+q}_T(M).
\end{equation}

For $M=\C_Y$, the constant sheaf, the R.H.S. of (\ref{threetwo}) becomes
the standard definition of the equivariant cohomology of $Y$. Moreover, for 
any
equivariant complex $M$ on $Y$, there is a natural isomorphism:
\[
H^\bullet_T(M)=\mathop{\fam0 Ext}\nolimits^\bullet_{D_G(Y)}(\C_Y,M).
\]
This isomorphism is one of the main reasons for introducing the equivariant 
derived
category.

Let $T_c$ be a maximal compact subgroup of $T$. The inclusion
$T_c\hookrightarrow T$ induces a system of maps $Y^n_{T_c}\hookrightarrow
Y^n_T$ and homotopy equivalence $BT^n_c\approx BT^n$, $n=1,2,\ldots$.
It follows that, for any $T$-equivariant complex $M$, one has
\begin{equation}
H^\bullet_{T_c}(M)\cong H^\bullet_T(M).
\end{equation}

There is an equivalent definition of equivariant cohomology based on the
simplicial approach
. Given an object $M\in D^b(Y_T)$, put
\[
H^\bullet_T(M)=H^\bullet(Y_T,M)
\]
where the R.H.S. stands for the hyper-cohomology of a double-complex with
the second differential coming from the combinatorial differential of the
simplicial scheme. This definition turns out to be equivalent to
(\ref{threetwo}).
\subsection{}
Assume now that $T$ is a connected complex reductive Lie group. Then,
one can find a diagram (\ref{next}) consisting of algebraic
varieties and algebraic maps. Moreover, one can find that diagram in such
a way that all the spaces $BT^n$ in the diagram are {\it smooth
projective} varieties (e.g. if $T$ is a torus then $BT^n$ can be chosen
to be the direct product of $\dim T$ copies of $\CPN$).

Further, let $Y$ be projective variety with an algebraic $T$-action.
Then, the projection: $Y^n_T\longrightarrow BT^n$ is a projective
morphism. Hence, Deligne's theorem on the degeneration of spectral
sequences yields the following (cf.~Decomposition theorem~\cite{BBD}):
\begin{Th}{\it
Let $M$ be a $T$-equivariant semisimple perverse sheaf of geometric origin
(see~\cite{BBD}) on $Y$. Then, all the spectral sequences (\ref{threeone})
collapse, so that $H^*_T(M)$ is free graded $H^*(BT)$-module and:
$\rk_{H^*(BT)}H^*_T(M)=\dim_{\C}H^*(M)$. Moreover, the cohomology
$H^*_T(M)$ is pure.
}\end{Th}

We will need the following two corollaries of Theorem 8.4.1 (cf.~\cite{Gi2}).
Let $I_0$ be the augmentation ideal in $H^*(BT)$.
\begin{corol}
The restriction morphism: $H^*_T(M)\rightarrow H^*(M)$ induced by the
inclusion $Y\hookrightarrow Y_T$ (as a fibre) vanishes on the submodule
$I_0\cdot H^*_T(M)$ and yields an isomorphism:
$H^*_T(M)/I_0\cdot H^*_T(M)\stackrel{\sim}{\to}H^*(M)$.
\end{corol}

If $M'$ is another perverse sheaf on a $T$-variety $Y'$, satisfying the
conditions of the Theorem, then we have:
\begin{corol}[Kunneth formula]
There is a natural isomorphism:
\[
H^*_T(M\boxtimes M')\cong H^*_T(M)\mathbin{\otimes_{H^*(BT)}}H^*_T(M').
\qquad\Box
\]
\end{corol}
\subsection{}
Assume from now on that $T$ {\it is a torus\/}. Then, the cohomology
$H^*(BT)$ is known to be isomorphic to $\C[\t]$, the polynomial algebra on the
Lie algebra of the torus $T$. We let $I_t$ denote the maximal ideal in 
$\C[\t]$
consisting of all polynomials vanishing at a point $t\in\t$.

Define the family $H_t(M)$, $t\in\t$, of specialized equivariant cohomology
groups of a complex $M$ by
\[
H_t(M)=H^*_T(M)/I_t\cdot H^*_T(M).
\]

For $t\neq0$ the space $H_t(M)$ has no natural grading (for, the ideal $I_t$
is not a homogeneous ideal). Instead, there is a natural increasing
filtration $W_\bullet$ on $H_t(M)$, called the {\it canonical filtration},
which is inherited from the filtration on $H^*_T(M)$ by degree, i.e.:
\begin{equation}
W_i(H_t(M))=\mathop{\fam0 image}(H^0_T(M)\oplus\cdots\oplus H^i_T(M)).
\end{equation}

If $t=0$, then the space $H_0(M)$ has a natural grading and there is a 
canonical
isomorphism: $H^*_0(M)\cong \mathop{\fam0 gr}^W H_{t'}(M)$ for any $t'\neq0$.
If, in addition, the complex $M$ satisfies the assumption of Theorem 8.4.1,
then Corollary 8.4.2 yields an isomorphism: $H^*_0(M)\cong H^*(M)$ ($=$ the
ordinary cohomology of $M$). Thus, for any $t\neq0$, we obtain a natural
isomorphism:
\begin{equation}
\mathop{\fam0 gr}_\bullet\nolimits^W H_t(M)\cong H^*(M).
\end{equation}

\subsection{}
Let $Y^T$ denote the subvariety of $T$-fixed points in $Y$
and $i:Y^T\hookrightarrow Y$ the inclusion. An element $t\in\t$ is
called {\it regular} if the zeroes of the vector field on $Y$ generated
by $t$ coincide with the fixed point subvariety $Y^T$.

A connection between fixed points and the equivariant cohomology is
provided by the following
\renewcommand{\theprop}{\thesubsection}
\begin{Loc}
For any $T$-equivariant complex $M$ on $Y$, the natural push-forward morphism:
\[
i_!:H_t(Y^T,i^!M)\longrightarrow H_t(Y,M)
\]
is an isomorphism.\hfill $\Box$
\end{Loc}
\end{proof}

\end{document}